\documentclass{amsproc}
\usepackage{amsmath}
\usepackage{amscd}
\usepackage{mathrsfs}
\usepackage[all]{xy}
\usepackage{eufrak}
\usepackage{graphicx}
\usepackage[dvips,colorlinks=true,linkcolor=blue]{hyperref}

\makeindex

\newtheorem{theorem}{Theorem}[section]
\newtheorem{proposition}[theorem]{Proposition}
\newtheorem{lemma}[theorem]{Lemma}
\newtheorem{corollary}[theorem]{Corollary}
\theoremstyle{definition}
\newtheorem{definition}[theorem]{Definition}
\newtheorem{example}[theorem]{Example}

\theoremstyle{remark}
\newtheorem{remark}[theorem]{Remark}
\newtheorem{hypothesis}[theorem]{Hypothesis}

\newtheorem{notation}[theorem]{Notation}
\numberwithin{equation}{subsection}

\newcommand{\turnup}[1]{\rotatebox[origin=c]{90}{\ensuremath#1}}

\newcommand{\dem}{\emph{Proof : }}
\newcommand{\CVD}{$\Box$}

\newcommand{\simto}{\xrightarrow[]{\sim}}

\newcommand{\dlog}{\partial_{T,\log}}

\renewcommand{\O}{\mathcal{O}}
\newcommand{\R}{\mathcal{R}}
\renewcommand{\a}{\mathcal{A}}

\newcommand{\B}{\mathrm{B}}

\newcommand{\W}{\mathbf{W}}
\newcommand{\V}{\mathrm{V}}
\newcommand{\F}{\mathrm{F}}
\newcommand{\Fb}{\bar{\mathrm{F}}}
\newcommand{\CW}{\widetilde{\mathbf{CW}}}
\newcommand{\lb}{\bs{\lambda}}
\newcommand{\J}{\mathrm{J}_p}
\newcommand{\T}{\widetilde{T}}
\newcommand{\G}{\mathrm{G}}
\renewcommand{\T}{\mathrm{T}}
\newcommand{\GT}{\bs{\mathcal{G}}}
\newcommand{\Hom}{\mathrm{Hom}}
\newcommand{\Ker}{\mathrm{Ker}}
\newcommand{\Vect}{\underline{\mathrm{Vect}}^{\mathrm{fin}}}
\newcommand{\Rep}{\underline{\mathrm{Rep}}}
\newcommand{\MLCS}{\mathrm{MLS}}
\newcommand{\MLC}{\mathrm{MLC}}
\newcommand{\E}{\mathrm{E}}
\newcommand{\Ed}{\mathcal{E}^{\dag}}

\newcommand{\bs}[1]{\boldsymbol{#1}}
\newcommand{\GE}{\mathrm{G}_{\mathrm{E}}}

\newcommand{\IE}{\mathcal{I}_{\mathrm{E}}}
\newcommand{\PAS}{\mathrm{\mathbf{P}}}
\renewcommand{\d}{\partial_T}
\newcommand{\e}{\mathbf{e}}
\renewcommand{\L}{\mathrm{L}}
\newcommand{\Mt}{\widetilde{\mathrm{M}}}
\newcommand{\M}{\mathrm{M}}
\newcommand{\ph}[1]{\langle #1\rangle}

\begin{document}
\setcounter{tocdepth}{2}

\title[Lubin-Tate Groups and p-adic Differential Equations]{
Rank One Solvable $p$-adic Differential Equations and\\
Finite Abelian Characters via Lubin-Tate Groups.}

\author{Andrea Pulita\\ June 2005}

\address{Equipe de Th\'eorie des Nombres, Universit\'e de Paris 6, 175 rue
 du Chevaleret,
75013, Paris, France}
\thanks{pulita@math.jussieu.fr}
\email{pulita@math.jussieu.fr} 
%\date{July 5, 2004.}

\begin{abstract}
We introduce a new class of exponentials of Artin-Hasse type,
called $\bs{\pi}$-exponentials. These exponentials depend on the
choice of a generator $\bs{\pi}$ of the Tate module of a
Lubin-Tate group $\mathfrak{G}$ over $\mathbb{Z}_p$. They arise
naturally as solutions of solvable differential modules over the
Robba Ring. If $\mathfrak{G}$ is isomorphic to
$\widehat{\mathbb{G}}_m$ over $\mathbb{Z}_p$, we develop methods
to test their over-convergence, and get in this way a stronger
version of the Frobenius Structure Theorem for differential
equations. We define a natural transformation of the
Artin-Schreier  complex into the Kummer complex. This provides an
explicit generator of the Kummer unramified extension of
$\Ed_{K_{\infty}}$, whose residue field is a given Artin-Schreier
extension of $k(\!(t)\!)$, where $k$ is the residue field of $K$.
We then compute explicitly the group, under tensor product, of
isomorphism classes of rank one solvable differential equations.
Moreover, we get a canonical way to compute the rank one
$\varphi$-module over $\Ed_{K_\infty}$ attached to a rank one
representation of
$\mathrm{Gal}(k(\!(t)\!)^{\mathrm{sep}}/k(\!(t)\!))$, defined by
an Artin-Schreier  character.
\end{abstract}

\maketitle

\tableofcontents

\section*{Introduction}
%\addcontentsline{toc}{section}{Introduction}
\subsubsection{}\label{presentation of the problems}
The aim of this paper is to make the theory of rank one solvable
differential equations over the Robba ring $\R_K$ (cf.
\ref{Robba}) as explicit as possible, where $(K,|.|)$ is a
complete ultrametric field with residue field $k$. It is known
(cf. \cite{An}, \cite{Me} and \cite{Ked}) that, under some
restrictions on $K$ and $k$, a solvable $p$-adic differential
module over $\R_K$ becomes unipotent, after pull back on a
covering of $\R_K$, coming from a separable extension of
$\E:=k(\!(t)\!)$. In particular, in \cite{Me} the aim is to
express this module via extension of rank one modules, and get a
$p$-adic analogue of Turritin's classical Theorem for
$K(\!(T)\!)$-differential modules.

Let $\d:=T\frac{d}{dT}$. We shall answer to the following
questions:
\begin{enumerate}
 \item When is a given differential equation
 \begin{equation}
 L=\d-g(T)\;,\quad g(T)=\sum
 a_iT^i \in\R_K
 \end{equation}
 solvable? Can we read the solvability of $L$ from the coefficients $a_i$ of $g(T)$?
 \item What is the irregularity of $L$?
 \item Can we explicitly describe the group (under tensor
 product) $\mathrm{Pic}^{\mathrm{sol}}(\R_K)$ of isomorphism classes of rank one
 solvable differential equations over $\R_K$?
 \item How does this differential equation change under Artin-Schreier  extensions? In
 particular what is the family of rank one solvable modules
 becoming trivial after a given separable extension of $\E=k(\!(t)\!)$?
 \item Let $\mathrm{E}^{\mathrm{sep}}$ be the separable closure of $\E$.
 What is explicitly the rank one $\varphi$-module attached to an
 Artin-Schreier character (or rank one representation) of
 $\GE:=\mathrm{Gal}(\mathrm{E}^{\mathrm{sep}}/\E)$ via the theory of Fontaine-Katz?
 In particular what is the solvable equation attached to this
 $\varphi$-module?
\end{enumerate}

\subsubsection{\textbf{Robba exponentials}}The first example of irregular
solvable differential equation was given by Dwork with the
function
%\begin{equation}
$\exp(\pi T^{-1})$, %\;,
%\end{equation}
which is the Taylor solution at $\infty$ of the irregular operator
$\d +\pi T^{-1}$, where $\pi$ is a solution of the equation
$X^{p-1}=-p$. Dwork  shows that the exponential
$\vartheta(T^{-1}):=\exp(\pi (T^{-p}-T^{-1}))$ is over-convergent
(i.e. converges for $|T|>1-\varepsilon$, for some
$\varepsilon>0$). This provides the so called ``Frobenius
Structure'' isomorphism between $\d+\pi T^{-1}$, and $\d+\pi
T^{-p}$.

\subsubsection{} In \cite{RoIV}, Robba  generalizes the example
of Dwork by producing a class of exponentials, here called
$E_m(T)$, commonly known as Robba's exponentials. Namely Robba
shows that, for all number $\pi_0$ such that
$|\pi_0|=|p|^{\frac{1}{p-1}}$, there exists a sequence
$\alpha_1,\alpha_2,\ldots$ such that, for all $m\geq 1$, the
exponential
\begin{equation}
E_m(T^{-1})=\exp\Bigl(\pi_0\Bigl(\frac{T^{-p^m}}{p^m}+
\alpha_1\frac{T^{-p^{m-1}}}{p^{m-1}}+\cdots+\alpha_mT^{-1}\Bigr)\Bigr)
\end{equation}
converges in the disk $|T|>1$, and hence the operator
$L=\d+\pi_0(T^{-p^m}+
\alpha_1T^{-p^{m-1}}+\cdots+\alpha_mT^{-1})$, with $E_m(T^{-1})$
as solution, is solvable at $\rho=1$. Moreover Robba shows the
necessity of the condition $|\pi_0\alpha_i|=|p|^{\frac{1}{p^i}}$,
for all $i\geq 0$. This construction leads Robba to define the
$p$-adic irregularity of a solvable differential equation as the
slope at $1^-$ of the radius of convergence (cf. \ref{definition
of irregularity}).

But Robba's construction is not sufficient for two reasons. The
first one is that the numbers $\alpha_i$ are obtained as
intersection of a decreasing sequences of disks, and then the
field $K$ must be spherically complete. The second reason is that
Robba was not able to prove the over-convergence of
$E_m(T^{-p})/E_m(T^{-1})$, since the $\alpha_i$'s are essentially
unknown.

\subsubsection{}These problems are solved by S.Matsuda in
\cite{Ma}. He simplifies remarkably the proof of Robba by using
the Artin-Hasse exponential. The idea of introducing the
Artin-Hasse exponential is due to Dwork (cf. \cite[21.1]{Dw}), and
Robba (cf. \cite[10.12]{RoIV}). Matsuda shows that, if $\xi_{m+1}$
is a primitive $p^{m+1}$-th root of $1$, and if
$\xi_{m+1-j}:=\xi_{m+1}^{p^j}$, then we can choose
$\alpha_i=(\xi_i-1)/(\xi_0-1)$. Then
\begin{equation}\label{matsuda's exponentials}
E_m(T^{-1})=\exp\Bigl((\xi_1-1)\frac{T^{-p^m}}{p^m}+
(\xi_2-1)\frac{T^{-p^{m-1}}}{p^{m-1}} + \cdots +
(\xi_{m+1}-1)T^{-1}\Bigr)\;.
\end{equation}
Matsuda proves also that, if $p\neq 2$, then the exponential
$E_m(T^{-p})/E_m(T^{-1})$ is over-convergent. He obtains these
results by a quite complicates, but elementary, explicit
estimation of the valuation of the coefficients of this
exponential.

For the first time we see, in the paper of Matsuda, the algebraic
nature of these analytic exponentials. Indeed, if
$\alpha:\GE\to\Lambda^\times$ is a character of $\GE$ into a
finite extension $\Lambda/\mathbb{Q}_p$, such that $\alpha(\GE)$
is finite, then Matsuda shows that the irregularity of the
differential equation, attached to the $\varphi-\nabla$-module
over $\Ed_K$ defined by $\alpha$, is equal to the Swan conductor
of $\alpha$.
\subsubsection{}Independently from Matsuda, D.Chinellato, under
the direction of Dwork, obtains a new algorithm showing the
existence of the $\alpha_i$s (cf. \cite{Chi}).
\subsubsection{}Even with the great progress given by Matsuda,
Andr\'e, Kedlaya, Crew, Mebkhout, Tsuzuki and others, the
questions given in \ref{presentation of the problems} are still
open, and are the object of this paper.

\subsubsection{} We generalize, and improve, the techniques of Matsuda and Chinellato,
by invoking the Lubin-Tate theory. We recall that the Artin-Hasse
exponential $E(-,T)$ is the group morphism $E(-,T):\W(B)\to
1+TB[[T]]$, functorial on the ring $B$, sending the Witt vector
$\lb=(\lambda_0,\lambda_1,\ldots)\in\W(B)$ into the series
\begin{equation}
E(\lb,T)=\exp\Bigl(\phi_0T+\phi_1\frac{T^p}{p}+\phi_2\frac{T^{p^2}}{p^2}+\cdots\Bigr)\;,
\end{equation}
where $\ph{\phi_0,\phi_1,\ldots}\in B^{\mathbb{N}}$ is the phantom
vector of $\lb$ (cf. \eqref{phantom components}). If
$B=\O_{K^{\mathrm{alg}}}$, then this exponential has bounded
coefficients, and hence converges at least for $|T|<1$. Given a
Frobenius automorphism of $\mathbb{Z}_p[[X]]$, that is a series
$P(X)\in X\mathbb{Z}_p[[X]]$ lifting $X^p\in\mathbb{F}_p[[X]]$, we
consider a sequence $\{\pi_j\}_{j\geq 0}$ in
$\O_{K^{\mathrm{alg}}}$, such that $P(\pi_0)=0$, and
$P(\pi_{j+1})=\pi_j$, for all $j\geq 0$. Then we provide, for all
$m\geq 0$, a Witt vector $[\pi_m]\in\W(\O_{K^{\mathrm{alg}}})$,
whose phantom vector is $\ph{\pi_m,\ldots,\pi_0,0,0,\ldots}$. In
this way, we obtain a large class of exponentials of ``Robba''
type:
\begin{equation}
E_m(T):=E([\pi_m],T)= \exp\Bigl(\pi_mT+\pi_{m-1}\frac{T^p}{p} +
\cdots + \pi_0 \frac{T^{p^m}}{p^m} \Bigr)\;.
\end{equation}
We show then that \emph{the radius of convergence of these
exponentials is $1$ if and only if $P(X)$ is a Lubin-Tate series
(cf. \eqref{Lubin-Tate series}), which thus defines a Lubin-Tate
group $\mathfrak{G}_P$}. In this case $\bs{\pi}:=(\pi_j)_{j\geq
0}$ is a generator of the Tate module of $\mathfrak{G}_P$ (cf.
\ref{E_m converge exactly for T<1 iff P=LT}). If $\mathfrak{G}_P$
is the formal multiplicative group $\widehat{\mathbb{G}}_m$, that
is if $P(X)=(X+1)^p-1$, then we recover Matsuda's exponentials
\eqref{matsuda's exponentials}. On the other hand, if
$P(X)=pX+X^p$, we recover, for $m=0$, Dwork's exponential. Observe
that, in the case considered by Dwork, the formal group
$\mathfrak{G}_P$ is isomorphic, but not equal, to
$\widehat{\mathbb{G}}_m$.

Furthermore, we show that \emph{$E_m(T^p)/E_m(T)$ is
over-convergent, for all $m\geq 0$, if and only if
$\mathfrak{G}_P$ is isomorphic (but not necessary equal) to
$\widehat{\mathbb{G}}_m$.} This is the reason of the
over-convergence of the exponentials $E_m(T^p)/E_m(T)$ of Matsuda
and Dwork.

From this starting point we develop the explicit link between
abelian characters of
$\mathrm{Gal}(k(\!(t)\!)^{\mathrm{sep}}/k(\!(t)\!))$ and rank one
solvable differential equations over $\R_K$, and examine various applications.

\subsubsection{\textbf{Organization of the paper}} In Sections
\ref{notations}, \ref{generalities on rank one eqq}, \ref{Witt
vectors and covectors}, \ref{Notation in AS theory}, and
\ref{Lubin-Tate theory} we give the definitions and recall some
facts used in the sequel.

In Section \ref{Construction of Witt vectors} we define some
canonical Witt vectors with coefficients in $\mathbb{Z}_p[[X]]$,
and show their properties with respect to the Artin-Hasse
exponential. In section \ref{exponentials} we introduce a new
class of exponentials called $\bs{\pi}$-exponentials (cf.
\eqref{definition of et_d(lb,T)}), and show their main properties
with respect to the convergence/over-convergence.

In Section \ref{Dw's theta function}, we give the first
application. Fix a Lubin-Tate group $\mathfrak{G}_P$ isomorphic to
$\widehat{\mathbb{G}}_m$, and a generator $\bs{\pi}=(\pi_j)_{j\geq
0}$ of the Tate module. Let $L$ be a complete discrete valued
field, with residue field $k_L$, and let $\varphi: L\to L$ be a
lifting of the Frobenius $x\mapsto x^p$ of $k_L$. Let
$L_m:=L(\xi_m)$, where $\xi_m$ is a primitive $p^{m+1}$-th root of
$1$. It is well known that we have the Henselian bijection
$$\{\textrm{Finite unramified extensions of } L\}\stackrel{\sim}{\longrightarrow}
\{\textrm{Finite separable extensions of } k_L\}.$$ We shall
describe an inverse of this map. Let $k'/k_L$ be a finite cyclic
abelian extension of degree $d$, and let $L'/L$ be the
corresponding unramified extension. If $(d,p)=1$, and if $k_L$
contains the $d$-th roots of $1$, then $k'/k_L$ is of Kummer type,
and hence $L'=L(\theta)$, where $\theta$ is the Teichm\"uller
representative of a Kummer generator $\bar{\theta}\in k'$.

On the other hand, if $d=p^m$, then $k'$ is of Artin-Schreier type
(cf. \ref{R(nu_0,...,nu_m)}), and it is generated, over $k_L$, by
(the entries of) a Witt vector
$\bar{\bs{\nu}}\in\W_m(k_L^{\mathrm{sep}})$, which is solution of
an equation of the type
$\Fb(\bar{\bs{\nu}})-\bar{\bs{\nu}}=\bar{\lb}$, where
$\bar{\lb}\in\W_m(k_L)$ is a so called Witt vector ``defining''
$k'$. In this case $L_m'/L_m$ is again a Kummer extension, since
all cyclic extensions of $L_m$ whose degree is $p^m$ are Kummer.
Now choose an arbitrary lifting $\lb\in\W_m(\O_L)$ of $\bar{\lb}$,
and solve the equation $\varphi(\bs{\nu})-\bs{\nu}=\lb$,
$\bs{\nu}\in\W_m(\widehat{L}^{\mathrm{unr}})$. Then a Kummer
generator $\theta$ of $L_m'$ is given by the value at $T=1$ of a
certain $\bs{\pi}$-exponential, called $\theta_{p^m}(\bs{\nu},T)$,
(cf. \ref{L' is generate by theta!!}).

The Artin-Schreier  theory and Kummer theory are given by some
\emph{complexes} computing the Galois cohomology. Roughly
speaking, we shall obtain a natural transformation of functors
which ``deforms'' the Artin-Schreier  complex into the Kummer
complex and induces a quasi isomorphism:
\begin{equation}
\xymatrix{ 1 \ar[r]& (L_m)^\times\ar[r]^-{\;\;x\mapsto
x^{p^{m+1}}} &
(L_m)^\times\ar[r]& 1\;\phantom{.}\\
 0 \ar[r]& \W_m(k_L)
\ar@{..>}[u]_{\theta}\ar[r]_{\Fb-1}& \W_m(k_L)
\ar[r]\ar@{..>}[u]_{e^{p^m}} & 0 \;. }
\end{equation}
Actually, such a natural transformation can not exist, because the
Artin-Schreier  complex is in characteristic $p$, and the Kummer
complex is in characteristic $0$. As a matter of fact, we lift the
Artin-Schreier  complex to characteristic $0$ and deform it into
the Kummer complex, by using \emph{the value at $T=1$ of some
over-convergent $\bs{\pi}$-exponentials called $\theta_{p^m}(-,T)$
and $\mathrm{e}_{p^m}(-,T)^{p^{m+1}}$} (see diagram \eqref{lifting
characters diagram for L}). This provides a well defined morphism
between the cohomologies.

Under some assumptions on $K$ (cf. \eqref{hypothesis on K for
AS-Kummer theory}), even if the field $L=\Ed_K$ is not complete,
we show that this diagram exists for $\Ed_K$ and its finite
unramified extensions (cf. \ref{lift-diagr-for Ed}). The
commutative diagram is then:
\begin{equation}\label{big diagram for Edag}
\xymatrix{ 1\ar[r]&\bs{\mu}_{p^{m+1}}\ar[r]&(\Ed_{K_m})^{\times}
\ar[r]^{f\mapsto f^{p^{m+1}}}&
(\Ed_{K_m})^{\times}\ar[r]^-{\delta_{\mathrm{Kum}}}&
\mathrm{H}^1(\G_{\Ed_{K_m}},\bs{\mu}_{p^{m+1}})\ar[r]&1\\
&\W_{m}(\O_{K}^{\sigma=1})\ar@{^{(}->}[r]\ar@{->>}[u]\ar@{->>}[d]&
\W_{m}(\O_{K}^{\dag})\ar[r]_{\varphi-1}\ar[u]_{\theta_{p^{m}}(-,1)}
\ar@{->>}[d]&
\W_{m}(\O_{K}^{\dag})\ar[u]_{\mathrm{e}_{p^{m}}(-,1)^{p^{m+1}}}
\ar@{->>}[d]&&\\
0\ar[r]&\mathbb{Z}/p^{m+1}\mathbb{Z}
\ar[r]\ar@/^{3pc}/@{..>}[uu]^{\wr}& \W_{m}(\E)\ar[r]_-{\Fb-1}&
\W_{m}(\E)\ar[r]_-{\delta}&
\mathrm{H}^1(\GE,\mathbb{Z}/p^{m+1}\mathbb{Z})
\ar@{..>}[uu]_{\overline{\mathrm{e}}:=
\overline{\mathrm{e}_{p^{m}}(-,1)^{p^{m+1}}}}\ar[r]&0}
\end{equation}
where
$\mathrm{G}_{\E}=\mathrm{Gal}(k(\!(t)\!)^{\mathrm{sep}}/k(\!(t)\!))$
and $\mathrm{G}_{\Ed_{K_m}}=
\mathrm{Gal}(\mathcal{E}_{K_m}^{\dag,\mathrm{alg}}/\Ed_{K_m})$. We
specify the kernel and the image of the morphism
$\overline{\mathrm{e}}$ between the cohomologies. If
$\overline{\bs{f}}(t)\in\W_m(k(\!(t)\!))$ is a Witt vector
defining an Artin-Schreier  separable extension of $k(\!(t)\!)$,
then (up to add the $p^{m+1}$-th roots of $1$) a generator of the
corresponding unramified  extension of $\Ed_{K_m}$ is given by
$\theta_{p^m}(\bs{\nu},1)$, where
$\bs{\nu}\in\W_m(\widehat{\mathcal{E}}_K^{\mathrm{unr}})$ is a
solution of the equation $\varphi(\bs{\nu})-\bs{\nu}=\bs{f}(T)$,
and $\bs{f}(T)\in\W_m(\O_K^{\dag})$ is an arbitrary lifting of
$\overline{\bs{f}}(t)$.

\subsubsection{}
Let $K_m:=K(\pi_m)=K(\xi_m)$, $K_\infty:=\cup_mK_m$, and let $k_m$
be the residue field of $K_m$. In Sections \ref{classific-rk-1},
\ref{Classification of rank one differential equats}, and
\ref{Proofs of the stat} we classify all solvable rank one
differential equations over $\R_{K_\infty}$. The key point is the
following equality, arising from the diagram \eqref{big diagram
for Edag}, and useful for describing the Kummer generator
$\theta_{p^m}(\bs{\nu},1)$:
\begin{equation}
\theta_{p^m}(\bs{\nu},1)^{p^{m+1}}=\mathrm{e}_{p^m}(\bs{f}(T),1)^{p^{m+1}}.
\end{equation}
The expression $\mathrm{e}_{p^m}(\bs{f}(T),1)$ has no meaning,
because $\mathrm{e}_{p^m}(-,Z)$ is not over-convergent as a
function of $Z$. We make sense of this symbol in some cases: in
Sections \ref{Classification of rank one differential equats} we
define a class of exponentials of the form
\begin{equation}\label{desired theta}
\mathrm{e}_{p^m}(\bs{f}^-(T),1)=
\exp\left(\pi_m\phi_0^-(T)+\pi_{m-1}\frac{\phi_1^-(T)}{p}+\cdots+
\pi_0\frac{\phi_{m}^-(T)}{p^m}\right)\;,
\end{equation}
where $\bs{f}^-(T)\in\W_m(T^{-1}\O_{K}[T^{-1}])$, and
$\ph{\phi_0^-(T),\cdots,\phi_m^-(T)}\in
(T^{-1}\O_{K}[T^{-1}])^{m+1}$ is its phantom vector. This
exponential is $T^{-1}$-adically convergent and defines a series
in $1+T^{-1}\O_{K_m}[[T^{-1}]]$, whose $p^{m+1}$-th power lies in
$\R_{K_m}$. These Witt vectors correspond to totally ramified
Artin-Schreier  extensions of $\E:=k(\!(t)\!)$. The exponential
\eqref{desired theta} is then the desired Kummer generator
$\theta_{p^m}(\bs{\nu},1)$.

\subsubsection{} We state then the explicit
bijection between the abelian Galois theory for $\E=k(\!(t)\!)$,
and the theory of rank one differential equations over
$\R_{K_\infty}$. Matsuda, in \cite{Ma}, has pointed out, under
some restrictions, that such a correspondence should exist. We go
further by removing any restrictions, improving his methods, and
by making the correspondence more explicit (cf. \ref{Introductive
foundamental Theorem}, \ref{rewritten in terms of P-exp}). Namely
we introduce the fundamental exponential
$\mathrm{e}_{p^m}(\bs{f}^-(T),1)$. We show that every rank one
differential module $M$ over $\R_{K_\infty}$ comes, by scalar
extension, from a module $M_{]0,\infty]}$, over
$K_{\infty}[T^{-1}]$, whose Taylor solution at $\infty$ is  of the
form
\begin{equation}
T^{a_0}\cdot\mathrm{e}_{p^m}(\bs{f}^-(T),1)\;,
\end{equation} for some $m\geq
0$, $a_0\in\mathbb{Z}_p$, and
$\bs{f}^-(T)\in\W_m(T^{-1}\O_{K_m}[T^{-1}])$. Moreover, the
isomorphism class of $M$ depends only on the class of $a_0$ in
$\mathbb{Z}_p/\mathbb{Z}$, and on the Artin-Schreier  character
$\alpha$ defined by the reduction of $\bs{f}^-(T)$ in
$\W_m(k_m(\!(t)\!))$. Suppose that $a_0$ belongs to
$\mathbb{Z}_{(p)}:=\mathbb{Q}\cap\mathbb{Z}_p$. Then $a_0$
corresponds to the moderate extensions of $\E=k(\!(t)\!)$,
generated by $t^{a_0}$. On the other hand, $\bs{f}^{-}(T)$
corresponds to the Artin-Schreier  extension given by (the kernel
of) the Artin-Schreier  character defined by the reduction
$\bs{f}^-(T)$. We recover in this way the well known bijection
\begin{equation}
\left\{\begin{array}{c}
        \textrm{Rank one }\\
        \textrm{characters of
        }\mathcal{I}_{k_\infty(\!(t)\!)}
        \end{array}
 \right\}\stackrel{\sim}{\longrightarrow}
  \left\{\begin{array}{c}
        \textrm{Isomorphism classes of rank one}\\
        \textrm{solvable differential equations}\\
        \textrm{over }\R_{K_\infty}\textrm{with rational residue}
        \end{array}
 \right\}\;,
\end{equation}
where $\mathcal{I}_{k_\infty(\!(t)\!)}$ is the inertia subgroup of
$\mathrm{Gal}(k_\infty(\!(t)\!)^{\mathrm{sep}}/k_\infty(\!(t)\!))$.

\subsubsection{}The central point is that the following $\bs{\pi}$-exponential is
over-convergent
\begin{equation}
\frac{\mathrm{e}_{p^m}(\bs{f}_{(\Fb)}^-(T),1)}{\mathrm{e}_{p^m}(\bs{f}^-(T),1)}=
\mathrm{e}_{p^m}(\bs{f}_{(\Fb)}^-(T)-\bs{f}^-(T),1)\;,
\end{equation}
where $\bs{f}_{(\Fb)}^-(T)$ is an arbitrary lifting of the
reduction $\Fb(\overline{\bs{f}^-}(t))\in\W_m(k_m(\!(t)\!))$. If a
lifting of the $p$-th power map $\varphi:\R_K\to\R_K$ is given,
then this result implies the usual Frobenius Structure Theorem.
Observe that we do not need the existence of $\varphi$ (cf.
\ref{we non need the existece of phi}), because actually the
isomorphism class of a given module $M$ depends only on the
reduction of $\bs{f}^-(T)$ in characteristic $p$. This represents
a progress in two directions, with respect to the analogous
Theorem of \cite{Ch-Ch}: firstly we do not suppose $k$ perfect,
and secondly we get a precise description of the isomorphism class
of $M$. In particular, we find that, if $a_0=0$, then ``the
order'' (cf. \ref{frob str}) of the Frobenius structure is $1$
(cf. \ref{final remarks}).

In Section \ref{Proofs of the stat} we prove these Theorems
essentially by reducing the study to the ``elementary''
$\bs{\pi}$-exponentials corresponding to simpler Witt vector
called $s$-co-monomials. These exponentials are studied in detail
in sections \ref{exponentials}.

\subsubsection{} We give then some complements (Sections
\ref{DESC-CAR-GR},\ref{iiii},\ref{OVER-CAL(E)},\ref{COMP-IRR},\ref{tannakian
group section}). In particular, in the Section \ref{iiii}, we
compute the group of rank one solvable equations killed by a given
Artin-Schreier  extension and then answer the question $(4)$ of
\ref{presentation of the problems}. In Section \ref{OVER-CAL(E)}
we extend the definition of $\bs{\pi}$-exponentials to a larger
class of differential equations, and
%give an
%\emph{incomplete description} of the equations over
%$\mathcal{E}_K$, which shall correspond to \emph{continuous} (but
%not necessary finite) characters of the inertia
%$\mathcal{I}_{k(\!(t)\!)}$. In these sections
we provide an algorithm (see proof of \ref{criteria of solvability
lemma}), which gives a \emph{criterion of solvability}, (cf.
\ref{crit fo solv}). In particular we show that \emph{there is no
irregular rank one equations if $K/\mathbb{Q}_p$ is unramified}
(cf. \ref{no equation over abs unr}). This answers the question
$(1)$ of \ref{presentation of the problems}. Then we \emph{compute
the irregularity} in some classical cases (cf. \ref{COMP-IRR}). We
describe the Tannakian group of the category whose objects are
successive extensions of rank one solvable modules. We remove the
hypothesis ``$K$ is spherically complete'' present in the
literature. In section \ref{PHI-MOD} we compute the
$\phi-\nabla$-module attached to a character with finite image of
$\GE$. This answers question $(5)$ of \ref{presentation of the
problems}.

\specialsection*{Acknowledgments} %\textsc{Acknowledgments :}
The author wish to express his gratitude to professor Gilles
Christol for his guidance and constant encouragement. A particular
thanks goes to Prof. B.Chiarellotto for his support and for
general remarks. We also thank Prof. Y.Andr\'e, P.Colmez,
J.Oesterle, N.Tsuzuki and friends O.Brinon, D.Chinellato and
A.Marmora for useful occasional discussions. We thank the referee
for the help he has given to us in the improvement of the
redaction.

%
%
%
%
%
%
%
%
%
%
%
%
%
%
%
%
%
%     F I R S T   P A R T   -------  F I R S T    P A R T  ----- F I R S T    P A R T
%
%
%
%
%
%
%
%
%
%
%
%
%
%
%
%
%
%
%
%
%
%
%
%
%
\specialsection{\textbf{Definitions and Notations}}
\subsection{General notations}\label{notations} Let $p>0$ be a fixed prime number. Let
$(K,|.|)$ be a complete valued field containing
$(\mathbb{Q}_p,|.|)$. For every valued extension field $L/K$, we
denote by $\mathcal{O}_{L}=\{x\in L\;|\;|x|\leq 1\}$ the ring of
integers of $L$, by $\mathfrak{p}_{L}=\{x\in L\;|\;|x|<1\}$
\index{pl@$\mathfrak{p}_L$}its maximal ideal and by
$k_L=\mathcal{O}_{L}/\mathfrak{p}_{L}$ \index{kl@$k$, $k_L$,
$k^{\mathrm{alg}}$, $K^{\mathrm{alg}}$}its residue field. We set
$k:=k_K$, and take $K^{\mathrm{alg}}$ to be a fixed algebraic
closure of $K$, and $k^{\mathrm{alg}}=k_{K^{\mathrm{alg}}}$ will
be its residue field. $\Omega/K$ \index{Omega@$\Omega$}will be a
spherically complete extension field containing
$K^{\mathrm{alg}}$, satisfying $|\Omega|=\mathbb{R}_{\geq 0}$, and
whose residue field $k_{\Omega}/k$ is not algebraic. We set
\index{omega@$\omega$}
\begin{equation*}
\omega:=|p|^{\frac{1}{p-1}}.
\end{equation*}
We denote by
$\d:=T\frac{d}{dT}$\index{derivation@$\partial_T:=T\frac{d}{dT}$}
the usual derivation. For all rings $R$ we denote by $R^{\times}$
the group of invertible elements in $R$.
\subsubsection{\textbf{Analytic functions and the Robba ring.}}
For every (non vacuous) interval $I\subseteq
[0,\infty[\subset\mathbb{R}$ we set%denote by $\a_K(I)$
\index{AKI@$\a_K(I)=$ analytic functions} %(or simply
%$\a(I)$) the ring of convergent analytic functions on $I$, that is
%the ring formed by the series $f(T)=\sum_{i\in\mathbb{Z}}a_i T^i$,
%$a_i\in K$, such that $\lim_{i\to\pm\infty}|a_i|\rho^i=0$, for all
%$\rho\in I$.
$\a_K(I):=\{\sum_{i\in\mathbb{Z}}
a_iT^i\;|\;\sup_i(|a_i|\rho^i)<\infty\;,\forall\rho\in I\}$, The
topology of $\a_K(I)$ is defined by the family of absolute values
\begin{equation}
|f(T)|_\rho:=\max_{i\in\mathbb{Z}}|a_i|\rho^i\;,
\qquad\forall\;\rho\in I.
\end{equation}
Let $\R_K:=\cup_\varepsilon\a_K(]1-\varepsilon,1[)$
\index{RK@$\R_K=$ Robba ring} be the Robba ring. $\R_K$ is
complete with respect to the limit topology. Let
$\mathcal{E}_K:=\{\sum_{i\in\mathbb{Z}}a_iT^i\;|\;
\sup_i|a_i|<\infty\;,\lim_{i\to-\infty}a_i=0 \}$
\index{EK@$\mathcal{E}_K$ Amice ring} be the Amice ring.
$\mathcal{E}_K$ is endowed with the Gauss norm $|.|_1$ and is
complete. We denote by
$\O_{\mathcal{E}_K}:=\{f\in\mathcal{E}\;|\;|f|_1\leq 1\}$ its ring
of integers. If the valuation on $K$ is not discrete we may have
$|a_i| < \sup_i|a_i|$, for all $i\in\mathbb{Z}$.
%
%The Robba ring $\R_K$ \index{RK@$\R_K=$ Robba ring}(or simply
%$\R$) is the ring of germs of convergent analytic functions at
%$1^-$, that is $\R_K=\cup_\varepsilon\a_K(]1-\varepsilon,1[)$,
%with the limit topology. Let $\mathcal{E}_K$
%\index{EK@$\mathcal{E}_K$ Amice ring}be the Amice ring. The
%elements of $\mathcal{E}_K$ are series
%%$f(T):=\sum_{i\in\mathbb{Z}}a_iT^i$, such that $\lim_{i\to
%-\infty}|a_i|=0$, and $\sup_i|a_i|<+\infty$. If the valuation on
%$K$ is not discrete we may have $|a_i| < \sup_i|a_i|$, for all
%$i\in\mathbb{Z}$. $\mathcal{E}_K$ is complete with respect to the
%topology given by the absolute value $|f(T)|_1:=\sup_i|a_i|$. We
%denote by $\O_{\mathcal{E}_K}:=\{f\in\mathcal{E}\;|\;|f|_1\leq
%1\}$ its ring of integers.
\begin{definition}\label{Robba}
For all algebraic extensions $H/K$ \index{H@$H=$ algebraic
extension of $K$}we set
\begin{equation}\label{a_H=a_K otimes H}
\a_H(I):=\a_K(I)\otimes_K H\;,\quad\R_{H}:=\R_K\otimes_K H\;,\quad
\mathcal{E}_H:=\mathcal{E}_K\otimes_K H\;.
\end{equation}
\end{definition}
Since $K$ is algebraically closed in $\a_K(I)$ (resp. $\R_K$,
$\mathcal{E}_K$) then $\a_H(I)$ (resp. $\R_H$, $\mathcal{E}_H$) is
a domain. All $p$-adic differential equations over $\a_H(I)$
(resp. $\R_H$, $\mathcal{E}_H$) come, by scalar extension, from an
equation over $\a_{L}(I)$ (resp. $\R_L$, $\mathcal{E}_L$) with
$L/K$ finite. This will justify the definition \ref{Picsol}.
\begin{definition}\label{defidefi of f^- + a_0 + f^+}
For all formal series $f(T)=\sum_{i\in\mathbb{Z}}a_iT^i$ we define
\begin{equation}
f^-(T):=\sum_{i\leq -1}a_{i}T^{i} \quad,\qquad f^+(T):=\sum_{i\geq
1}a_{i}T^{i}\;,
\end{equation}
we have $f(T)= f^{-}(T) + a_0 + f^{+}(T)$.
\end{definition}
\begin{definition} \label{Ed_{K,T}}
For all algebraic extension $H/K$, let \index{EdagH@$\Ed_H:=\R_H
\cap \mathcal{E}_H$} $\Ed_{H,T}:=\R_H \cap \mathcal{E}_H$. We
denote by $\O^\dag_{H,T}:=\O_{\mathcal{E}_H}\cap\R_H$. If no
confusion is possible we will write $\Ed_{H}$ (resp.$\O^\dag_{H}$)
instead of $\Ed_{H,T}$ (resp. $\O^\dag_{H,T}$).
\end{definition}

\begin{remark} \label{residue field is alwais E}
The quotients $\O_{\mathcal{E}_K}/\{f\!\in\!
\O_{\mathcal{E}_K}:|f|_1<1\}$ or
$\O_{K}^{\dag}/\{f\!\in\!\O_{K}^{\dag}:|f|_1<1\}$ is $k(\!(t)\!)$
if and only if the valuation on $K$ is discrete. Nevertheless, if
the valuation is not discrete, the rings $\O_{\mathcal{E}_K}$ and
$\O^{\dag}_{K}$ are always local, their maximal ideals
$\mathfrak{p}_{\O_{\mathcal{E}_K}}$ and
$\mathfrak{p}_{\O^{\dag}_{K}}$ are formed by series
$f=\sum_ia_iT^i$ such that $|a_i|<1$, for all $i\in\mathbb{Z}$,
observe that, since the valuation is not discrete, this condition
do not implies that $|f|_1<1$. The residue fields
$\O^{\dag}_{K}/\mathfrak{p}_{\O^{\dag}_{K}}$, and
$\O_{\mathcal{E}_K}/\mathfrak{p}_{\O_{\mathcal{E}_K}}$ are
actually always equals to $k(\!(t)\!)$.
\end{remark}

\subsection{Generalities on rank one differential equations}
\label{generalities on rank one eqq} \label{matrix of derivation}
Let $\mathrm{B}$ be one of the rings $\a_K(I)$, $\R_K$, $\Ed_K$,
$\mathcal{E}_K$. Let $\d-g(T)$, $g(T)\in\mathrm{B} $ be a first
order linear differential operator. The differential module
defined by $\d-g(T)$ is the free rank one module $M$ over
$\mathrm{B}$, endowed with the action of the derivation $\d:M\to
M$ given, in the chosen basis $\e$, by $\d(\e)=g(T)\cdot\e$. We
will say that $g(T)$ is the matrix of the derivation $\d$ in the
basis $\e$. In the sequel we will work with both derivations $\d$
and $d/dT$. We set
\begin{eqnarray}\label{g_[s]-g_s}
g_s(T)=\textrm{the matrix of }\d^s&;& g_{[s]}(T)=\textrm{the
matrix of }\bigl(d/dT\bigr)^{s}\;;
\end{eqnarray}
Then one has
\begin{equation}\label{eq:g_{[s]}}
g_{[s+1]}=\frac{d}{dT}(g_{[s]}(T))+ g_{[s]}(T)g_{[1]}(T)\;,\quad
g_{[0]}(T):=1\;.
\end{equation}
\index{g_s@$g_s=$matrix of $\d^s$} \index{g_{[s]}@$g_{[s]}=$matrix
of $d/dT$}
%\begin{eqnarray}\label{g_[s]-g_s}
%g_s(T)=\textrm{the matrix of }\d^s&;&
%g_{[s]}(T)=\textrm{the matrix of }\bigl(d/dT\bigr)^{s}\;.
%\end{eqnarray}
%\index{g_s@$g_s=$matrix of $\d^s$} \label{g_{[s]}@$g_{[s]}=$matrix
%of $d/dT$}\label{eq:g_{[s]}} We have $g_{s+1}=\d(g_{s}(T))+
%g_{s}(T)g_1(T)$, $g_0(T):=1$, and
%$g_{[s+1]}=\frac{d}{dT}(g_{[s]}(T))+ g_{[s]}(T)g_{[1]}(T)$,
%$g_{[0]}(T):=1$.
Let $\mathrm{C}$ be a $\mathrm{B}$-differential algebra. A
solution of $\d-g(T)$, $g(T)\in\mathrm{B}$, with values in
$\mathrm{C}$ is, by definition, an element $y\in\mathrm{C}$
satisfying $d(y)=g(T)\cdot y$. If $M$ is the rank one module
defined by $d-g(T)$, then the solution $y$ define a morphism of
$\mathrm{B}$-modules $\e\mapsto y:M\to \mathrm{C}$ commuting with
the derivation. \label{change of basis} The operator corresponding
to the basis $f\cdot\e$, $f\in\mathrm{B}^\times$, is
$\d-(g(T)+\frac{\d(f)}{f})$. On the other hand, the tensor product
of the modules defined by $\d-g(T)$ and $\d-\tilde{g}(T)$ is the
module defined by the operator $\d-(g(T)+\tilde{g}(T))$.\label{d +
g + gtilde} \label{tensor product} Then we will identify the
group, under tensor product, of isomorphism classes of (free) rank
one differential modules (here called $\mathrm{Pic}(\mathrm{B})$)
with the group
$$\mathrm{B}/\dlog(\mathrm{B}^{\times})\;,$$ where
$\dlog:\mathrm{B}^\times\to\mathrm{B}$ is the morphism of groups
$f\mapsto\d(f)/f$\index{dlog@$\dlog(f):=\frac{d(f)}{f}$}.
\subsubsection{\textbf{Taylor solution and radius of convergence}}

\label{Taylor solutions and radius of convergence} \emph{ }

Let $I\subseteq\mathbb{R}_{\geq 0}$ be some interval. In this
subsection, $M$ will be a \emph{rank one} $\a_K(I)$-differential
module defined by the operator $\d-g(T)$.

\label{formal solution} Let $x\in\Omega$, $|x|\in I$. We regard
$\Omega[[T-x]]$ as an $\a_K(I)$-differential algebra by the Taylor
map $f(T)\mapsto\sum_{k\geq
0}(\frac{d}{dT})^k(f)(x)\frac{(T-x)^k}{k!}:\a_K(I)\to\Omega[[T-x]]$.
The Taylor solution of $\d-g(T)$ at $x$ is (recall that
$g(T)=Tg_{[1]}(T)$)
\begin{equation}\label{s_x(T)}
s_x(T):=\sum_{k\geq 0} g_{[k]}(x)\frac{(T-x)^k}{k!}\;.
\end{equation}
\index{Taylor solution} Indeed $\d (s_x(T))=g(T)s_x(T)$.
The radius of convergence of $s_x(T)$ at $x$ is, by the usual
definition,
$Ray(M,x)=\liminf_s(|g_{[k]}(x)|/|k!|)^{-\frac{1}{k}}$.\index{Ray(M,x)@$Ray(M,x)$, $Ray(M,\rho)$}
\begin{definition}\label{eq:radius}
The radius of convergence of $M$ at $\rho\in I$ is
\begin{equation*}
Ray(M,\rho):=\min\Bigl(\rho\;,\;\liminf_k(|g_{[k]}|_\rho/|k!|)^{-1/k}\Bigr)=
\min\Bigl(\rho\;,\;\omega\bigl[\limsup_k(|g_{[k]}|_\rho)^{1/k}\bigr]^{-1}\Bigr).
\end{equation*}
\end{definition}
The second equality follows from the fact that the sequence
$|k!|^{1/k}$ is convergent to $\omega$, and $|g_{[k]}|_\rho^{1/k}$
is bounded by $\max(|g_{[1]}|_\rho,\rho^{-1})$. The presence of
$\rho$ in the minimum makes this definition invariant under change
of basis in $M$.
\begin{theorem}[Transfer]\label{transfer}
\index{Transfer Theorem} For all $\rho\in I$ we have
\begin{equation}
Ray(M,\rho)=\min
\bigl(\;\;\rho\;,\;\;\inf_{x\in\Omega,|x|=\rho}Ray(M,x)\;\bigr).
\end{equation}
Assume now that $I=[0,\rho]$. Then
$Ray(M,\rho)=\min(\rho,\min_{x\in\Omega,|x|\leq\rho}Ray(M,x))$. In
particular $Ray(M,\rho)\leq\min(\rho,Ray(M,0))$.
\end{theorem}
\emph{Proof : } Since for $\rho=|x|$ we have
$|g_{[s]}(T)|_\rho\geq |g_{[s]}(x)|$, hence by definition
\ref{eq:radius}, $Ray(M,\rho)\leq \min(\rho,Ray(M,x))$. Let
$t_\rho\in\Omega$ be such that
$\{x\in\Omega\;|\;|x-t_\rho|<\rho\}\cap K=\emptyset$, then
$|g_{[s]}|_\rho = |g_{[s]}(t_\rho)|$, for all $s\geq 0$
(\cite[9.1]{Ch-Ro}), hence $Ray(M,\rho) =
\min(\rho,Ray(M,t_\rho))$. The last assertion follows similarly.
$\Box$
\begin{lemma}[Small Radius]\label{small radius}
\index{Small radius lemma} Let $\rho\in I$. Then
\begin{equation}\label{eq:minoration radius}
Ray(M,\rho)\geq \omega\rho\cdot \min(1,|g(T)|_\rho^{-1}).
\end{equation}
Moreover $Ray(M,\rho)<\omega\rho$ if and only if $|g(T)|_\rho >
1$, and in this case we have
\begin{equation}\label{eq:small radius}
Ray(M,\rho)=\omega\rho\cdot|g(T)|_{\rho}^{-1}.
\end{equation}
\end{lemma}
\emph{Proof : } By induction on \eqref{eq:g_{[s]}},
$|g_{[s]}|_\rho\leq
\max(\rho^{-1},|g_{[1]}|_\rho)^{s}=\rho^{-s}\max(1,|g|_\rho)^s$
(cf. \eqref{g_[s]-g_s}), and equality holds if $|g_{[1]}|_\rho >
\rho^{-1}$. Then apply definition \ref{eq:radius}. $\Box$

\begin{definition}\index{solvable at $\rho$} $M$ is called
\emph{solvable} at $\rho\in I$, if $Ray(M,\rho)=\rho$.
\end{definition}
\begin{theorem}[\cite{Ch-Dw}]
The map $\rho\mapsto Ray(M,\rho):I\to \mathbb{R}_\geq$ is
continuous and locally of the form $r \cdot \rho^{\beta+1}$, for
suitable $r\in\mathbb{R}_{\geq}$, and $\beta\in\mathbb{N}$. More
precisely there exist a partition $I=\cup_{n\in\mathbb{Z}} I_n$,
$\sup I_{n}= \inf I_{n+1}$, and two sequences
$\{r_n\}_{n\in\mathbb{Z}}$, $\{\beta_n\}_{n \in \mathbb{Z}}$, such
that $\beta_n\in\mathbb{Z}$, $Ray(M,\rho)=r_n\rho^{(\beta_n+1)}$,
$\forall\;\rho\in I_n$, and (cf. \ref{graphic})
\begin{equation}\label{log-concavity}
\index{logconcavity@$\log$-concavity} \beta_{n}\geq \beta_{n+1}.
\end{equation}
\end{theorem}
\emph{Proof : } The existence of the partition follows from the
Small Radius Lemma \ref{small radius} and Theorem \ref{radius of
frobenius}. For more details see \cite[8.6]{Astx} and
\cite[2.5]{Ch-Dw}. $\Box$
\begin{definition}
We will call the property \eqref{log-concavity} the
\emph{$\log$-concavity} of the function $\rho\mapsto Ray(M,\rho)$.
We will call $\beta_n$ the slope of $M$ in the interior of $I_n$.
More generally if $\rho = \sup I_n=\inf I_{n+1}$, we set
$\index{sl-(M,rho)@$\mathrm{sl}^-(M,\rho)$,
$\mathrm{sl}^+(M,\rho)$, $\mathrm{sl}_F(M)$}
\mathrm{sl}^-(M,\rho):=\beta_n$ and
$\mathrm{sl}^+(M,\rho):=\beta_{n+1}$.
%\begin{equation}
%\index{sl-(M,rho)@$\mathrm{sl}^-(M,\rho)$,
%$\mathrm{sl}^+(M,\rho)$, $\mathrm{sl}_F(M)$}
%\mathrm{sl}^-(M,\rho):=\beta_n\quad,\qquad
%\mathrm{sl}^+(M,\rho):=\beta_{n+1}.
%\end{equation}
\end{definition}
\begin{remark}The Taylor solution of $M\otimes N$ is the product of the
Taylor solutions of $M$ and $N$. Hence, by
\ref{transfer}\label{tensor radius}, $Ray(M\otimes
N,\rho)\geq\min(Ray(M,\rho),Ray(N,\rho))$. \label{min tensor with
radius} If $Ray(M,\rho)\neq Ray(N,\rho)$, then we have
$Ray(M\otimes N,\rho)=\min(Ray(M,\rho),Ray(N,\rho))$.
\end{remark}

\subsubsection{\textbf{Solvability, Slopes and
Irregularities}}\label{solvability} In this subsection, $M$ is the
rank one module over $\R_K$, defined by $\d+g(T)$,
$g(T):=\sum_{i\in\mathbb{Z}}a_iT^i\in\R_K$.

\begin{lemma} \label{algebricity}
There exists $d > 0$ such that $M$ is isomorphic to the module
defined by $\d + \sum_{i\geq -d} a_i T^i$. In other words there
exists $f(T)\in\R_K^\times$ such that
$\partial_{T,\log}(f)=\sum_{i < -d}a_iT^i$.
\end{lemma}
\emph{Proof : } By hypothesis $g(T)\in\a_K(]1-\varepsilon,1[)$,
for some $\varepsilon>0$. Then $\sum_{i\neq
0}a_iT^i/i\in\a_K(]1-\varepsilon,1[)$. In particular $\lim_{i\to
-\infty}|a_i/i|\rho^i=0$, for all
$\rho\in]1-\varepsilon,+\infty[$. Let $d > 0$ be such that
$\sup_{i < -d}(|a_i/i|\rho^i)<\omega$ for a fixed
$\tilde{\rho}\in]1-\varepsilon,1[$. Then $\sup_{i <
-d}(|a_i/i|\rho^i)<\omega$, for all $\rho\geq \tilde{\rho}$. Then
$f(T)=\exp(-\sum_{i < -d} a_i T^{i}/i)$ lies in
$\mathcal{R}_K$.$\Box$

\begin{definition}\label{MLCS}\index{solvable over $\R$}
Let $M$ be a differential module over $\mathcal{R}_K$. The module
$M$ is called \emph{solvable} if and only if $\lim_{\rho\to
1^-}Ray(M,\rho)=1$. We will denote the category of solvable
differential modules over $\R_K$ by
$\MLCS(\R_K)$.\index{MLS(R)@$\MLCS(\R_{K})$}
\end{definition}

\begin{lemma}\label{module over ]0,1[}
\label{ray in 0} \label{graphic} Let $M\in\MLCS(\R_K)$ be defined
in some basis by the operator $\d-g(T)$, $g(T)\in\R_K$. Then
\begin{enumerate}
 \item There exist $0<\varepsilon<1$ and a last slope
 $\beta:=\mathrm{sl}^-(M,1)\geq 0$ such
 that
 \begin{equation}\label{exists last slope}
 Ray(M,\rho)=\rho^{\beta + 1}\;,\quad\textrm{ for all }\rho\in
 ]1-\varepsilon,1[\;.
 \end{equation}
 \item There exists $\varepsilon'$ such that $|g(T)|_\rho \leq 1$,
 for all $\rho\in]1-\varepsilon',1[$.
 \item If $g(T)=\sum_{-d}^{\infty}a_iT^i$, $d>0$,
 then $|a_{-d}|\leq \omega$ and, for $\rho$ close to $0$,
\begin{equation}
Ray(M,\rho)=\omega |a_{-d}|^{-1}\rho^{d+1}.
\end{equation}
 \item Moreover, if $d>0$, and if $|a_{-d}|=\omega$, then
 $\beta=d$.
 \item If $d\leq 0$, then $Ray(M,\rho)=\rho$, for all
 $\rho\in]0,1[$ and $\beta=0$.
\end{enumerate}
\end{lemma}
\emph{Proof : } The slopes are positive natural numbers, hence the
decreasing sequence $\{\beta_n\}_n$ becomes constant for
$n\to\infty$. Then $\beta=\min_{n\in\mathbb{Z}}\{\beta_n\}$. The
second assertion follows from the small radius lemma \ref{small
radius}. Let now $g(T)=\sum_{i\geq -d}a_iT^i$, with $d>0$. We
study the function $\rho\mapsto Ray(M,\rho)/\rho$. Let
\begin{equation}\label{R(M,r), r}
\index{R(M,r)@$R(M,r)=\log(Ray(M,\rho)/\rho)$}
R(M,r):=\log(Ray(M,\rho))-\log(\rho)\quad,\qquad r:=\log(\rho)\;.
\end{equation}
Then $R(M,r)\leq 0$, for all $r\leq 1$, and the function $r\mapsto
R(M,r):\;]-\infty,1[\longrightarrow ]-\infty,0]$ is of the
following form
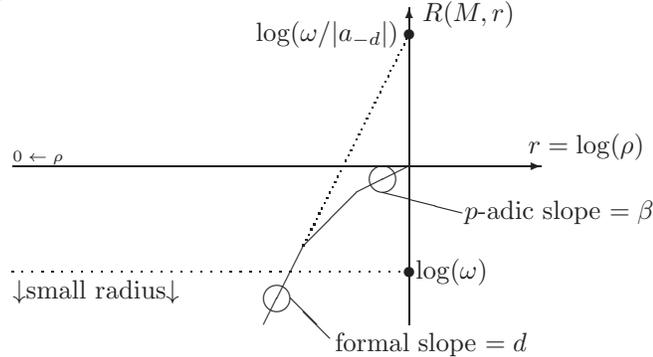
\begin{center}
\begin{picture}(200,130)
% AXES
\put(150,10){\vector(0,1){120}} \put(0,70){\vector(1,0){200}}
% COORDINATE
\put(195,75){$r=\log(\rho)$} \put(155,125){$R(M,r)$}
\put(0,72){\begin{tiny}$0\leftarrow\rho$\end{tiny}}
%\put(152,72){\begin{tiny}$\rho=1$\end{tiny}}
% POLYGON
\put(95,10){\line(1,2){15}} % first line
\put(110,40){\line(1,1){20}} % second line
\put(130,60){\line(2,1){20}} % third line
% COMMENTS
\put(140,65){\circle{10}}     % upper circle
\put(140,60){\line(4,-1){30}}
\put(100,20){\circle{10}}     % lower circle
\put(105,20){\line(1,-1){15}}
% COMMENTS
\put(122,0){formal slope $=d$} \put(171,50){$p$-adic slope
$=\beta$} \put(0,20){$\downarrow$small radius$\downarrow$}
% COMMENTS
\put(147.5,117.5){$\bullet$} % first interesting point
\put(92,117.5){$\log(\omega/|a_{-d}|)$}
\put(147.5,27.5){$\bullet$} % second interesting point
\put(152.5,27.5){$\log(\omega)$}
\qbezier[40](110,40)(130,80)(150,120)
\qbezier[40](0,30)(75,30)(150,30)
\end{picture}
\end{center}
Since $d>0$, one has $|g(T)|_\rho=|a_{-d}|\rho^{-d}>1$, for $\rho$
close to $0$. Hence, near $0$, we can apply the Small Radius Lemma
(cf. \ref{formal and p-adic}): we have $Ray(M,\rho)=\omega
|a_{-d}|^{-1}\rho^{d+1}$. Since $\lim_{\rho\to 1^-}Ray(M,\rho) =
1$, hence by $\log$-concavity and continuity, we must have $\omega
|a_{-d}|^{-1}\geq 1$ (or equivalently $\log(\omega/|a_{-d}|)\geq
0$ as in the picture) and if $|a_{-d}|=\omega$, then, again by
continuity and log-concavity, this graph is a line, and $\beta =
d$. If $d\leq 0$, then $|g(T)|_\rho\leq 1$, for all $\rho< 1$,
hence the Small Radius Lemma gives $R(M,r)\geq \log(\omega)$ for
all $r\leq 0$. Since $R(M,r)\!\to 0$ for $r\!\to 0$ (solvability),
then by $\log$-concavity and continuity this implies $R(M,r)=0$,
$\forall\; r\leq 0$. $\Box$

\begin{remark} \label{formal and p-adic}
\index{Irr(M)@$\mathrm{Irr}(M)=\mathrm{sl}^-(M,1)$,
$\mathrm{Irr}_F(M)$} We maintain the notation of Lemma \ref{module
over ]0,1[} part $(3)$. We recall that
$\mathrm{sl}^+(M,0)=\min(0,d)$ is equal to the classical formal
slope $\mathrm{Irr}_F(M)$ of $M$ as $K(\!(T)\!)$-differential
module (cf. \cite{VS}). This is actually true in all ranks.
\end{remark}

\begin{definition} \label{definition of irregularity}
Let $M$ be a solvable rank one differential module over $\R$. The
$p$-adic \emph{irregularity} of $M$ is the natural number
$\mathrm{Irr}(M):= \mathrm{sl}^-(M,1)$.
\end{definition}

\begin{remark}
If $M$ is defined by an operator $\d + g(T)$,
$g(T)=\sum_{-d}^\infty a_iT^i \in \R_K$, then by log-concavity and
continuity we have $\mathrm{Irr}_F(M)\geq \mathrm{Irr}(M)$.
\end{remark}

\begin{definition}\label{Picsol}
If $K'/K$ is a finite extension, then we denote by
$\mathrm{Pic}^{\mathrm{sol}}(\R_{K'})$ the group, under tensor
product, of isomorphism classes of solvable rank one differential
modules over $\R_{K'}$. For all algebraic extensions $H/K$, we set
\begin{equation}\index{Picsol@$\mathrm{Pic}^{\mathrm{sol}}(\R_H)$}
\mathrm{Pic}^{\mathrm{sol}}(\R_{H}):=\bigcup_{K\subset K'\subset
H\;,\;K'/K\textrm{ finite}}\mathrm{Pic}^{\mathrm{sol}}(\R_{K'}).
\end{equation}
\end{definition}

\begin{corollary}
\label{min tensor radius} \label{tensor product with different
radius} We have $\mathrm{Irr}(M\otimes
N)\leq\max(\mathrm{Irr}(M),\mathrm{Irr}(N))$, for all
$M,N\in\MLCS(\R_K)$. Moreover the equality holds if
$\mathrm{Irr}(M) \neq \mathrm{Irr}(N)$. \CVD
\end{corollary}
\begin{proposition}\label{division of the problem}
Let $\d-g(T)$, $g(T)=\sum_{i\in\mathbb{Z}}a_i T^i\in\R_K$, be a
solvable differential equation. Then $\d-g^-(T)$, $\d-a_0$,
$\d-g^+(T)$ are all solvable (cf. \eqref{defidefi of f^- + a_0 +
f^+}).
\end{proposition}
\dem Let us call $M_{]1-\varepsilon,\infty]}$, $M_0$, $M_{[0,1[}$
the differential modules defined by $\d-g^-(T)$, $\d-a_0$,
$\d-g^+(T)$ respectively. Then
$M=M_{]1-\varepsilon,\infty]}\otimes M_0\otimes M_{[0,1[}$. By the
Small Radius Lemma \ref{small radius}, the equation $\d-g^-(T)$
(resp. $\d-g^+(T)$) has a convergent solution at $\infty$ (resp.
at $0$), hence $Ray(M_{]1-\varepsilon,\infty]},\rho)=\rho$, for
large values of $\rho$ and $Ray(M_{[0,1[},\rho)=\rho$, for $\rho$
close to $0$. While $Ray(M_0,\rho)= R_0\cdot\rho$, for all $\rho$
(cf. \ref{residue}). Hence the slopes of
$M_{]1-\varepsilon,\infty]}$ (resp. $M_{[0,1[}$, $M_0$) in the
interval $]1-\varepsilon,1[$ are strictly positive (resp. strictly
negative, resp. equal to $0$) as in the picture (cf.
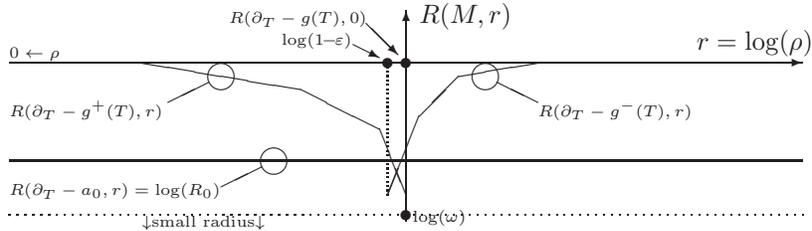
\eqref{R(M,r), r}).
\begin{center}
\begin{picture}(300,80)
% AXES
\put(150,0){\vector(0,1){80}} \put(0,60){\vector(1,0){300}}
% COORDINATE
\put(260,65){$r=\log(\rho)$} \put(155,75){$R(M,r)$}
\put(0,62){\begin{tiny}$0\leftarrow\rho$\end{tiny}}
%\put(152,72){\begin{tiny}$\rho=1$\end{tiny}}
% POLYGON of \partial+g^+(T)
\put(50,60){\line(6,-1){60}} % first line
\put(110,50){\line(2,-1){30}} % second line
\put(140,35){\line(2,-5){10}} % third line
% POLYGON of \partial+g^-(T)
\put(170,55){\line(6,1){30}} % first line
\put(170,55){\line(-1,-1){15}} % second line
\put(155,40){\line(-2,-5){12}} % third line
% POLYGON of \partial + a_0
\put(0,23){\line(1,0){300}} % first line
% INTERESTING POINTS
\put(147.5,57.5){$\bullet$} % first interesting point
\put(83,75){\begin{tiny}$R(\d-g(T),0)$\end{tiny}}
 \put(135,75){\vector(1,-1){12}}

\put(140.5,57.5){$\bullet$} % second interesting point
\put(103,67){\begin{tiny}$\log(1\!\!-\!\!\varepsilon)$\end{tiny}}
\put(130,67){\vector(2,-1){10}}
\qbezier[30](143,10)(143,35)(143,60)
% COMMENTS
 \put(80,55){\circle{10}}     %  circle of R+
 \put(60,45){\line(2,1){15.5}}
 \put(0,40){\begin{tiny}$R(\d-g^+(T),r)$\end{tiny}}
 \put(180,55){\circle{10}}     %  circle of R-
 \put(200,45){\line(-2,1){15.5}}
 \put(200,40){\begin{tiny}$R(\d-g^-(T),r)$\end{tiny}}
 \put(100,23){\circle{10}}     %  circle R0
 \put(80,13){\line(2,1){15.5}}

 \put(0,10){\begin{tiny}$R(\d-a_0,r)=\log(R_0)$\end{tiny}}
\put(147.5,0){$\bullet$} % second interesting point
\put(152.5,0){\begin{tiny}$\log(\omega)$\end{tiny}}
\qbezier[100](0,2.5)(150,2.5)(300,2.5)
\put(50,-2){\begin{tiny}$\downarrow$small
radius$\downarrow$\end{tiny}}
\end{picture}
\end{center}
By \ref{tensor radius}, we have
$Ray(M,\rho)=\inf(Ray(M_{]1-\varepsilon,\infty]},\rho),Ray(M_{[0,1[},\rho),Ray(M_{0},\rho))$,
 for all
$1-\varepsilon<\rho<1$, with the exception of a finite numbers of
$\rho$. By continuity of the radius, we have equality even for
these isolated values of $\rho$. Since $\lim_{\rho\to
1^-}Ray(M,\rho)=1$, this implies $Ray(M_{[0,1[},\rho)=\rho$ for
all $\rho<1$, $Ray(M_0,\rho)=\rho$ for all $\rho$, and
$Ray(M_{]1-\varepsilon,\infty]},\rho)=\rho$ for all $\rho\geq 1$.
\CVD

The classification of the equations of the type $\d-a_0$, $a_0\in
K$, is well known (see \ref{moderate characters}), while the
solvable equations of the form $\d-g^+(T)$ are always trivial:

\begin{proposition}\label{positive queue}
Let $\d-g^+(T)$, $g^+(T)=\sum_{i\geq
1}a_iT^i\in\a_H([0,1[)\subset\R_K$ be solvable at $1^-$ (cf.
\ref{MLCS}). Let $M$ be the module attached to $\d-g^+(T)$, then
\begin{enumerate}
    \item We have $g^+(T)\in T\O_H[[T]]$. Hence $M$ comes, by scalar
    extension, from a differential module $M_{[0,1[}$
    over $\O_H[[T]]$;
    \item We have $Ray(M_{[0,1[}\;,\;\rho)=\rho$, for all $\rho<1$;
    \item $M_{[0,1[}$ is trivial as $\O_H[[T]]$-module;
    \item The exponential $\exp(\sum_{i\geq 1}a_iT^i/i)$ lies in
    $1+T\O_H[[T]]$.
 \end{enumerate}
\end{proposition}
\emph{Proof : } We have $|a_i|\leq 1$, because the Small Radius
Lemma \ref{small radius}. Since $\d-g^+(T)$ has a convergent
solution at $0$ (namely this Taylor solution is $\exp(\sum_{i\geq
1}a_iT^i/i)$), then $Ray(M_{[0,1[},\rho)=\rho$ for all $\rho$
close to $0$. Since $\lim_{\rho\!\to\!1^-}Ray(M_{[0,1[},\rho)=1$,
then by log-concavity we must have $Ray(M_{[0,1[},\rho)=\rho$, for
all $\rho<1$. By the transfer Theorem \ref{transfer} the Taylor
solution $\exp(\sum_{i\geq 1}a_iT^i/i)$ converges in the disk
${|T|<1}$ and belongs to $\O_K[[T]]$ (because a non trivial
solution of a differential equation has no zeros in its disk of
convergence).\CVD
\begin{corollary}\label{Reduction to a K[T^-1]-lattice}
Every rank one solvable differential module over $\R_K$ has a
basis in which the matrix lies in $\O_K[T^{-1}]$.
\end{corollary}
\dem By \ref{positive queue} there exists a basis in which the
matrix lies in $\R_K\cap\O_K[[T^{-1}]]$. The base change matrix to
obtain this basis is an exponential convergent in $[0,1[$. Now, by
\ref{algebricity} we recover the good basis. This last base change
matrix is again an exponential convergent in
$]1-\varepsilon,\infty[$.\CVD

\subsubsection{\textbf{Frobenius structure and $p$-th ramification}}
\begin{definition} \label{what is an absolute Frobenius}
An \emph{absolute Frobenius on $K$}\index{absolute Frobenius} is a
$\mathbb{Q}_p$-endomorphism $\sigma:K \to K$\index{sigma@$\sigma$}
such that $|\sigma(x)-x^p|<1$, for all $x\in \O_{K}$. \\ If an
absolute Frobenius $\sigma:K\to K$ is given, an \emph{absolute
Frobenius on $\R_K$} is then a continuous endomorphism of rings
$\varphi: \R_K \to \R_K$ extending $\sigma$, and such that
\begin{equation}\label{a_i(varphi)}
\varphi(T)-T^p=\sum a_i(\varphi)T^i\;, \textrm{ with }
|a_i(\varphi)|<1 \textrm{ for all }
i\in\mathbb{Z}\;,\;a_i(\varphi)\in K\;.
\end{equation}
\end{definition}
By continuity $\varphi$ is given by $\sigma$ and by the choice of
$\varphi(T)$. Namely if we set $(\sum a_iT^i)^{\sigma}:=\sum
\sigma(a_i)T^i$, then $\varphi(f(T))=f^\sigma(\varphi(T))$, for
all $f\in\R_K$. The simplest absolute Frobenius is given by the
choice $\varphi(T)=T^p$ and we denote it by
$\varphi_\sigma$\label{varphi_sigma^*}\index{phi^*@$\varphi_p$,
$\varphi_\sigma$}.

Let $\varphi:\R\to\R$ be an absolute Frobenius. By scalar
extension (and change of derivation), we have a functor:
$\varphi^*:\MLC(\R)\leadsto\MLC(\R)$. If $M\in\MLC(\R)$ is defined
by the operator $\d-g(T)$, then $\varphi^*(M)$ is defined by the
operator $\d-(\dlog(\varphi(T))\cdot g^\sigma(\varphi(T)))$. The
isomorphism class of $\varphi^*(M)$ does not depend on the choice
of $\varphi$ (cf. \cite[7.1]{Astx}). \label{independence on
varphi}\label{Independence on varphi}

\subsubsection{\textbf{$p$-th ramification}}\label{varphi^*} Let $\sigma$ be an absolute
Frobenius on $K$. For all analytic functions
$f(T):=\sum_ia_iT^i\in\a(I)$, we set \label{varphi_p^*}
$\varphi_p(f(T)):=f(T^p)$. Observe that $\varphi_p$ is not an
absolute Frobenius. We set
$\varphi_{\sigma}(f(T)):=f^\sigma(T^p)$. The $p$-th ramification
map $\varphi_p:\a(I^p)\to\a(I)$ defines, as before, a functor
denoted by $\varphi_p^*: \MLC(\a_K(I^p))\leadsto\MLC(\a_K(I))$.
\begin{theorem}[\cite{Astx}] \label{radius of frobenius}
Let $M\in\MLC(\a_K(I^p))$. Then for all $\rho\in I$
\begin{equation*}
Ray(\varphi_\sigma^*(M),\rho)=Ray(\varphi_p^*(M),\rho)\geq
\rho\min\Bigl(\Bigl(\frac{Ray(M,\rho^p)}{\rho^p}\Bigr)^{1/p}\;,\;
|p|^{-1} \frac{Ray(M,\rho^p)}{\rho^p}\Bigr),
\end{equation*}
and equality holds if $Ray(M,\rho)\neq\omega^p\rho$.
\end{theorem}
\emph{Proof : } Since $f(T)\mapsto f^\sigma(T)$ is an isometry, we
have the first equality. The second one follows from a quite
complex, but elementary, explicit computation.$\Box$
\begin{example}
The radius of the operator $\d- 1/p$ is equal to
$\omega|p|\rho=\omega^p\rho$ (cf. \ref{small radius}), but its
image by Frobenius is the trivial module.
\end{example}
\begin{corollary}
Let $M\in\MLCS(\R_K)$, let $\varphi:\R_K\to\R_K$ be an absolute
Frobenius, then
$\mathrm{sl}^-(\varphi^*(M),1)=\mathrm{sl}^-(M,1)$. \CVD
\end{corollary}
\begin{definition}[Frobenius structure]
\label{frob str} \index{Frobenius structure} Let $M$ be a module
over $\R_K$. We will say that $M$ has a \emph{Frobenius structure
of order $h$}, if $M$ is isomorphic to $(\varphi^*)^h(M)$.
\end{definition}
\begin{remark} \label{frobenius implies solvability}
If $M$ has a Frobenius structure, then it is solvable by Theorem
\ref{radius of frobenius} applied to ``antecedents'' of $M$. (see
\cite[8.6 and 7.7 infra]{Astx}).
\end{remark}

\begin{remark}\label{irreg of frob}
By equation \ref{radius of frobenius} we have
$\mathrm{Irr}(\varphi^*(M))=\mathrm{Irr}(\varphi_p^*(M))=\mathrm{Irr}(M).$
\end{remark}

\subsubsection{\textbf{Moderate characters}}\label{moderate characters}
\label{moderate module} Let $a_0\in K$. We denote by $\M(a_0,0)$
the module defined by the constant operator $\d-a_{0}$ (cf.
\ref{matrix of derivation}). We will call \emph{moderate}
\index{moderate character} every \emph{solvable} differential
module (over $\R_K$) of the form $\M(a_0,0)$. \label{robba shown
that} By \cite[5.4]{RoIV}, $\M(a_0,0)$ is solvable if and only if
$a_0\in\mathbb{Z}_p$. Moreover the equation
$\partial_{T,\log}(f(T))=a_0$ has a solution $f(T)\in\R_K^\times$
if and only if $a_0\in\mathbb{Z}$, and in this case
$f(T)=T^{a_0}$. This shows that the group under tensor product of
moderate differential modules is isomorphic to
$\mathbb{Z}_p/\mathbb{Z}$. On the other hand it is well known that
an $\M(a_0,0)$ has a Frobenius structure if and only if
$a_0\in\mathbb{Z}_{(p)}$.
\begin{lemma} \label{residue}
Let $\alpha(a_{0}):=
\limsup_s(|a_{0}(a_{0}-1)(a_{0}-2)\cdots(a_{0}-s+1)|^{\frac{1}{s}})$.
Then $Ray(\M(a_0,0),\rho)=\rho\cdot R_0\leq \rho$, for all
$\rho>0$, with $R_0:=\min(1,\omega\cdot\alpha(a_0)^{-1})$.
\end{lemma}
\emph{Proof : } A direct computation gives
$g_{[s]}(T)=\alpha_s(a_{0})T^{-s}$, with
$\alpha_s(a_{0}):=a_{0}(a_{0}-1)\cdots(a_{0}-s+1)$
(cf.\eqref{g_[s]-g_s}). Then apply \ref{eq:radius}.$\Box$

\subsection{Notations on Witt Vectors and covectors}\label{Witt
vectors and covectors} Let $R$ be a ring. Notations concerning the
ring $\W(R)$ \index{W(R)@$\W(R)$, $\W_m(R)$} of Witt vectors will
follow \cite{Bou}, except for the indexation ``$m$'' of the ring
$\W_m(R)$ of Witt vector of finite length. We set
$\W_m(R):=\W(R)/\V^{m+1}\W(R)$ (see \ref{FV=p}). We denote by
\begin{equation}\label{phantom components}
\phi_n:=\phi_n(X_0,\ldots,X_n):=X_0^{p^n}+pX_1^{p^{n-1}}+\cdots+p^nX_n
\end{equation}
the Witt polynomial. Vectors in $R^{\mathbb{N}}$ and in $R^{m+1}$
will be distinguished from Witt vectors by the notation
\index{phi_0@$\ph{\phi_0,\phi_1,\ldots}=$notations for phantom
vectors} $\ph{\phi_0,\phi_1,\ldots}$ instead of
$(\phi_0,\phi_2,\ldots)$. For all Witt vector
$\bs{r}=(r_0,r_1,\ldots)\in\W(R)$, the vector
$\phi(\bs{r})=\ph{\phi_0(r_0),\phi_1(r_0,r_1),\ldots}$ is called
the phantom vector of $\bs{r}$. The map
$\bs{r}\mapsto\phi(\bs{r}):\W(R)\to R^{\mathbb{N}}$ is a ring
morphism.
\begin{lemma}[\protect{\cite[Lemme 3 $\S 1, N^02$]{Bou}}]\label{phi is inj}
Let $\lb\mapsto\phi(\lb):\W(R)\stackrel{\phi}{\to} R^{\mathbb{N}}$
be the phantom component map. If $p\in R$ is not a zero divisor,
then $\phi$ is injective. If $p\in R$ is invertible then $\phi$ is
bijective. $\Box$
\end{lemma}
\begin{lemma}[\protect{\cite[Lemme 2 $\S 1, N^02$]{Bou}}]
\label{congruece on the phantom components} Let $\sigma:R\to R$ be
a ring morphism satisfying $\sigma(a)\equiv a^p \pmod{pR}$, for
all $a\in R$. Then a vector $\ph{\phi_0,\ldots,\phi_{m}}\in
R^{m+1}$ is the phantom vector of a Witt vector if and only if
\begin{equation}
\phi_{i}\equiv \sigma(\phi_{i-1})\mod p^{i} R\;, \quad\textrm{for
all }\;i=1,\ldots,m\;.\;\Box
\end{equation}
\end{lemma}
\begin{remark}
All assertions concerning relations between Witt vectors or
properties of $\bs{\pi}$-exponentials (see below) will be proved
by translating these relations or properties in terms of phantom
components.
\end{remark}
\subsubsection{\textbf{Frobenius and Verschiebung}}\label{frob and versc}
We denote by $\F:\W(R)\to\W(R)$ \index{F1@$\F:\W(R)\to\W(R)$,
$\F:\W_{m+1}(R)\to\W_{m}(R)$}
 and $\V:\W(R)\to\W(R)$ \index{V@$\V:\W(R)\to\W(R)$, $\V:\W_m(R)\to\W_{m+1}(R)$}
 the usual Frobenius and Verschiebung morphisms.
%\begin{center}
%\begin{tabular}{|l|l|l|}
%\hline
%  % \hline or \cline{col1-col2} \cline{col3-col4} ...
%  Name& Witt Vector & Phantom Vector  \\
%  \hline
%%  $\bs{r}$ & $(r_0,r_1,\ldots)$ & $\ph{\phi_0,\phi_1,\ldots}$  \\
%  $\V(r_0,r_1,\ldots)$ & $(0,r_0,r_1,\ldots)$ & $\ph{0,p\phi_0,p\phi_1,\ldots}$  \\
%  $\F(r_0,r_1,\ldots)$ & $(r_0^p+pr_1,?,?,\ldots)$& $\ph{\phi_1,\phi_2,\ldots}$\\
%  $\F\V(r_0,r_1,\ldots)=p\cdot(r_0,r_1,\ldots)$&$(pr_0,?,?,\ldots)$& $\ph{p\phi_0,p\phi_1,\ldots}$\\
%\hline
%\end{tabular}
%\end{center}
\label{FV=p} We denote again by $\F:\W_{m+1}(R)\to\W_{m}(R)$,
$\V:\W_m(R)\to\W_{m+1}(R)$ the reduction of $\F$ and $\V$ to
$\W_m(R)$. We have again $\F\V(\bs{r})=p\cdot\bs{r}$ in $\W_m(R)$.
We recall that
$\phi(\V(r_0,r_1,\ldots))=\ph{0,p\phi_0,p\phi_1,\ldots}$ and
$\phi(\F(r_0,r_1,\ldots))=\ph{\phi_1,\phi_2,\ldots}$.

If $R$ has characteristic $p$, then
$\F(r_0,r_1,\ldots)=(r_0^p,r_1^p,\ldots)$. Hence it is possible to
reduce the morphism $\F$ of $\W(R)$ to a morphism of $\W_m(R)$
into itself, by setting
$\Fb(r_0,\ldots,r_m)=(r_0^p,\ldots,r_m^p)$. We denote this
morphism by\index{F3@$\Fb:\W_m(R) \to \W_m(R)$ Frobenius (exist
only if $p=0$ in $R$)} $\Fb : \W_m(R) \to \W_m(R)$.
\subsubsection{\textbf{Completeness}}
Let $R$ be a topological ring. We identify topologically $\W_m(R)$
with $R^{m+1}$, via the function
$(r_0,\ldots,r_m)\mapsto\ph{r_0,\ldots,r_m}$. The operations on
$\W_m(R)$ are continuous, because defined by polynomials.
\begin{lemma}\label{completeness of W_m(R^0)}
If $R$ has a basis $\mathcal{U}_R$ of neighborhood of $0$ formed
by ideals, then $R$ is complete if and only if $\W_m(R)$ is
complete for all $m\geq 0$.
\end{lemma}
\emph{Proof : } It is evident for $m=0$. Let $m\geq 1$ and
$\{\bs{r}_n\}_n$, $\bs{r}_n:= (r_{n,0},\ldots,r_{n,m})$, be a
Cauchy sequence in $\W_m(R)$. The sequence $r_{0,n}$ is Cauchy in
$R$ and we denote $r_0:=\lim_nr_{0,n}$. The translate sequence
$\bs{r}^{1}_n:=\bs{r}_n-(r_0,0,\ldots,0)$ is Cauchy, so we can
suppose $r_0=0$. For every ideal $I\in\mathcal{U}_R$ there exists
$n_I$ such that $\bs{r}_{n_1}^1-\bs{r}_{n_2}^1=
(S_{0,n_1,n_2},,\ldots,S_{m,n_1,n_2})\in \W_m(I)$, for all
$n_1,n_2\geq n_I$. Let us write
$S_{1,n_1,n_2}=r_{n_1,1}^1-r_{n_2,1}^1+
P(r_{n_1,0}^1,r_{n_2,0}^1)$. By \cite[$\S 1$ $n^0 3$ $a)$]{Bou}
the polynomial $S_{k,n_1,n_2}$ is isobaric without constant term.
Since $r_{n,0}^1\in I$, for $n\geq n_{I}'$, sufficiently large and
since $I$ is an ideal, hence $r_{n_1,1}^1-r_{n_2,1}^1\in I$, for
all $n_1,n_2\geq n_{I}'$. So the sequence $r_{n,1}^1$ is Cauchy
and converges to $r_1\in R$. Moreover the sequence
$\bs{r}^{2}_n:=\bs{r}_n^1-(0,r_1,0,\ldots,0)$ is such that both
$r_{n,0}^{2}$ and $r_{n,1}^{2}$ go to $0$. This process can be
iterated indefinitely.\CVD
\begin{corollary} \label{completeness of W_m(R)}
If $(R,|.|)$ is an ultrametric valued ring, then $R$ is complete
if and only if $\W_m(R)$ is complete for all $m\geq 0$.
\end{corollary}
\subsubsection{\textbf{Length}}\label{ell}
Let $R$ be a ring of characteristic $p$. If the vector $\bs{r} \in
\W_m(R)$ is such that $r_0=\ldots=r_{k-1}=0$ and $r_{k}\neq 0$,
then we define the length of $\bs{r}$ as
%\begin{equation}
$\index{l@$\ell$} \ell(\bs{r}) := m-k$, %\;,
%\end{equation}
and $\ell(\bs{0}):=-\infty$. If $R$ is not of characteristic $p$,
then we will define $\ell(\bs{r})$ as the length of the image of
$\bs{r}$ in $\W_m(R/pR)$.
\subsubsection{\textbf{Covectors}}
We recall that the covectors module $\mathbf{CW}(R)$
\index{CW@$\mathbf{CW}(R)$, $\CW(R)$} is the additive group
defined by the following inductive limit (\cite[$\S 1$ ex. $23$
pag.47]{Bou}): $\mathbf{CW}(R):=
\varinjlim(\W_m(R)\stackrel{\V}{\to}\W_{m+1}(R)
\stackrel{\V}{\to}\cdots)$. \label{covectors tilde} In the sequel
we must work with a slightly different sequence. Let $R$ be a ring
of characteristic $p$. Then $\V\Fb=\Fb\V$ and
$\V\Fb(r_0,\ldots,r_m)=(0,r_0^p,\ldots,r_m^p)$. Let $\CW(R)$ be
the following inductive limit:
\begin{equation}
\CW(R):=\varinjlim(\W_m(R) \stackrel{\V\Fb}{\to}
\W_{m+1}(R)\stackrel{\V\Fb}{\to} \cdots)\;.
\end{equation}
If $R$ is a perfect field of characteristic
$p$, then $\mathbf{CW}(R)$ is isomorphic to $\CW(R)$. This results
from the following commutative diagram:
\begin{equation}\label{covectors tilde: the case perfect}
\xymatrix{ R\ar[r]^-{\V} \ar@{=}[d] \ar@{}[dr]|{\odot} &\W_1(R)
\ar[r]^{\V} \ar[d]^{\Fb}_{\wr} \ar@{}[dr]|{\odot}&\W_2(R)
\ar[r]^-{\V}\ar[d]^{\Fb^2}_{\wr} & \cdots\ar[r]&
\mathbf{CW}(R)\;\phantom{.}\ar[d]_{\wr} \\
R\ar[r]^-{\V\Fb} & \W_1(R)\ar[r]^{\V\Fb} & \W_2(R)\ar[r]^-{\V\Fb}
& \cdots\ar[r]&\CW(R)\;. }
\end{equation}
\begin{remark} \label{CW(R)=CW(R^p)}
If $R$ is a field of characteristic $p$, then $\CW(R)=\CW(R^p)$.
\end{remark}

\subsection{Notations in Artin-Schreier  theory}\label{Notation in
AS theory}
\begin{definition}
\label{G_R, I_R, P_R}Let $R$ be a field of characteristic $p>0$
and let $R^{\mathrm{sep}}/R$ be a fixed separable closure of $R$.
 We denote by
$\G_R=\mathrm{Gal}(R^{\mathrm{sep}}/R)$.
\index{G_R@$\G_R=\mathrm{Gal}(R^{\mathrm{sep}}/R)$,
$\mathcal{I}_R=$inertia of $\G_R$, $\mathcal{P}_R=$wild inertia of
$\G_R$} If $R$ is a complete discrete valuation field, we denote
by $\mathcal{I}_R$ the inertia group and by $\mathcal{P}_R$ the
pro$-p$-sylow subgroup of $\mathcal{I}_R$.
\end{definition}
We have
$\textrm{H}^1(\G_R,\mathbb{Z}/p^{m}\mathbb{Z})\stackrel{\sim}{\to}
\Hom^{\mathrm{cont}}(\G_R,\mathbb{Z}/p^{m}\mathbb{Z})$ (cf.
\cite[Ch.$X$, $\S 3$]{Se}). The situation is then expressed by the
following commutative diagram:\index{delta@$\delta=$ cohomological
morphism of the Artin-Schreier  complex}
\begin{equation}\label{artin-screier-diagram}
\begin{scriptsize}
\xymatrix{ 0\ar[r]&\mathbb{Z}/p^{m+1}\mathbb{Z}\ar[r]
\ar@{}[dr]|{\odot}\ar[d]^{\imath} &
\W_m(R)\ar[r]^{\Fb-1}\ar[d]^{\V}\ar@{}[dr]|{\odot}&
\W_m(R)\ar[r]^-{\delta}\ar@{}[dr]|{\odot}\ar[d]^{\V}&
\Hom^{\mathrm{cont}}(\G_R,\mathbb{Z}/p^{m+1}\mathbb{Z})
\ar[d]^{\jmath}\ar[r]&0\\
0\ar[r]&\mathbb{Z}/p^{m+2}\mathbb{Z}\ar[r]&
\W_{m+1}(R)\ar[r]^{\Fb-1}&\W_{m+1}(R)\ar[r]^-{\delta}&
\Hom^{\mathrm{cont}}(\G_R,\mathbb{Z}/p^{m+2}\mathbb{Z})\ar[r]&0
}
\end{scriptsize}
\end{equation} where $\imath:1\mapsto p$ is the usual inclusion,
and $\jmath$ is the composition with $\imath$. For $\lb\in
\W_m(R)$, the character $\alpha=\delta(\lb)$
\index{alpha@$\alpha=\delta(\lb)=$ Artin-Schreier  character}
sends the automorphism $\gamma$ to the element
$\alpha(\gamma):=\gamma(\bs{\nu})-\bs{\nu}\in
\mathbb{Z}/p^{m+1}\mathbb{Z}$, where $\bs{\nu}\in
R^{\mathrm{sep}}$ is a solution of the equation
$\Fb(\bs{\nu})-\bs{\nu}=\lb$. Taking the inductive limit, we get
the following exact sequence:
\begin{equation}\label{artin-schreier-diagram-covectors}
0\to\mathbb{Q}_p/\mathbb{Z}_p
\to\bs{\mathrm{CW}}(R)\xrightarrow[]{\Fb-1}\bs{\mathrm{CW}}(R)\to
\Hom^{\textrm{cont}}(\G_R,\mathbb{Q}_p/\mathbb{Z}_p)\to 0\;,
\end{equation}
where the word ``$\textrm{cont}$'' means that all characters
$\G_R\to\mathbb{Q}_p/\mathbb{Z}_p$ factorize on a finite quotient
of $\G_R$. Indeed
$\varinjlim_m\Hom(\G_R,\mathbb{Z}/p^m\mathbb{Z})$ can be seen as
the subset of $\Hom(\G_R,\mathbb{Q}_p/\mathbb{Z}_p)$ formed by the
elements killed by a power of $p$.
\begin{remark}\label{VF is the same of F}
If the vertical arrows $\V$ are replaced by $\V\Fb$ in the diagram
\eqref{artin-screier-diagram}, then the morphisms $\imath$ and
$\jmath$ remain the same. Indeed $\delta(\lb)=\delta(\Fb(\lb))$,
because $\Fb(\lb) = \lb + (\Fb-1)(\lb)$, for all $\lb\in\W_s(R)$.
Hence we have also
\begin{equation}
0\to\mathbb{Q}_p/\mathbb{Z}_p
\to\CW(R)\xrightarrow[]{\Fb-1}\CW(R)\to
\Hom^{\textrm{cont}}(\G_R,\mathbb{Q}_p/\mathbb{Z}_p)\to 0\;.
\end{equation}
\end{remark}
\begin{remark}\label{R(nu_0,...,nu_m)}
Let $\lb=(\lambda_0,\ldots,\lambda_m)\in\W_m(R)$. The kernel of
$\alpha:=\delta(\lb)$ is the subgroup of $\G_R$ whose
corresponding extension field is $R(\{\nu_0,\ldots,\nu_m\})$,
(i.e. the smallest field containing the set
$\{\nu_0,\ldots,\nu_m\}$), where $\bs{\nu}=(\nu_0,\ldots,\nu_m)\in
\W_m(R^{\mathrm{sep}})$ is solution of
$\Fb(\bs{\nu})-\bs{\nu}=\lb$. All cyclic separable extensions of
$R$, whose degree is a power of $p$, are of this form for a
suitable $m\geq 0$, and $\lb$.
\end{remark}
\subsubsection{}
Let $\kappa$ be a field of characteristic $p>0$, and let
$R:=\kappa(\!(t)\!)$. The Galois group of an abelian extension of
$\kappa(\!(t)\!)$ is the product of its $p$-torsion part
(controlled by the Artin-Schreier  theory) and its moderate part
(controlled by Kummer theory).
\begin{definition} \label{definition of PASS}
We set \index{PAS@$\PAS(\kappa):=
\Hom^{\mathrm{cont}}(\mathcal{P}_R,
 \mathbb{Q}_p/\mathbb{Z}_p)$}
 $\PAS(\kappa):= \Hom^{\mathrm{cont}}(\mathcal{P}_R,
 \mathbb{Q}_p/\mathbb{Z}_p)=
 \Hom^{\mathrm{cont}}(\mathcal{I}_R,
 \mathbb{Q}_p/\mathbb{Z}_p)$.
\end{definition}
\begin{remark} \label{first remark on PAS}
We will see that $\PAS(\kappa)\cong
\frac{\bs{\mathrm{CW}}(t^{-1}\kappa[t^{-1}])}{(\Fb-1)
\bs{\mathrm{CW}}(t^{-1}\kappa[t^{-1}])}$. This group describes the
abelianization of the pro-$p$-Sylow of the quotient
$\mathcal{P}_R$. On the other hand
$\bs{\mathrm{CW}}(\kappa)/(\Fb-1)\bs{\mathrm{CW}}(\kappa)=
\Hom^{\textrm{cont}}(\mathrm{G}_R/\mathcal{I}_R,
\mathbb{Q}_p/\mathbb{Z}_p)$.\end{remark}

\subsection{Notations in Lubin-Tate theory}\label{Lubin-Tate
theory} For notations and results on Lubin-Tate theory we refer to
\cite{L-T}. In this paper we will treat only Lubin-Tate groups
over the field $\mathbb{Q}_p$. We recall briefly only the facts
used in this paper. Let\index{w@$\mathrm{w}=$ uniformizer of
$\mathbb{Z}_p$} $\mathrm{w}:=p\cdot u \in p\mathbb{Z}_p$,
$u\in\mathbb{Z}_p^{\times}$, be a uniformizing element. Let
$\mathfrak{F}_\mathrm{w}$ be the family of formal power series
$P(X) \in \mathbb{Z}_p[[X]]$ satisfying\index{P(X)@$P(X)$
Lubin-Tate series}
\begin{equation}\label{Lubin-Tate series}
P(X)   \equiv   \mathrm{w} X
\pmod{X^2\mathbb{Z}_p[[X]]}\quad,\quad P(X)   \equiv        X^p
\pmod{\mathrm{w}\mathbb{Z}_p[[X]]}.
\end{equation}
A series in $\mathfrak{F}_\mathrm{w}$ will be called a
\emph{Lubin-Tate series}. For all $P\in\mathfrak{F}_\mathrm{w}$,
there exists a unique formal group law $\mathfrak{G}_P(X,Y)\in
\mathbb{Z}_p[[X,Y]]$
\index{G_P(X,Y)@$\mathfrak{G}_P(X,Y)=$Lubin-Tate formal group}
such that $P(\mathfrak{G}_P(X,Y))=\mathfrak{G}_P(P(X),P(Y))$ (i.e.
$P(X)$ is an endomorphism of $\mathfrak{G}_P(X,Y)$).
\begin{lemma} \label{[a]_{P,widetilde{P}}}
Let $P,\widetilde{P}\in\mathfrak{F}_\mathrm{w}$. For all
$a\in\mathbb{Z}_p$ there exists a unique formal series
$[a]_{P,\widetilde{P}}(X)\in \mathbb{Z}_p[[X]]$ such that
\begin{enumerate}
\item $[a]_{P,\widetilde{P}}(X)\equiv
aX\pmod{X^2\mathbb{Z}_p[[X]]}$, \item
$[a]_{P,\widetilde{P}}(\mathfrak{G}_P(X,Y))=
\mathfrak{G}_{\widetilde{P}}([a]_{P,\widetilde{P}}(X),[a]_{P,\widetilde{P}}(Y))$.
\end{enumerate}
In other words $[a]_{P,\widetilde{P}}(X)$ is a morphism of group
laws. \CVD
\end{lemma}
We set $\index{a_P@$[a]_{P,\widetilde{P}}(X)$, $[a]_P(X)$}
[a]_{P}(X):=[a]_{P,P}(X)$. By the uniqueness, we have that
$P(X)=[\mathrm{w}]_P(X)$. The setting $x*y:=\mathfrak{G}_P(x,y)$
defines a new group law on $\mathfrak{p}_K$, denoted by
$\mathfrak{G}_P(\mathfrak{p}_K)$.
%We denote this group by
%$\mathfrak{G}_P(\mathfrak{p}_K)$.
\label{P^(k)} Let $P^{(k)}$ denote the series $P\circ
P\circ\cdots\circ P$, $k$-times. Following \cite{L-T} let
 \begin{equation}
 \Lambda_{P,m}=\Ker(P^{(m)})=\Ker([\mathrm{w}^{m}]_P)=
 \{x\in\mathbb{C}_p\;|\;P^{(m)}(x)=0\textrm{ and }
 |x|<1\}
 \end{equation}
be the set of $[\mathrm{w}]^m$-torsion points of
$\mathfrak{G}_{P}(\mathfrak{p}_{\mathbb{C}_p})$, and
$\Lambda_P:=\cup_m\Lambda_{P,m}$. We have $\Lambda_{P}\subset
\mathbb{Q}_p^{\mathrm{alg}}$. Moreover
$\mathbb{Q}_p(\Lambda_{P,m})/\mathbb{Q}_p$ is Galois and depend
only on $\mathrm{w}$. The formal group law $\mathfrak{G}$ makes
$\Lambda_{P,m}$ a group.

\begin{theorem}[\protect{\cite[Th.2]{L-T}}] \label{reciprocity and filtrations}
We have the following properties:
\begin{enumerate}
 \item We have $\Lambda_{P,m}\cong \mathbb{Z}/p^m\mathbb{Z}$, for all $m\geq 0$,
 and then $\Lambda_P\cong \mathbb{Q}_p/\mathbb{Z}_p$.
 \item \index{gamma@$\gamma=$ automorphism} Let $\gamma\in
\mathrm{Gal}(\mathbb{Q}_p(\Lambda_P)/\mathbb{Q}_p)$. There exists
an unique unit $u_\gamma\in\mathbb{Z}_p^\times$ such that
$$\gamma(x) = [u_\gamma]_P(x)\;,\quad\forall x\in\Lambda_{P}.$$
 \item The map $\gamma\mapsto u_\gamma$ is an isomorphism of
$\mathrm{Gal}(\mathbb{Q}_p(\Lambda_P)/\mathbb{Q}_p)$ onto the
group $\mathbb{Z}_p^{\times}$. The same map gives an isomorphism
$$\mathrm{Gal}(\mathbb{Q}_p(\Lambda_P)/\mathbb{Q}_p(\Lambda_{P,m}))
\stackrel{\sim}{\to} 1+\mathrm{w}^m\mathbb{Z}_p\;,\quad \forall\;
m\geq 1.$$
 \item Let $u\in\mathbb{Z}_p^{\times}$, then $[u]_{P}(x)=(u^{-1},
\mathbb{Q}_p(\Lambda_{P,m})/\mathbb{Q}_p)(x)$, for all
$x\in\Lambda_{P,m}$, where $(u^{-1},
\mathbb{Q}_p(\Lambda_{P,m})/\mathbb{Q}_p)\in
\mathrm{Gal}(\mathbb{Q}_p(\Lambda_{P})/\mathbb{Q}_p)$ is the norm
residue symbol.
\end{enumerate}
\end{theorem}
\begin{remark}\label{(xi_m-1)}
The simplest Lubin-Tate series is $P(X)=\mathrm{w} X+X^p$. If
$\mathrm{w}=p$, then a non trivial zero $\pi_0$ of $P$ is the
``$\pi$'' of Dwork. If again $\mathrm{w}=p$ and $P(X)=(X+1)^p-1$,
then $\mathfrak{G}_{P}\cong\widehat{\mathbb{G}}_m$, and all
torsion points are of the form $\xi-1$, with $\xi^{p^k}=1$, for
some $k\geq 0$. This was the choice made by Matsuda \cite{Ma}.
\end{remark}
\begin{theorem}[\protect{\cite[Prop. 8.3.22]{Haz}}]\label{moduly space of Lubin Tate Groups}
Let $\mathfrak{G}$ and $\widetilde{\mathfrak{G}}$ be two
Lubin-Tate groups relative to the uniformizers $\mathrm{w}$ and
$\tilde{\mathrm{w}}$ respectively. Then $\mathfrak{G}$ is
isomorphic to $\widetilde{\mathfrak{G}}$ (as formal groups over
$\mathbb{Z}_p$) if and only if $\mathrm{w}=\tilde{\mathrm{w}}$.
$\Box$
\end{theorem}

\subsubsection{\textbf{Tate module}}
The multiplication by $[\mathrm{w}]_P$ sends $\Lambda_m$ into
$\Lambda_{m-1}$. The Tate module of $\mathfrak{G}_P$ is, by
definition, $\T(\mathfrak{G}_P):=\varprojlim_{m}\Lambda_{P,m}$.
\index{T(G_P)@$\T(\mathfrak{G}_P):=\varprojlim_{m}\Lambda_{P,m}$}
A generator $\bs{\pi}= (\pi_{P,j})_{j\geq 0}$ \index{pi@$\bs{\pi}=
(\pi_{P,j})_{j\geq 0}=$ generator of $\T(\mathfrak{G}_P)$, where
$\pi_j = \mathrm{w}^{j+1}$-torsion point of $\mathfrak{G}_P$} of
the Tate module $\T(\mathfrak{G}_P)$ is a sequence
$(\pi_{P,j})_{j\geq 0}$, $\pi_j\in \Lambda_P$, such that
$P(\pi_{P,0})=0$, $\pi_{P,0}\neq 0$ and
$P(\pi_{P,j+1})=\pi_{P,j}$, for all $j\geq 0$. If no confusion is
possible, we will write $\pi_j$ instead of $\pi_{P,j}$. The Newton
polygon of $P$ shows that $P$ has exactly $p-1$ non trivial zeros
of value $\omega=|p|^{\frac{1}{p-1}}$, and inductively
$P(X)-\pi_{j-1}$ has $p$ zeros of valuation
$\omega^{\frac{1}{p^j}}$. Hence
$|\pi_j|=\omega^{1/p^{j}}$, %\;,\quad\forall\;
for all $j\geq 0$, and the Galois extension
$\mathbb{Q}_p(\Lambda_{P,m})=\mathbb{Q}_p(\pi_{m-1})$ is totally
ramified. On the other hand the field $K(\pi_{m-1})$ is not always
totally ramified.

\begin{definition} \label{L_infty}
We set $K_m:=K(\pi_m)$\index{K_m@$K_m:=K(\pi_m)$, $K_\infty=\cup_m
K_m$, $k_m$, $k_\infty$} (resp. $K(\Lambda_P)$), and denote by
$k_{m}$ (resp. $k_{\mathrm{w}}$) its residue field. Moreover, if
$\mathrm{w}=p$, we put $K_\infty:=K(\Lambda_P)$ and
$k_\infty:=k_p$. For all algebraic extensions $L/K$, $L_\infty$
will be the smallest field containing $L$ and $K_\infty$.
\end{definition}
\begin{example}
If $P(X)=(X+1)^p-1$, then $\mathfrak{G}_P=\widehat{\mathbb{G}}_m$,
and $\Lambda_m=\{\xi_m-1\;|\;\xi_m^{p^{m}}=1\}$ is the set of
$p^m$-th root of $1$. \index{xi@$\xi_m=p^{m+1}$-th primitive root
of $1$} A generator of $\T(\widehat{\mathbb{G}}_m)$ is a family
$(\xi_j-1)_{j\geq 0}$ satisfying $\xi_j^{p^i}=\xi_{j-i}$, for all
$0\leq i\leq j$.
\end{example}
\begin{definition}\label{equiv classe-}
Let $P,\widetilde{P}\in\mathfrak{F}_\mathrm{w}$ be two Lubin-Tate
series. We will say that $x\in\Lambda_{P}$ and
$y\in\Lambda_{\widetilde{P}}$ are equivalent if
$y=[1]_{P,\widetilde{P}}(x)$ (cf. \ref{[a]_{P,widetilde{P}}}).
\end{definition}
\begin{remark} \label{relations on equivalence classe}
Since $[1]_{P,\widetilde{P}}(x) = x + (\textrm{things divisible by
}x^2)$, it follows that $|x-[1]_{P,\widetilde{P}}(x)|\leq |x|^2$.
In particular, if $\mathrm{w}=p$ and if $\pi_m$ is fixed, then
there exists a unique $p^{m+1}$-th root of $1$, say $\xi_m$, such
that $|(\xi_m-1)-\pi_m|\leq\omega^{\frac{2}{p^m}}$ (cf.
\ref{(xi_m-1)}).
\end{remark}

\specialsection{\textbf{$\bs{\pi}$-exponentials and applications}}

\subsection{Construction of Witt vectors}\label{P-sequences and
Witt vectors} \label{Construction of Witt vectors}
 Let $P(X)\in \mathbb{Z}_p[[X]]$ be a series, with $P(0)=0$, satisfying
\begin{equation}\label{P(X)iii}
P(X)\equiv X^p\mod{p\mathbb{Z}_p[[X]]}.
\end{equation}
We consider the Frobenius
$\sigma_P:\mathbb{Z}_p[[X]]\to\mathbb{Z}_p[[X]]$ given by
$\sigma_P(h(X)):=h(P(X))$.
\begin{lemma}[\protect{\cite[Ch.IX,$\S 1$,ex.$14$,a)]{Bou}}]\label{lemma 6.1}
There is a unique ring morphism
\begin{equation}
[-]:\mathbb{Z}_p[[X]]\xrightarrow[\protect{h(X)\mapsto[h(X)]}]{}\W(\mathbb{Z}_p[[X]])
\end{equation}
such that $\phi_j\circ [-]=\sigma_{P}^j$. In other words, for all
$h(X)\in\mathbb{Z}_p[[X]]$, the Witt vector $[h(X)]$ is the unique
one whose phantom vector is equal to
\begin{equation}
\ph{\;h(X)\;,\;h(P(X))\;,\;h(P(P(X)))\;,\;\ldots\;}.
\end{equation}
Moreover $[-]$ is also the unique ring morphism satisfying the
relation
\begin{equation}\label{F^m lambda =0}
\F([h(X)])=[h(P(X))].
\end{equation}
\end{lemma}
\emph{Proof : } By Lemma \ref{congruece on the phantom
components}, the ring morphism
$h(X)\mapsto\ph{h(X),h(P(X)),\ldots}:\mathbb{Z}_p[[X]]\to(\mathbb{Z}_p[[X]])^{\mathbb{N}}$
has its values in the image of the phantom component map
$\phi:\W(\mathbb{Z}_p[[X]])\hookrightarrow(\mathbb{Z}_p[[X]])^{\mathbb{N}}$.
Since, by \ref{phi is inj}, $\phi$ is injective, the lemma is
proved. $\Box$
\begin{definition}\label{[h(pi_m)]}
Let $\B$ be a complete topologized $\mathbb{Z}_{p}$-ring, and let
$b\in \B$ be a topologically nilpotent element. The specialization
$X\mapsto b:\mathbb{Z}_p[[X]]\to \B$ provides, by funtoriality, a
morphism $\W(\mathbb{Z}_p[[X]])\to\W(\B)$. For brevity, we denote
by $[h(b)]$ the image of $h(X)$ via the morphism
\begin{equation}\label{s_P^b}\index{hb@$[h(b)]$}
\mathbb{Z}_{p}[[X]]\xrightarrow[]{[-]}\W(\mathbb{Z}_p[[X]])\xrightarrow[]{X\mapsto
b}\W(\B).
\end{equation}
We will denote again by $[h(b)]$ its image in $\W_m(\B)$.
\end{definition}
\begin{remark}
The phantom vector of $[h(b)]$ is
\begin{equation}\label{phi_j of s_sigma}
\ph{\;h(b)\;,\;h(P(b))\;,\;h(P(P(b)))\;,\;\ldots\;}.
\end{equation}

In general there is no morphism $\mathbb{Z}_p[b]\to\W(\B)$ sending
$h(b)$ into $[h(b)]$, the notation $[h(b)]$ is imprecise, but more
handy.
\end{remark}

\begin{lemma}[key lemma]
Let $(\B,|.|)$ be a $\mathbb{Z}_p$-ring, complete with respect to
an absolute value $|.|$, extending the absolute value of
$\mathbb{Z}_p$. Let $h(X)=\sum_{i\geq 0}a_iX^i\in
\mathbb{Z}_p[[X]]$, and let
$[h(b)]=(\lambda_0,\lambda_1,\ldots)\in\W(\B)$, with $|b|<1$. Then
the following statements are equivalent:
\begin{enumerate}
    \item $|a_0|=|p|^{r}$ ,
    \item  $|\lambda_0|,\ldots,|\lambda_{r-1}|<1$, and $|\lambda_r|=1$.
\end{enumerate}
\end{lemma}
\emph{Proof : } Let $\lb=(\lambda_0,\lambda_1,\ldots)=[h(b)]$. We
denote by $\overline{\B}$ the residue field. The condition $(2)$
is equivalent to $\bar{\lambda}_r\neq 0$, and $\bar{\lambda}_i=0$,
for all $i<r$, or, if $k\geq 0$ is given, it is equivalent to
$\bar{\lambda}_r^{p^k}\neq 0$, and $\bar{\lambda}_i^{p^k}=0$, for
all $i<r$. This last condition is equivalent to the condition
$(2)$ for the vector $\F^k(\lb)$. Now the phantom vector of
$\F^k(\lb)$ is $\ph{h(P^{(k)}(b)),h(P^{(k+1)}(b)),\ldots}$ (cf.
\ref{P^(k)}). Moreover $|P(b)|\leq \sup(|b|^{p},|p||b|)$, hence,
for all $\varepsilon>0$, there exists $k\geq 0$ such that
$|P^{(i)}(b)|<\varepsilon$, for all $i\geq k$. If $\varepsilon$ is
small enough, then $|h(P^{(i)}(b))|=|a_0|$, for all $i\geq k$. Let
$(\nu_0,\nu_1,\ldots):=\F^k(\lb)$, then, since
$p^j\nu_j=h(P^{(j)}(b))-(\nu_0^{p^j}+\cdots+p^{j-1}\nu_{j-1}^p)$,
we see, by induction, that $|a_0| = |p|^r$ if and only if
$|\nu_r|=1$ and $|\nu_j| = |p|^{r-j}$, for all $j\leq r-1$.$\Box$

\begin{definition}\label{varpi_m}
We fix now a sequence $\bs{\varpi}:=\{\varpi_j\}_{j\geq 0}$ in
$\mathbb{Q}_p^{\mathrm{alg}}$ satisfying $|\varpi_0|<1$,
$P(\varpi_0)=0$, and $P(\varpi_{j+1})=\varpi_j$, for all $j\geq
1$.
\end{definition}
\begin{remark}
The ring $\mathbb{Z}_p[\varpi_m]$ is complete, for all $m\geq 0$.
Indeed $\varpi_m$ is algebraic and integral over $\mathbb{Z}_p$,
hence $\mathbb{Z}_p[\varpi_m]$ is a free module over
$\mathbb{Z}_p$.
\end{remark}
\begin{remark}
%Let $P(X)$ be the series fixed at \eqref{P(X)iii}, and let
%$\bs{\pi}$ be the sequence defined at \ref{varpi_m}.
If $P$ is a Lubin-Tate series, and if $\varpi_0\neq 0$, then
$\bs{\varpi}$ is a generator of the Tate module
$\mathrm{T}(\mathfrak{G}_P)$, while if $P(X)\equiv
X^p\mod{X^{p+1}\mathbb{Z}_p[[X]]}$, then $\varpi_j=0$ for all
$j\geq 0$. Observe that, taking $h(T):=T$ and $b:=\varpi_m$ in the
lemma \ref{lemma 6.1}, then
$[\varpi_m]\in\W(\mathbb{Z}_p[\varpi_m])$ is the unique Witt
vector whose phantom vector is
$\ph{\varpi_m,\varpi_{m-1},\ldots,\varpi_{0},0,\ldots}$. The
uniqueness follows from the injectivity of the phantom map
$\phi:\W(\mathbb{Z}_p[\varpi_m])\hookrightarrow(\mathbb{Z}_p[\varpi_m])^{\mathbb{N}}$.
\end{remark}

\begin{proposition}
For all $\mathbb{Z}_{p}[\varpi_{m}]$-algebra $\B$ of
characteristic $0$
\begin{equation}\label{[pi_m]W(B)subset [pi_m+1]W(B)}
 [\varpi_{j}]\W(\B)\subset [\varpi_{j+1}]\W(\B)\;,\qquad j=0,\ldots,
 m-1\;.
\end{equation}
Moreover, for all $\lb\in\W(\B)$, and all $j=0,\ldots, m-1$
\begin{equation}\label{V([pi_m+1]lb)=[pi_m]V(lb)}
\F([\varpi_{j+1}])=[\varpi_j]\qquad;\qquad
\V([\varpi_{j}]\cdot\lb)=[\varpi_{j+1}]\cdot\V(\lb).
\end{equation}
Hence $\F([\varpi_{j+1}]\W(\B))\subset [\varpi_{j}]\W(\B)$ and
$\V([\varpi_j]\W(\B))\subset[\varpi_{j+1}]\W(\B)$.

If now $\varpi_0\neq 0$, then the kernel of the morphism
$\lb\mapsto[\varpi_m]\lb$ is the ideal $\V^{m+1}\W(\B)$. The
induced morphism $\W_m(\B)\to\W(\B)$ is a functorial isomorphism
of $\W_m(\B)$ into the ideal $[\varpi_m]\W(\B)$ (as
$\W(\B)$-modules), which commutes with
$\V:\W_m(\B)\to\W_{m+1}(\B)$ and $\F:\W_{m+1}(\B)\to\W_m(\B)$
\begin{equation}\label{W_m simto [pi_m]W}
 \xymatrix{
 \W(\B)\ar@{->>}[d]\ar[r]^-{\lb\mapsto[\varpi_m]\lb}&\W(\B)\\
 \W_m(\B)\ar[r]^-{\sim}&[\varpi_m]\cdot\W(\B)\;.\ar@{^(->}[u]
 }
\end{equation}
\end{proposition}
\emph{Proof : } Let $h(X):=P(X)/X$, then
$[\varpi_{j}]\lb=[P(\varpi_{j+1})]\lb=[\varpi_{j+1}\cdot
h(\varpi_{j+1})]\lb=[\varpi_{j+1}]\cdot[h(\varpi_{j+1})]\lb$. This
shows the inclusion \eqref{[pi_m]W(B)subset [pi_m+1]W(B)}. All
other assertions are easily verified on the phantom components.
$\Box$

\begin{corollary} \label{base lemma!}
Let $\ph{\phi_0,\ldots,\phi_m}\in\mathbb{Z}_p[\varpi_m]^{m+1}$. If
there exists a formal series $h(X) = \sum_{i\geq 0} a_i X^i \in
\mathbb{Z}_p[[X]]$ satisfying
\begin{equation}
h(\varpi_{m-j})=\phi_j\;,\quad\textrm{for all }0\leq j\leq m,
\end{equation}
then $\ph{\phi_0,\ldots,\phi_m}$ is the phantom vector of
$[h(\varpi_m)]:=(\nu_0,\ldots,\nu_m)\in\W_m(\mathbb{Z}_p[\varpi_m])$.
Moreover, $|a_0|=|p|^r$, for some $r\geq 0$, if and only if
$|\nu_0|,\ldots,|\nu_{r-1}|<1$ and $|\nu_r|=1$. $\Box$
\end{corollary}

\subsubsection{\textbf{Artin-Hasse exponential and Robba exponentials}}

\begin{definition}[\protect{\cite[ex.58]{Bou}}]
\label{prod_j geq 0 E(lambda_j
T^p^j)=exp(phi_0T+phi_1T^p/p+cdots)}\index{E(T)@$E(T)=\exp(\sum_{j=0}^\infty
\frac{T^{p^j}}{p^j})$, $E(\lb,T)$} Let $\B$ be a
$\mathbb{Z}_{(p)}$-ring, and let
\begin{equation}
E(T):=\exp(T+\frac{T^p}{p}+\frac{T^{p^2}}{p^2}+\cdots)\in
1+T\mathbb{Z}_{(p)}[[T]]\;.
\end{equation}
For all $\lb:=(\lambda_0,\lambda_1,\ldots)\in\W(\B)$, the
Artin-Hasse exponential relative to $\lb$ is
\begin{equation}
E(\lb,T):=\prod_{j\geq 0}E(\lambda_j\cdot T^{p^j})=\exp(\phi_0
T+\phi_1\frac{T^{p}}{p}+\phi_2\frac{T^{p^2}}{p^2}+\cdots)\in 1+T
\B[[T]]\;,
\end{equation}
where $\ph{\phi_0,\phi_1,\ldots}$ is the phantom vector of $\lb$.
\end{definition}
\begin{remark}
The Artin-Hasse exponential is then a group morphism
\begin{equation}
E(-,T):\W(\B)\to 1+T\B[[T]],
\end{equation}
functorial on the $\mathbb{Z}_{(p)}$-ring $\B$.
\end{remark}
\begin{proposition}\label{E_m converge exactly for T<1 iff P=LT}
Let $[\varpi_m]\in\W(\mathbb{Z}_p[\varpi_m])$ be as in \ref{varpi_m}. The exponential
\begin{equation}\index{E_m(T)@$E_m(T):=\exp(\pi_{m}T+
\cdots+\pi_0\frac{T^{p^m}}{p^m})$}
E_m(T):=E([\varpi_m],T)=\exp(\varpi_mT+\varpi_{m-1}\frac{T^p}{p}+\cdots+\varpi_0\frac{T^{p^m}}{p^m})
\end{equation}
converges exactly in the disk $|T|<1$, for all $m\geq 0$, if and
only if $P(X)$ is a Lubin-Tate series, and
$\bs{\varpi}:=(\varpi_j)_{j\geq 0}$ is a generator of the Tate module
$\mathrm{T}(\mathfrak{G}_P)$.
\end{proposition}
\emph{Proof : } Assume that the radius of convergence of
$E([\varpi_m],T)$ is equal to $1$, for all $m\geq 0$. Then, for
$m=0$, the radius of convergence of $\exp(\varpi_0T)$ is $1$,
hence $|\varpi_0|=\omega$. The Newton polygon of $P(X)$ implies
that $P(X)\equiv\mathrm{w} X\mod{X\mathbb{Z}_p[[X]]}$, for some
$\mathrm{w}$, with $|\mathrm{w}|=|p|$, hence $P(X)$ is a
Lubin-Tate series. Conversely, assume that $P(X)$ is a Lubin-Tate
series, and that $\bs{\varpi}:=(\varpi_j)_{j\geq 0}$ is a
generator of $\mathrm{T}(\mathfrak{G}_P)$. Consider the
differential operator
$L:=\d+\varpi_mT^{-1}+\varpi_{m-1}T^{-p}+\cdots+\varpi_0T^{-p^m}$.
Then $E_m(T^{-1})$ is the Taylor solution at $+\infty$ of $L$.
Since $|\varpi_0|=\omega$, by Lemma \ref{graphic}, we have
$Ray(L,\rho)=\rho^{p^m+1}$, for all $\rho<1$. In particular, the
irregularity of $L$ is $p^m$. Then $E_m(T^{-1})$ is not convergent
for $|T|<1$, because otherwise, by transfer at $\infty$,
$E_m(T^{-1})\in\R$, and $L$ will be trivial. $\Box$
\begin{theorem}\label{over-convergence of E_m(T)/E_m(T)}
Let $P(X)=\mathrm{w} X+\cdots$ be a Lubin-Tate series, and let
$\bs{\varpi}:=(\varpi_j)_{j\geq 0}$ be a generator of
$\mathrm{T}(\mathfrak{G})$. Then the formal series
$E_m(T^p)/E_m(T)$ is over-convergent (i.e. convergent for
$|T|<1+\varepsilon$, for some $\varepsilon>0$) if and only if
$$|\mathrm{w}-p|\leq |p|^{m+2}\;.$$
In particular, $E_m(T^p)/E_m(T)$ is over-convergent for all $m\geq
0$ if and only if $\mathfrak{G}_P$ is isomorphic to the formal
multiplicative group $\widehat{\mathbb{G}}_m$ (cf. Theorem
\ref{moduly space of Lubin Tate Groups}).
\end{theorem}
\emph{Proof : } This Theorem will follow easily from the theory of
$\bs{\pi}$-exponentials (cf. the proof of \ref{overconvergence of
theta} infra), and is placed here for expository reasons.$\Box$

\subsection{$\bs{\pi}$-exponentials}\label{exponentials}

We maintain the notations of Section \ref{P-sequences and Witt
vectors}. In this section we fix a uniformizing element
$\mathrm{w}$ of $\mathbb{Z}_p$, a $\mathbb{Q}_p$-Lubin-Tate series
$P\in\mathbb{Z}_p[[X]]$, $P\in\mathfrak{F}_\mathrm{w}$, and a
generator $\bs{\pi}=(\pi_j)_{j\geq 0}$ of the Tate module. We fix
three natural numbers $n,m,d$ such that
\begin{equation}
d=n\cdot p^m > 0\;,\; \textrm{ and } (n,p)=1.
\end{equation}
\begin{definition}
Let $\B$ be a $\mathbb{Z}_p[\pi_m]$-algebra.  Let
$\lb=(\lambda_0,\ldots,\lambda_m)\in\W_m(\B)$, and let
$\ph{\phi_0,\ldots,\phi_m}\in\B^{m+1}$ be its phantom vector. We
set
\index{e_d(lambda,T)@$\mathrm{e}_d(\lb,T)=\exp(\pi_m\phi_0T^{n}+
\cdots + \pi_0\phi_m\frac{T^{d}}{p^m})$}
\begin{equation}\label{definition of et_d(lb,T)}
\mathrm{e}_d(\lb,T):=E([\pi_m]\lb,T^n)=\exp\bigl(\pi_m\phi_0T^{n}+
\pi_{m-1}\phi_1\frac{T^{np}}{p}+\cdots+
\pi_0\phi_m\frac{T^{d}}{p^m}\bigr)\;.
\end{equation}
We will call $\mathrm{e}_d(\lb,T)\in 1+\pi_m T\B[[T]]$ the
\emph{$\bs{\pi}$-exponential} attached to $\lb$.
\end{definition}

\begin{proposition}\label{rules}
The map $\lb\mapsto\mathrm{e}_d(\lb,T)$ defines a group morphism
\begin{equation}\label{et_d(-,T) has values in 1+pi_m T}
\W_m(\B)\longrightarrow 1+\pi_m T\B[[T]]\;.
\end{equation}
Moreover, for all $\lb,\bs{\nu}\in\W_m(\B)$, we have
\begin{equation}\label{et_d is a product}
\mathrm{e}_d(\lb,T)=\prod_{j=0}^{m}E_{m-j}(\lambda_{j}T^{np^j}),
\end{equation}
\begin{equation*}
\begin{array}{rclcrcl}
E_m(T)&=&\mathrm{e}_{p^m}((1,0,\ldots,0),T)&,\quad&\mathrm{e}_d(\lb,T)&=&\mathrm{e}_{p^m}(\lb,T^n)\;,\\
\mathrm{e}_d(\lb,T^p)&=&\mathrm{e}_{p\cdot d}(\V(\lb),T)&,\quad&
\mathrm{e}_d(\lb+\bs{\nu},T)&=&\mathrm{e}_d(\lb,T)\cdot
\mathrm{e}_d(\bs{\nu},T)\;.
\end{array}
\end{equation*}
Furthermore, if $\B=\O_L$ is the ring of integers of some finite
extension $L/K$, and if, for some $r\geq 1$, there exists a
Frobenius $\sigma$ on $\O_L$ lifting the $p^r$-th power map
$x\mapsto x^{p^r}$ of $k_L$, and satisfying $\sigma(\pi_j)=\pi_j$,
$\forall\; 0\leq j\leq m$, then we have
\begin{equation}
 \mathrm{e}_d^\sigma(\lb,T) = \mathrm{e}_d(\sigma(\lb),T)\; ,
\end{equation}
where
$\sigma(\lambda_0,\ldots,\lambda_m)=(\sigma(\lambda_0),\ldots,\sigma(\lambda_m))$
and, for all $f(T)=\sum a_iT^i$, we set $f^\sigma(T)=\sum
\sigma(a_i)T^i$ (cf. \ref{varphi^*}).
\end{proposition}
\emph{Proof : } All the assertions are easily verified on the
phantom components. $\Box$

\subsubsection{\textbf{Study of the differential module attached to a
$\bs{\pi}$-exponential}} \label{Study of the module attached to a
pi-exponential}

We maintain the notations of Section \ref{exponentials}. As usual
$d=np^m>0$, with $(n,p)=1$. In this subsection $H/K$ is an
algebraic extension (not necessary complete) and
\begin{equation}
H_m:=H(\pi_m)\;.
\end{equation}
\begin{remark}
The Witt vectors we are considering have a finite number of
entries. Hence the exponential $\mathrm{e}_d(\lb,T)$ has its
coefficients in a \emph{finite} (thus complete) extension of $K$.
This will solve all problems concerning convergence.
\end{remark}
\begin{definition}\label{eq: L_d}
Let $\lb=(\lambda_0,\ldots,\lambda_m)\in\W_m(\O_H)$, and let
$\ph{\phi_0,\ldots,\phi_m}\in \O_H^{m+1}$ be its phantom vector.
We define
\begin{eqnarray*}
\L_d(\lb)&:=&\d-\partial_{T,\log}(\mathrm{e}_{d}(\lb,T^{-1}))=\d+
n\cdot\bigl(\sum_{j=0}^m\pi_{m-j}\cdot\phi_j \cdot
T^{-np^j}\bigr)\;.
\end{eqnarray*}
We denote by $\Mt_d(\lb)$ the differential module over $\R_{H_m}$
defined by $\L_d(\lb)$.\index{Mdlambda@$\Mt_d(\lb)=$module defined
by $\L_d(\lb)$}
\end{definition}
\begin{lemma}
$\L_d(\lb)$ is solvable at $\rho=1$, and hence
$\Mt_d(\lb)\in\mathrm{Pic}^{\mathrm{sol}}(\R_{H_m})$.
\end{lemma}
\emph{Proof : } The Taylor solution at $+\infty$ of $\L_d(\lb)$ is
$\mathrm{e}_d(\lb,T^{-1})\!\in\! 1+\pi_mT^{-1}\O_{H_m}[[T^{-1}]]$,
which has bounded coefficients and so converges for $|T|>1$. By
transfer (cf. \ref{transfer}), $Ray(\L_d(\lb),\rho)=\rho$, for all
$\rho>1$. By continuity of the radius, $Ray(\L_d(\lb),1)=1$.$\Box$
\begin{proposition}
The map $\lb\mapsto\mathrm{e}_d(\lb,T^{-1})$ defines a group
morphism
\begin{equation}
\W_m(\O_{H})\longrightarrow 1+\pi_m T^{-1}\O_{H_m}[[T^{-1}]]\;.
\end{equation}
More precisely, for all $\lb,\bs{\nu}\in\W_m(\O_H)$, one has:
\begin{equation}\label{rules 2}
\begin{array}{rclcrcl}
\varphi_p^*(\Mt_{d}(\lb))&=&\Mt_{pd}(\V(\lb))&,&
\Mt_d(\lb+\bs{\nu})&=&\Mt_d(\lb) \otimes \Mt_d(\bs{\nu})\;,
\end{array}
\end{equation}
where $\varphi_p(f(T))=f(T^p)$ (cf. \ref{varphi_p^*}). Moreover,
if there exists an absolute Frobenius $\sigma$ on $H_m$ (cf.
\ref{what is an absolute Frobenius}) such that
$\pi_j^\sigma=\pi_j$, for all $0\leq j\leq m$, then
\begin{equation*}
\begin{array}{rclcrcl}
 \varphi_{\sigma}(\mathrm{e}_d(\lb,T)) & = & \mathrm{e}_d(\varphi_{\sigma}(\lb),T^p)
 &,\quad& \varphi_{\sigma}^*(\Mt_d(\lb))&=&\Mt_{pd}(\V(\sigma(\lb)))\;,
\end{array}
\end{equation*}
where $\varphi_\sigma(f(T))=f^{\sigma}(T^p)$, (cf.
\ref{varphi_sigma^*}), and $\sigma(\lambda_0,\ldots,\lambda_m)=
(\sigma(\lambda_0),\ldots,\sigma(\lambda_m))$.
\end{proposition}
\emph{Proof : } The first part is a direct consequence of
\ref{rules}. The last assertion is a consequence of
\ref{automorphism} and is placed here for expository reasons.
Observe that, in the sequel, we do not suppose the existence of
$\sigma$ on $H$. Indeed, our ``Frobenius structure Theorem'' does
not require the existence of $\varphi$ (cf. \ref{particular
frobenius}). $\Box$
\begin{remark}
In particular $\Mt_d$ defines a morphism of groups
\begin{equation}
\Mt_d:\W_m(\O_H)\longrightarrow\mathrm{Pic}^{\mathrm{sol}}(\R_{H_m})\;.
\end{equation}
\end{remark}
\begin{theorem} \label{reduction to k}
Let $\bs{\lb}:=(\lambda_0,\ldots,\lambda_m) \in \W_m(\O_{H})$ and
let $\ph{\phi_0,\ldots,\phi_m}\in \O_{H}^{m+1}$ be its phantom
vector. The following assertions are equivalent:
\begin{enumerate}
 \item $\Mt_d(\lb)$ is trivial (i.e. isomorphic to $\R_{H_m}$);
 \item The exponential $\mathrm{e}_d(\lb,T)$ is over-convergent
 (i.e. convergent in some disk $|T|<1+\varepsilon$, with
 $\varepsilon>0$);
 \item $|\lambda_0|,\ldots,|\lambda_m|<1$.
\end{enumerate}

Moreover, if $|\lambda_0|,\ldots,|\lambda_{r-1}|<1$ and
$|\lambda_r|=1$, $r\leq m$, then we have (cf. \ref{ell})
\begin{equation}
\mathrm{Irr}(\Mt_d(\lb))=n\cdot p^{\ell(\lb)}=d/p^r\;.
\end{equation}
\end{theorem}
\emph{Proof : } The equivalence $(1)\Leftrightarrow(2)$ is
evident. By \eqref{et_d is a product} and \ref{E_m converge
exactly for T<1 iff P=LT}, condition $(3)$ implies that
$\mathrm{e}_{d}(\lb,T^{-1})\in\R_{H_m}$, so $\Mt_d(\lb)$ is
trivial. The converse follows from the last assertion below. Let
then $|\lambda_0|,\ldots,|\lambda_{r-1}|<1$, and $|\lambda_r|=1$,
$r\leq m$. Clearly $|\phi_m(\lb)|=1$ if and only if
$|\lambda_0|=1$ (cf. \eqref{phantom components}). Then, if $r=0$,
we can apply \ref{graphic}, and hence
$\mathrm{Irr}(M)=\mathrm{Irr}_F(M)=d$. Let now $0<r\leq m$, then
$E_{m-j}(\lambda_jT^{-p^j})$ belongs to
$\mathcal{R}_{H_m}^{\times}$, for all $j=0,\ldots,r-1$. Then we
change basis by the function
$f(T):=\prod_{j=0}^{r-1}E_{m-j}(\lambda_jT^{-p^j})^{-1}\in
\R_{H_m}^\times$. By \ref{rules}, the new solution is
\begin{equation}
f(T)\cdot\mathrm{e}_d(\lb,T)=
\mathrm{e}_d((0,\ldots,0,\lambda_r,\ldots,\lambda_m),T)=
\mathrm{e}_{d/p^r}((\lambda_r,\ldots,\lambda_m),T^{p^r}).
\end{equation}
In other words, we have
$\Mt_d(\lambda_0,\ldots,\lambda_m)\stackrel{\sim}{\to}
\varphi_p^{*}(\Mt_{d/p^r}(\lambda_r,\ldots,\lambda_m))$ (cf.
\ref{varphi^*} and \eqref{rules 2}). By \ref{irreg of frob}, the
Theorem is proved by induction, since $|\lambda_r|=1$.$\Box$

\begin{remark} \label{the morphism Mt_d}
In particular, $\Mt_d$ passes to the quotient $\W_m(k_H)$, and
induces an injective additive map called $\M_d$:
\begin{equation}
\xymatrix{
\W_m(\O_{H})\ar[r]^-{\Mt_d}\ar@{->>}[d] &\mathrm{Pic}^{\mathrm{sol}}(\R_{H_m})\\
\W_m(k_H)\ar@{^{(}-->}[ur]_{\M_d} &  }
\end{equation}
\end{remark}

\begin{corollary} \label{explicit corollary link with Lubin-Tate}
Consider the morphism of groups
\begin{equation}
\mathbb{Z}_p[[X]]\xrightarrow[]{[-]}
\W_m(\mathbb{Z}_p[\pi_m])\subset\W_m(\O_{H_m})
\xrightarrow[]{\Mt_d}\mathrm{Pic}^{\mathrm{sol}}(\R_{H_m}).
\end{equation}
Let $h(X):=\sum_{i\geq 0} a_iX^i\in\mathbb{Z}_p[[X]]$ be such that
$|a_0|=|p|^r$ ($v_p(a_0)=r$). Then $\Mt_d([h(\pi_m)])$ has
irregularity $d/p^r$, and is trivial if and only if $r\geq m+1$. In
other words, the kernel of the composite map is the ideal
$p^{m+1}\mathbb{Z}_p[[X]]+X\mathbb{Z}_p[[X]]$.
\end{corollary}
\emph{Proof : } Combine \ref{base lemma!} and \ref{reduction to
k}.$\Box$
\subsubsection{\textbf{Dependence on the Lubin-Tate group and on
$\bs{\pi}$}} We maintain the notations of the previous sections.
As usual $d=np^m>0$, $(n,p)=1$.
\begin{theorem}[Dependence on the choice of $\bs{\pi}$]
\label{dependence from pi_m} Let $\bs{\pi}=(\pi_j)_{j\geq 0}$,
$\bs{\pi}'=(\pi_j')_{j\geq 0}$ be two generators of
$T(\mathfrak{G}_P)$. Denote by $\M'_d(-)$, $E'_j(T)$ and
$\mathrm{e}'_d(-,T)$ the constructions attached to $\bs{\pi'}$.
Then $\M_d(1,0,\ldots,0)$ and $\M'_d(1,0,\ldots,0)$ are isomorphic
over $\R_{H_m}$ if and only if $\pi_m=\pi'_m$. Moreover, in this
case, $\M_d(\lb)$ and $\M'_d(\lb)$ are isomorphic over
$\mathcal{R}_{H_m}$, for all $\lb\in\W_m(k_H)$.
\end{theorem}
\emph{Proof : } The solution at $\infty$ of $\R_{H_m}$ is
$\mathrm{e}_{d}((1,0,\ldots,0),T^{-1})=E_m(T^{-n})$. We shall show
that $E_m(T^{-n})/E'_m(T^{-n})\in\R^{\times}$, that is
$E_m(T^{-1})/E'_m(T^{-1})\in\R^{\times}$, if and only if
$\pi_m=\pi'_m$. We have
\begin{equation}
E_m(T^{-1})/E'_m(T^{-1})=
\exp\Bigl(\sum_{j=0}^{m}\pi_{m-j}(1-\frac{\pi'_{m-j}}{\pi_{m-j}})
\frac{T^{-p^j}}{p^j}\Bigr)\;.
\end{equation}
There exists $\gamma\in
\mathrm{Gal}(\mathbb{Q}_p(\Lambda_P)/\mathbb{Q}_p)$ such that
$\pi'_j=\gamma(\pi_j)$, for all $j\geq 0$ and, by the Lubin-Tate
Theorem \ref{reciprocity and filtrations},
$\gamma(\pi_j)=[u_\gamma]_{P}(\pi_j)$,
$u_\gamma\in\mathbb{Z}^{\times}$. We set \footnote{Note that the
symbol $[-]_P$ was defined in \ref{[a]_{P,widetilde{P}}} and is
different from $[-]$ defined at \ref{[h(pi_m)]}.}
\begin{equation}
h_{\gamma}(X):=1-[u_\gamma]_P(X)/X\;,
\end{equation}
in order to have
\begin{equation}
E_m(T^{-1})/E'_m(T^{-1})=\mathrm{e}_{d}([h_{\gamma}(\pi_m)],T^{-1}).
\end{equation}
Indeed, by construction (cf. corollary \ref{base lemma!} and
definition \eqref{definition of et_d(lb,T)}), we have
$\phi_j([h_\gamma(\pi_m)])=1-\pi'_{m-j}/\pi_{m-j}$. Since
$h_\gamma(0) = 1-u_{\gamma}$, hence, by the Reduction Theorem
\ref{reduction to k} and lemma \ref{explicit corollary link with
Lubin-Tate}, the series $E_m(T^{-1})/E'_m(T^{-1})$ lies in
$\mathcal{R}^{\times}$ if and only if $|1-u_\gamma|\leq
|p|^{m+1}$, i.e. $u_\gamma\in 1+p^{m+1}\mathbb{Z}_p$. Then, again
by the reciprocity law \ref{reciprocity and filtrations}, the
automorphism $\gamma$ is the identity on
$\mathbb{Q}_p(\Lambda_{P,m+1})=\mathbb{Q}_{p}(\pi_m)$. Hence
$\pi_m=\pi_m'$.$\Box$\\

We recall that two Lubin-Tate groups are isomorphic (as formal
groups over $\mathbb{Z}_p$) if and only if they are relative to
the same uniformizer $\mathrm{w}$ (cf. Theorem \ref{moduly space
of Lubin Tate Groups}).
\begin{theorem}[Independence on the Lubin-Tate group]\label{indipendence from P}
Let $P,\widetilde{P}\in\mathfrak{F}_\mathrm{w}$ be two Lubin-Tate
series, let $\bs{\pi}=(\pi_j)_{j\geq 0}$ and
$\widetilde{\bs{\pi}}=(\pi_{\widetilde{P},j})_{j\geq 0}$ be a
generator of $\T(\mathfrak{G}_P)$ and
$\T(\mathfrak{G}_{\widetilde{P}})$ respectively. Let us denote by
$\M_d^{(\widetilde{P})}(-)$, $E^{(\widetilde{P})}_m(T)$,
$e_d^{(\widetilde{P})}(-,T)$ the constructions attached to
$\widetilde{\bs{\pi}}$, and denote in the usual way the
constructions attached to $\bs{\pi}$. If
$\pi_{\widetilde{P},m}=[1]_{P,\widetilde{P}}(\pi_{P,m})$, then
$\M_d(\lb)\stackrel{\sim}{\to}\M_d^{(\widetilde{P})}(\lb)$ over
$\mathcal{R}_{H_m}$, for all $\lb\in\W_m(k_H)$.
\end{theorem}
\emph{Proof : } Let $\lb\in\W_m(k_H)$, and let
$\tilde{\lb}\in\W_m(\O_{H})$ be a lifting of $\lb$. We shall show
that $e_d(\tilde{\lb},T)/e_d^{(\widetilde{P})}(\tilde{\lb},T)$
belongs to $\R_{H_m}$. By equation \eqref{et_d is a product}, we
reduce to showing that
$E_{m-j}(T^{-1})/E^{(\widetilde{P})}_{m-j}(T^{-1})\in\R_{H_{m-j}}$,
for all $0\leq j\leq m$. Since
$\pi_{\widetilde{P},m}=[1]_{P,\widetilde{P}}(\pi_{P,m})$, then
$\pi_{\widetilde{P},j}=[1]_{P,\widetilde{P}}(\pi_{P,j})$, for all
$0\leq j\leq m$. We have
\begin{equation}
E_m(T^{-1})/E^{(\widetilde{P})}_m(T^{-1})=
\exp\bigl(\sum_{j=0}^{m}\pi_{m-j}(1-\frac{\pi_{\widetilde{P},m-j}}{\pi_{P,m-j}})
\frac{T^{-p^j}}{p^j}\bigr).
\end{equation}
Let us set, as usual,
$h_{P,\widetilde{P}}(X):=1-[1]_{P,\widetilde{P}}(X)/X$, in order
to have (cf. \ref{base lemma!} and definition \eqref{definition of
et_d(lb,T)}) $E_m(T^{-1})/E^{(\widetilde{P})}_m(T^{-1})=
\mathrm{e}_d([h_{P,\widetilde{P}}(\pi_m)],T^{-1})$. Since
$[1]_{P,\widetilde{P}}(X)\equiv X\mod X^2\mathbb{Z}_p[[X]]$, we
have $h_{P,\widetilde{P}}(X)\in X\cdot\mathbb{Z}_p[[X]]$, and, by
the Reduction Theorem \ref{reduction to k} and Lemma \ref{explicit
corollary link with Lubin-Tate}, this exponential lies in
$\R_{H_m}$. By the way, its inverse lies also in $\R_{H_m}$, so
$\M_d(\lb)\stackrel{\sim}{\to}\M_d^{(\tilde{P})}(\lb)$, over
$\R_{H_m}$.$\Box$

\begin{remark}
If $\mathrm{w}=p$, and if $P$ is given, then, by \ref{equiv
classe-} and \ref{relations on equivalence classe}, the
isomorphism class of $\M_d(\lb)$ is determined by the choice of a
sequence $\{\xi_j\}_{j\geq 0}$ of $p^{j+1}$-th roots of $1$ such
that $\xi_m^{p^j}=\xi_{m-j}$.
\end{remark}

\begin{corollary}\label{automorphism}
Let $\gamma:H(\Lambda_P)\to H(\Lambda_P)$ be a continuous
endomorphism of fields. Then
$\gamma(E_m(T^{-1}))/E_m(T^{-1})\in\mathcal{R}_{H_m}$ if and only
if $\gamma$ is the identity on $\mathbb{Q}_p(\pi_m)$, and in this
case, for all $\lb\in\W_m(\O_{H(\Lambda_{P})})$, we have
\begin{equation}
\mathrm{e}_d^{\gamma}(\lb,T)=\mathrm{e}_d(\gamma(\lb),T)\;,
\end{equation}
where, for all $f(T)=\sum a_iT^i$, we set $f^\gamma(T):=\sum
\gamma(a_i)T^i$.
\end{corollary}
\emph{Proof : } The proof follows the same lines as the proof of
\ref{dependence from pi_m}.$\Box$

\subsubsection{\textbf{Frobenius structure for $\bs{\pi}$-exponentials}}

\begin{theorem}\label{overconvergence of theta}\label{frobenius structure}
Let $r\geq 0$ and let $\bar{\lb}\in\W_m(k_H)$. Let
$\lb\in\W_m(\O_H)$ be a lifting of $\bar{\lb}$, and let
$\lb^{(\Fb)}\in\W_m(\O_H)$ be an arbitrary lifting of
$\Fb(\bar{\lb})\in\W_m(k_H)$. The following statements are
equivalent:
\begin{enumerate}
 \item The power series
 $\mathrm{e}_d(\lb^{(\Fb)},T^{p})/\mathrm{e}_d(\lb,T)$
 is over-convergent, for all choices of $\lb$, $\bar{\lb}$ and $\lb^{(\Fb)}$;
 \item The modules $\M_d(\bar{\lb})$ and $\M_{pd}(\V\Fb(\bar{\lb}))$ are
 isomorphic over $\R_{H_m}$, for all $\bar{\lb}\in\W_m(k_H)$;
 \item The power series $E_m(T^p)/E_m(T)$ is over-convergent;
 \item We have the inequality $|\mathrm{w}-p|\leq |p|^{m+2}$.
\end{enumerate}
\end{theorem}
\emph{Proof : } $(1)\Leftrightarrow (2)$ and $(1)\Rightarrow(3)$
are evident. Let us show $(3)\Leftrightarrow(4)$. Write
\begin{eqnarray*}
E_{m}(T^p)/E_m(T)&=&\exp\Bigl(\bigl(\sum_{j=0}^{m}\pi_{m-j}
\frac{T^{p^{j+1}}}{p^j}\bigr)-\bigl(\sum_{j=0}^{m}\pi_{m-j}
\frac{T^{p^j}}{p^j}\bigr) \Bigr)\\
&=&\exp(-p\pi_{m+1}T)\cdot\exp\Bigl(\sum_{j=0}^{m+1}\pi_{m+1-j}
\bigl(p-\frac{\pi_{m-j}}{\pi_{m-j+1}}\bigr)\frac{T^{p^j}}{p^j}\Bigr)\;,
\end{eqnarray*}
where $\pi_{-1}:=P(\pi_0)=0$.  Let $h_{\mathrm{Frob}}(X) := p -
P(X)/X$, in order to have (cf. \ref{base lemma!} and definition
\eqref{definition of et_d(lb,T)})
\begin{equation}\label{E_m(T^p)/E_m(T) = exp(-ppiT) exp(nu_Frob,T)}
E_{m}(T^p)/E_m(T)=\exp(-p\pi_{m+1}T)\cdot\mathrm{e}_{p^{m+1}}([h_{\mathrm{Frob}}(\pi_{m+1})],T).
\end{equation}
Since the function $\exp(-p\pi_{m+1}T)$ is over-convergent, the
quotient $E_{m}(T^p)/E_m(T)$ is over-convergent if and only if
$\mathrm{e}_{p^{m+1}}([h_{\mathrm{Frob}}(\pi_{m+1})],T)$ is. The
constant term of $h_{\mathrm{Frob}}(X)$ is $p-\mathrm{w}$. Hence,
as usual, by the Reduction Theorem \ref{reduction to k} and Lemma
\ref{base lemma!}, $E_m(T^p)/E_m(T)$ is over-convergent if and
only if $|p-\mathrm{w}|\leq|p|^{m+2}$. Let us now show
$(3)\Rightarrow(1)$. Since $(3)$ and $(4)$ are equivalent, we see
that $E_{j}(T^p)/E_{j}(T)$ is over-convergent, for all
$j=0,\ldots,m$. Let $\lb=(\lambda_0,\ldots,\lambda_m)$ and
$\lb^{(\Fb)}:=(\lambda_0^{(\Fb)},\ldots,\lambda_m^{(\Fb)})$. We
can suppose $\lambda_j^{(\Fb)}=\lambda_j^{p}$,
$\forall\;j=0,\ldots,m$. Indeed, the Witt vector
$\bs{\eta}:=\lb^{(\Fb)}-(\lambda_0^p,\ldots,\lambda_m^p)=
(\eta_0,\ldots,\eta_m)\in\W_m(\O_H)$ satisfies $|\eta_j|<1$,
$\forall\;j=0,\ldots,m$. Hence
%\begin{equation}
$\mathrm{e}_d(\lb^{(\Fb)},T^p)= \mathrm{e}_d(\bs{\eta},T^p)\cdot
\mathrm{e}_d((\lambda_0^p,\ldots,\lambda_m^p),T^p)$,
%\end{equation}
and the function $\mathrm{e}_d(\bs{\eta},T^p)$ is over-convergent
by the Reduction Theorem \ref{reduction to k}. Now, by equation
\eqref{et_d is a product}, we have
\begin{equation}
\frac{\mathrm{e}_d((\lambda_0^p,\ldots,\lambda_m^p),T^p)}{\mathrm{e}_d(\lb,T)}=
\prod_{j=0}^{m}\frac{E_{m-j}(\lambda_{j}^{p}
T^{np^{j+1}})}{E_{m-j}(\lambda_j T^{np^j})}.
\end{equation}
Since $E_{m-j}(T^p)/E_{m-j}(T)$ is over-convergent, for all
$j=0,\ldots,m$, all factors $E_{m-j}(\lambda_{j}^{p}
T^{np^{j+1}})/E_{m-j}(\lambda_j T^{np^j})$ are over-convergent.
$\Box$

\begin{remark}\label{we non need the existece of phi}
In this Theorem we do not need the existence of an absolute
Frobenius on $H$. This is due to the fact that the isomorphism
class of $\M_d(\lb)$ depends only on the reduction
$\bar{\lb}\in\W_m(k_H)$, and $k_H$ is endowed naturally with the
Frobenius given by the $p$-th power map.
\end{remark}
\begin{remark}\label{particular frobenius}
We will generalize this Theorem for all rank one differential
equations (cf. Theorem \ref{Introductive foundamental Theorem}).
Let us show how to recover, from \ref{frobenius structure}, the
Frobenius structure Theorem in the usual sense. Let $\lb \in
\W_m(\O_{H})$ be a lift of $\bar{\lb}\in\W_m(k_H)$. Suppose that
$\mathrm{w}=p$, in order to apply Theorem \ref{frobenius
structure}. Suppose that $\sigma:H_\infty \to H_\infty$ is an
absolute Frobenius (cf. \ref{what is an absolute Frobenius}) such
that $\pi_j^\sigma=\pi_j$, for all $j\geq 0$, and such that
$\sigma(H)\subseteq H$. By Corollary \ref{automorphism}, we have
$\varphi_\sigma(\mathrm{e}_d(\lb,T))=\mathrm{e}_d(\lb^\sigma,T^p)$,
and hence
$\varphi_{\sigma}^*(\Mt_{d}(\lb))=\Mt_{pd}(\V(\lb^\sigma))$. By
\ref{reduction to k}, the isomorphism class of
$\Mt_{pd}(\V(\lb^\sigma))$ depends only on the reduction
$\V(\overline{\lb^\sigma})=\V\Fb(\bar{\lb})\in\W_{m+1}(k_{H_\infty})$,
so $\Mt_{pd}(\V(\lb^\sigma))$ is isomorphic to
$\M_{pd}(\V\Fb(\lb))$ over $\R_{H_\infty}$. Then Theorem
\ref{frobenius structure} gives us the usual Frobenius structure.
Indeed,
 $$\varphi_{\sigma}^*(\Mt_d(\lb))
 \xrightarrow[\textrm{Cor.\ref{automorphism}}]{\sim}
 \Mt_{p\cdot d}(\V(\lb^\sigma))
 \xrightarrow[\textrm{Th.\ref{reduction to k}}]{\sim}
 \M_{p\cdot d}(\V\Fb(\bar{\lb}))
 \xrightarrow[\textrm{Th.\ref{frobenius structure}}]{\sim}
 \Mt_d(\lb)\;.$$
\end{remark}
\begin{remark}
Let $\varphi_p^*$ be the $p$-th ramification map (cf.
\ref{varphi_p^*}), and let $\lambda\in\O_H$. We observe that we
can not have an isomorphism
$\M_1(\lambda)\stackrel{\sim}{\to}(\varphi_p^*)^{h}(\M_1(\lambda))$,
for all $\lambda$ and all $h\geq 1$. In other words there exists
module which do not admits a ``ramification structure''. For
example, suppose that $\lambda$ is such that
$\bar{\lambda}^{p^r}\neq \bar{\lambda}$ in $k_H$, for all $r\geq
0$ (i.e. $\bar{\lambda}\notin\mathbb{F}_p^{\mathrm{alg}})$. Then
$\exp(\pi_0\lambda T^p)/\exp(\pi_0\lambda T)$ is not
over-convergent. Indeed, for all liftings
$\lambda^{(\Fb^r)}\in\O_H$ of $\bar{\lambda}^{p^r}$ we have
$|\lambda^{(\Fb^r)}-\lambda|=1$, then
$$\frac{\exp(\pi_0\lambda T^p)}{\exp(\pi_0\lambda
T)}=\frac{\mathrm{e}_1(\lambda,T^p)}{\mathrm{e}_1(\lambda,T)}=
\frac{\mathrm{e}_1(\lambda^{(\Fb^r)},T^p)}{\mathrm{e}_1(\lambda,T)}
\cdot\mathrm{e}_1(\lambda-\lambda^{(\Fb^r)},T^p)\;,$$ and while
$\mathrm{e}_1(\lambda^{(\Fb^r)},T^p)/\mathrm{e}_1(\lambda,T)$ is
over-convergent, the function
$\mathrm{e}_1(\lambda-\lambda^{(\Fb^r)},T^p)$ is not
over-convergent, since the reduction of
$\lambda^{(\Fb^r)}-\lambda$ is not $0$ in $k_H$ (cf.
\ref{reduction to k}).
\end{remark}

\subsection{Deformation of the Artin-Schreier complex
into the Kummer complex}

\label{Dw's theta function}
\begin{hypothesis}\label{hypothesis uni=p}
From now on, we will suppose $\mathrm{w}=p$ in order to have the
Theorem \ref{overconvergence of theta}. Then
$\mathfrak{G}_P\stackrel{\sim}{\to}\widehat{\mathbb{G}}_m$, the
formal multiplicative group (cf. \ref{moduly space of Lubin Tate
Groups}). We fix moreover a generator $\bs{\pi}=(\pi_j)_{j\geq 0}$
of $\T(\mathfrak{G}_P)$.
\end{hypothesis}
In this section $L$ will be a complete valued field, containing
$(\mathbb{Q}_p,|.|)$, and endowed with an absolute Frobenius
$\varphi:\O_L\to \O_L$ (i.e. a lifting of the map $x\mapsto x^p$
of $k_L$).

We set as usual $L_m:=L(\pi_m)$ and $L_\infty:=\cup_mL_m$. We
denote by $k_m$ (resp. $k_{\infty}$) the residue field of $L_m$
(resp. $L_\infty$). We fix an algebraic closure $L^{\mathrm{alg}}$
of $L$, then $k_L^{\mathrm{alg}}:=k_{L^{\mathrm{alg}}}$ is an
algebraic closure of $k_L$. Let $k_{L}^{\mathrm{sep}}$ (resp.
$k_{m}^{\mathrm{sep}}$, $k_{\infty}^{\mathrm{sep}}$) be the
separable closure of $k_L$ (resp. $k_{m}$, $k_{\infty}$) in
$k_{L}^{\mathrm{alg}}$ (we recall that \emph{$k_L$ is not supposed
to be perfect}). We denote by $\widehat{L}^{\mathrm{unr}}$ (resp.
$\widehat{L}_m^{\mathrm{unr}}$) the completion of the maximal
unramified extension of $L$ (resp. $L_m$) in $L^{\mathrm{alg}}$.
We set $\mathrm{G}_{k_L}:=\mathrm{Gal}(k_L^{\mathrm{sep}}/k_L)$,
$\mathrm{G}_{k_{m}}:= \mathrm{Gal}(k_{m}^{\mathrm{sep}}/k_{m})$,
$\mathrm{G}_L:=\mathrm{Gal}(L^{\mathrm{alg}}/L)$, and
$\mathrm{G}_{L_m}:=\mathrm{Gal}(L_m^{\mathrm{alg}}/L_m)$.

\begin{remark}\label{definition of k_L^0}
Let $k_{L}^{0}=k_L^{\mathrm{sep}}\cap k_m$ be the separable
closure of $k_L$ in $k_{m}$ and let
$L^0:=\W(k_{L}^0)\otimes_{\W(k_L)}L=\widehat{L}^{\mathrm{unr}}\cap
L_m$. The extension $k_{m}/k_{L}^{0}$ is purely inseparable (i.e.
for all $x\in k_{m}$ there exists $r\geq 0$ such that $x^{p^r}\in
k_{L}^0$), so $\mathrm{Gal}(k_{m}/k_{L}^{0})=1$, and we have a
canonical identification
$\G_{k_{m}}:=\mathrm{Gal}(k_m^{\mathrm{sep}}/k_m)
\xrightarrow[]{\sim}\mathrm{Gal}(k_{L}^{\mathrm{sep}}/k_{L}^{0})$.
Hence $\G_{k_{m}}$ is naturally contained in $\G_{k_L}$:
\begin{equation}
\begin{array}{ccccc}
  && L_m &\subseteq& \widehat{L}_m^{\mathrm{unr}} \\
  && \turnup{\subseteq}&&\turnup{\subseteq}\\
L &\subseteq& L^0 &\subseteq& \widehat{L}^{\mathrm{unr}} \\
\end{array}
\qquad\qquad
\begin{array}{ccccc}
  && k_{m} &\subseteq& k_{m}^{\mathrm{sep}}\;\phantom{.} \\
  && \turnup{\subseteq}&&\turnup{\subseteq}\;\phantom{.}\\
k_L &\subseteq& k_L^0 &\subseteq& k_L^{\mathrm{sep}}\;. \\
\end{array}
\end{equation}
All these extensions are normal. We will identify $\G_{k_{m}}$
with $\mathrm{Gal}(\widehat{L}_m^{\mathrm{unr}}/L_m)$, and
$G_{k_L}$ with $\mathrm{Gal}(\widehat{L}^{\mathrm{unr}}/L)$. In
this way $\G_{k_{m}}$ acts naturally on
$\widehat{L}^{\mathrm{unr}}$.
\end{remark}

\begin{remark}\label{trtrtrtrrrgft}
The absolute Frobenius $\varphi$ extends uniquely to all
unramified extensions of $L$, and hence it commutes with the
action of $\G_{k_L}$. It extends also (not uniquely) to an
absolute Frobenius $\widetilde{\varphi}$ of $L_m$. Indeed, since
$\varphi$ extends uniquely to $L^0$, then to prove the existence
of $\widetilde{\varphi}$ one can assume that $L=L^0$, and hence
$k_{m}$ is a purely inseparable extension of $k_L=k_{L}^0$. Since
the map $x\mapsto x^p$ of $k_L$ extends uniquely to $k_m$, then
every field morphism $\widetilde{\varphi}:L_m\to L_m$ extending
$\varphi$ is an absolute Frobenius of $L_m$. %
Such a $\widetilde{\varphi}$ exists since, by \cite[$\S
6,\textrm{n}^{\textrm{o}}1$, prop.$1$]{Bou-Alg-4-5}, $\varphi$
extends to a $\mathbb{Q}_p$-linear morphism
$\widetilde{\varphi}:L^{\mathrm{alg}}\to L^{\mathrm{alg}}$,
inducing an automorphism of $\mathbb{Q}_p(\pi_m)$.

In general there is no absolute Frobenius on $L_m$ satisfying
$\varphi(\pi_m)=\pi_m$. Indeed if $L$ is totally ramified over
$\mathbb{Q}_p$, and if $\varphi=\mathrm{Id}_L$, then the unique
extension of $\varphi$ to $L_m$, fixing $\pi_m$, is the identity.
On the other hand $L_m/\mathbb{Q}_p$ is not always totally
ramified, hence the identity of $L_m$ is not always an absolute
Frobenius\footnote{Indeed let $p=3$, $m=0$, and
$L:=\mathbb{Q}_p(\pi_{P,0})$, where $P(X)$ is the Lubin-Tate
series $P(X)=-3X+X^3$. %, with $\mathrm{w}=-p$.
If $\xi_1^3=1$ is a primitive root of unity, then
$L_0=\mathbb{Q}_p(\pi_{P,0},\xi_1)$ is not totally ramified since
the element $x:=\widetilde{\pi}/\pi_{P,0}$, where
$\widetilde{\pi}=(\xi_1-1)$, verifies $|x|=1$ and $x^{6}=-1$,
indeed
$x^9=(\frac{\widetilde{\pi}}{\pi_{P,0}})^9=(\frac{-3\widetilde{\pi}-3\widetilde{\pi}^2}{3\pi_{P,0}})^3=
-\xi_1^3x^3=-x^3$. But there is no element $\overline{x}$ in
$\mathbb{F}_3$ verifying $\overline{x}^6=-\overline{1}$.}.

In the sequel of the paper \emph{we will never use such a
$\widetilde{\varphi}$, hence we do not fix it.}

On the other hand, we need the existence of $\varphi$ because the
functor of Witt vectors of finite length $\W_m(-)$ is not
canonically endowed with an additive functorial Frobenius morphism
(see remark \ref{F_p instead of varphi} to improve this
situation).
\end{remark}

\begin{definition}
For $\lb :=(\lambda_0,\ldots,\lambda_m)\in\W_m(\O_L)$, we set
\begin{equation}\label{definition of theta_d}\index{theta_d@$\theta_d^{(\varphi)}(\lb,T)$}
\theta_d^{(\varphi)}(\lb,T):=
\frac{\mathrm{e}_{d}(\varphi(\lb),T^{p})}{\mathrm{e}_{d}(\lb,T)}.
\end{equation}
To simplify the notations, we will write $\theta_d(\lb,T)$ if no
confusion is possible.
\end{definition}
\begin{example}
Let $d=1$ and $P(X) = pX + X^p$ (cf. \ref{(xi_m-1)}). Then $\pi_0$
is the ``$\pi$ of Dwork'', and $\theta_1(1,T)=\exp(\pi_0(T^p-T))$
is the usual Dwork splitting function. While in general, if
$\lambda\in\O_L$, we have
 $\theta_1(\lambda,T) = \exp(\pi_0(\varphi(\lambda)T^p-\lambda
 T))$.
\end{example}
The following Theorem shows that the over-convergent function
$\lb\mapsto\theta(\lb,1)$ is a splitting function in a generalized
sense with respect to Dwork (cf. \cite[$\S4$, $a)$,p.55]{Dw}). In
a paper in preparation we shall analyse such functions in detail.

\begin{definition}
Set $\O_L^{\varphi=1}:=\{\lambda\in
\O_L\;|\;\varphi(\lambda)=\lambda\}$ and
$\overline{\O_L^{\varphi=1}}:=\O_L^{\varphi=1}/(\O_L^{\varphi=1}\cap\mathfrak{p}_L)$.
We see that $\overline{\O_L^{\varphi=1}}=\mathbb{F}_p$.
\end{definition}

\begin{theorem}\label{theta pass to the quotient}
Let $a^{p}=a\in\O_L$, and let $\lb\in\W_m(\O_{L}^{\varphi=1})$.
Then $\theta_d^{(\varphi)}(\lb,a)$ is a $p^{m+1}$-th root of $1$.
Moreover the group morphism
$$\theta_d^{(\varphi)}(-,a) :
\W_m(\O_{L}^{\varphi=1})\longrightarrow
\bs{\mu}_{p^{m+1}}\subset\mathbb{Z}_p[\pi_m]$$ factorizes on
$\W_m(\overline{\O_{L}^{\varphi=1}})=\W_m(\mathbb{F}_{p})=\mathbb{Z}/p^{m+1}\mathbb{Z}$
and defines an isomorphism
\begin{equation}\label{1 to xi_m^-1}
\overline{\theta}_d^{(\varphi)}(-,a):
\mathbb{Z}/p^{m+1}\mathbb{Z}\xrightarrow[]{\;\;\sim\;\;}
\bs{\mu}_{p^{m+1}}.
\end{equation}
More precisely the image of $1\in\mathbb{Z}/p^{m+1}\mathbb{Z}$ is
\emph{the inverse} of the unique primitive $p^{m+1}$-th root of
$1$, say $\xi_m$, satisfying
\begin{equation}
|a^n\pi_m-(\xi_m-1)|<|a^n\pi_m|.
\end{equation}
In particular, if $a=1$, then $\xi_m$ is the $p^{m+1}$-th root of
$1$ defined in remark \ref{relations on equivalence classe}.
\end{theorem}
\emph{Proof : } Let $\lb=(\lambda_0,\ldots,\lambda_m)\in
\W_m(\O_{L}^{\varphi=1})$. Let us show that
$\theta_d(\lb,a)^{p^{m+1}}=1$. Indeed
$T\mapsto\mathrm{e}_d(\lb,T)^{p^{m+1}}$ is over-convergent (cf.
\ref{et_d(lb,1)^p^m+1 is over-convergent}), so
$\theta_d(\lb,a)^{p^{m+1}}=\mathrm{e}_{d}(\varphi(\lb),a^{p})^{p^{m+1}}/
\mathrm{e}_d(\lb,a)^{p^{m+1}}=1$, since both numerator and
denominator do make sense and are equal. If $|\lambda_j|<1$, for
all $j=0,\ldots,m$, then $T\mapsto\mathrm{e}_d(\lb,T)$ is
over-convergent (cf. reduction Theorem \ref{reduction to k}),
hence both numerator and denominator of the expression
$\mathrm{e}_{d}(\lb,a^{p})/\mathrm{e}_d(\lb,a)$ do make sense and
are equal. Let us show the last assertion. By equation
\eqref{E_m(T^p)/E_m(T) = exp(-ppiT) exp(nu_Frob,T)}, we have
\begin{eqnarray}
\theta_d((1,0\ldots,0),T) &=&
E_{m}(T^{np})/E_{m}(T^n)\\
&=&\exp(-p\pi_{m+1}T^n)\cdot
\mathrm{e}_{pd}([h_\mathrm{Frob}(\pi_{m+1})],T)\nonumber.
\end{eqnarray}
By \eqref{et_d(-,T) has values in 1+pi_m T} this series lies in
$1+\pi_{m+1}T\mathbb{Z}_p[\pi_{m+1}][[T]]$. To show that this root
is $\xi_m^{-1}$ it is sufficient to show that
$|\theta_{d}((1,0,\ldots,0),a)^{-1}-\xi_m|<|\pi_m|=|\xi_m-1|$. We
work therefore modulo the following sub group
$$C:=\{1+\sum c_iT^i\;|\;c_i\in\mathbb{Z}_p[\pi_{m+1}]\;,\;
|c_i|<|\pi_{m}|\;, \textrm{ for all }i\geq 1 \}.$$ We have
$\exp(-p\pi_{m+1}T^n)\equiv 1\mod C$. Let us consider (cf.
definition \eqref{E_m(T^p)/E_m(T) = exp(-ppiT) exp(nu_Frob,T)})
\begin{equation}
[h_\mathrm{Frob}(\pi_{m+1})]=[p-P(\pi_{m+1})/\pi_{m+1}]=(\nu_0,\ldots,\nu_{m+1}).
\end{equation}
Then $\nu_0=p-(\pi_m/\pi_{m+1})$ and, since $p=\mathrm{w}$, by
Lemma \ref{base lemma!}, we have $|\nu_j|\leq |\pi_{m+1}|$, for
all $j=0,\ldots,m+1$. By equation \eqref{et_d is a product} we
have $\mathrm{e}_{pd}([h_\mathrm{Frob}(\pi_{m+1})],T)=
\prod_{j=0}^{m+1}E_{m+1-j}(\nu_{j}T^{np^{j}})$. Moreover, we know
that (cf. equation \eqref{et_d(-,T) has values in 1+pi_m T})
\begin{equation}
E_{m+1-j}(\nu_{j}T^{np^{j}})=1+(\textrm{things of valuation }\leq
|\pi_{m+1-j}\cdot\nu_j|)\;,
\end{equation}
for all $j=0,\ldots,m+1$. Then
\begin{equation}
\theta_d((1,0,\ldots,0),T)^{-1}\equiv E_{m+1}(\nu_0T^{n})^{-1}
\mod C.
\end{equation}
Since $|\nu_0^p|=|\pi_m|^{p-1}$, it follows from \eqref{et_d(-,T)
has values in 1+pi_m T} that only the first $p-1$ terms of
$E_{m+1}(\nu_0T^{n})^{-1}$ are greater than or equal to $|\pi_m|$,
that is
\begin{equation}
E_{m+1}(\nu_0T^{n})^{-1}\equiv 1+\pi_{m+1}\nu_0T^n+\cdots+
\frac{(\pi_{m+1}\nu_0T^{n})^{p-1}}{(p-1)!}\mod C.
\end{equation}
Since $\pi_{m+1}\nu_0=p\cdot\pi_{m+1}-\pi_m$, hence
$\theta_d((1,0,\ldots,0),T)^{-1} \equiv 1 + \pi_m T^n \mod
C$.$\Box$

\begin{remark}\label{et_d(lb,1)^p^m+1 is over-convergent}
Observe that
$T\mapsto\mathrm{e}_d(\lb,T)^{p^{m+1}}=\mathrm{e}_d(p^{m+1}\lb,T)$
is over-convergent for all $\lb\in\W_m(\O_L)$, because the
reduction of $p^{m+1}\lb$ in $\W_m(k_L)$ is $0$ (cf.
\ref{reduction to k}).
\end{remark}

\begin{remark}
We recall that we do not fix an absolute Frobenius on $L_m$ (cf.
remark \ref{trtrtrtrrrgft}).
\end{remark}

\begin{theorem}\label{lifting characters diagram for L Theorem}
The following diagram is well-defined, commutative and functorial,
on the complete (or algebraic) unramified extensions of $L$
\begin{equation}\label{lifting characters diagram for L}
\index{delta@$\delta_{\mathrm{Kum}}$}
\xymatrix{ 1\ar[r]&\bs{\mu}_{p^{m+1}}\ar[r]&(L_m)^{\times}
\ar[r]^{f\mapsto f^{p^{m+1}}}&(L_m)^{\times}
\ar[r]^-{\delta_{\mathrm{Kum}}}&
\mathrm{H}^1(\G_{L_m},\bs{\mu}_{p^{m+1}})\ar[r]&1\\
&\W_{m}(\O_{L}^{\varphi=1})\ar@{^{(}->}[r]\ar@{->>}[d]\ar[u]&
\W_{m}(\O_{L})\ar[r]_{\varphi-1}
\ar@{->>}[d]\ar[u]_{\theta_{p^{m}}(-,1)}&
\W_{m}(\O_{L})\ar@{->>}[d]
\ar[u]_{\mathrm{e}_{p^{m}}(-,1)^{p^{m+1}}}&&\\
0\ar[r]& \mathbb{Z}/p^{m+1}\mathbb{Z}\ar[r]
\ar@/^{3pc}/@{.>}[uu]^{\wr}& \W_{m}(k_L)\ar[r]_-{\Fb-1}&
\W_{m}(k_L)\ar[r]^-{\delta}&
\mathrm{H}^1(\mathrm{G}_{k_L},\mathbb{Z}/p^{m+1}\mathbb{Z})
\ar@{..>}[uu]_{\overline{e}:=\overline{\mathrm{e}_{p^{m}}(-,1)^{p^{m+1}}}}
\ar[r]&0}
\end{equation}
where $\G_{L_m}:=\mathrm{Gal}(L_m^{\mathrm{alg}}/L_m)$. More
explicitly $\theta_{p^m}(-,1)$ induces the identification
\begin{equation}\label{action of G on xi...}
1\mapsto\xi_m^{-1}\; :\; \mathbb{Z}/p^{m+1}\mathbb{Z}
\stackrel{\sim}{\longrightarrow}\mu_{p^{m+1}} \quad\textrm{ (cf.
\ref{theta pass to the quotient})\;,}
\end{equation}
where $\xi_m$ is the unique $p^{m+1}$-th root of $1$ satisfying
$|(\xi_m-1)-\pi_m|<|\pi_m|$ (cf. \ref{relations on equivalence
classe}). Moreover $\overline{e}$ sends
$\mathrm{H}^1(\G_{k_L},\mathbb{Z}/p^{m+1}\mathbb{Z})$ in
$\mathrm{H}^1(\G_{k_{m}},\bs{\mu}_{p^{m+1}})\subseteq
\mathrm{H}^1(\G_{L_m},\bs{\mu}_{p^{m+1}})$ via the canonical
diagram \begin{equation}\label{action of e bar on G} \xymatrix{
\G_{k_L} \ar[d]_{\alpha} & \G_{k_{m}} \ar@{_{(}->}[l]
\ar@{..>}[d]^{\overline{e}(\alpha)}  \\
\mathbb{Z}/p^{m+1}\mathbb{Z}
\ar[r]_{1\mapsto\xi_m^{-1}}^{\sim}&\bs{\mu}_{p^{m+1}} }.
\end{equation}
In other words the Artin-Schreier  character
$\gamma\mapsto\alpha(\gamma):
\G_{k_L}\to\mathbb{Z}/p^{m+1}\mathbb{Z}$ is sent by $\overline{e}$
into the Kummer character $\gamma\mapsto
\overline{e}(\alpha)(\gamma)=
\xi_m^{-\alpha(\gamma)}:\G_{k_{m}}\to\bs{\mu}_{p^{m+1}}$. In
particular $\overline{e}(\alpha)=1$ if and only if
$G_{k_m}\subseteq \mathrm{Ker}(\alpha)$.
\end{theorem}
\emph{Proof : } Let $L'/L$ be an unramified extension, and let
$\lb'=(\lambda'_0,\ldots,\lambda'_m)\in\W_m(\O_{L'})$. If $L'/L$
is not complete, but algebraic, then the series
$\theta_{p^m}(\lb',T)$, and $\mathrm{e}_{p^m}(\lb',T)^{p^{m+1}}$,
are convergent at $T=1$, since the finite extension
$L(\{\lambda'_i\}_i)/L$ is complete. By \ref{theta pass to the
quotient}, to show the commutativity it is enough to prove that
$\overline{e}$ is well-defined. Let $\lb\in\W_m(\O_L)$ be such
that $\delta(\bar{\lb})=0$ (cf. diagram
\eqref{artin-screier-diagram}). By definition, there exist
$\bs{z},\bs{\eta}\in \W_{m}(\O_L)$, such that
$\bs{\eta}=(\eta_0,\ldots,\eta_m)$, with $|\eta_j|<1$, for all
$j=0,\ldots,m$, and $\lb=\varphi(\bs{z})-\bs{z}+\bs{\eta}$. Hence
\begin{equation}
\mathrm{e}_{p^m}(\lb,1)^{p^{m+1}}=
\theta_{p^m}(\bs{z},1)^{p^{m+1}}\cdot
\mathrm{e}_{p^m}(\bs{\eta},1)^{p^{m+1}}\;.
\end{equation}
Then $\mathrm{e}_{p^m}(\lb,1)^{p^{m+1}}\in(\O_{L_m})^{p^{m+1}}$.
In other words, even if the symbol $\mathrm{e}_{p^m}(\lb,1)$ has
no meaning, the number $\mathrm{e}_{p^m}(\lb,1)^{p^{m+1}}$ is the
$p^{m+1}$-th power of the number $\theta_{p^m}(\bs{z},1)\cdot
\mathrm{e}_{p^m}(\bs{\eta},1)$ of $L_m$. Hence
$\delta_{\mathrm{Kum}}(\mathrm{e}_{p^m}(\lb,1)^{p^{m+1}})=1$.

Let us show that the map $\overline{e}$ works as indicated in the
diagram \eqref{action of e bar on G}. Let
$\alpha=\delta(\bar{\lb})$, and let $\lb\in\W_m(\O_L)$ be an
arbitrary lifting of $\bar{\lb}\in\W_m(k_{L})$. By
\ref{varphi(nu)-nu=lambda has a solution in Lunr}, an easy
induction on $m$ shows that there exists
$\bs{\nu}\in\W_m(\O_{\widehat{L}^{\mathrm{unr}}})$ such that
\begin{equation}\label{equation: varphi(nu)-nu=lambda witt}
\varphi(\bs{\nu})-\bs{\nu} = \lb\;.
\end{equation}
By definition (cf. \eqref{artin-screier-diagram}), for all
$\gamma_1\in\G_{k_L}$, we have
$\alpha(\gamma_1)=\gamma_1(\overline{\bs{\nu}})
-\overline{\bs{\nu}}\in \mathbb{Z}/p^{m+1}\mathbb{Z}$. On the
other hand, by definition, $\bar{e}(\alpha)$ is the Kummer
character of $\mathrm{G}_{L_m}$ defined by
$\mathrm{e}_{p^m}(\bs{\lb},1)^{p^{m+1}}$, and is given by
$\bar{e}(\alpha)(\gamma)=\gamma(y)/y$, for all
$\gamma\in\mathrm{G}_{L_m}$, where $y$ is an arbitrary root of the
equation $Y^{p^{m+1}}=\mathrm{e}_{p^m}(\lb,1)^{p^{m+1}}$. We let
$y:=\theta_{p^m}(\bs{\nu},1)$. Then
\begin{equation}\label{popo}
\overline{e}(\alpha)(\gamma)=\gamma(y)/y=
\gamma(\theta_{p^m}(\bs{\nu},1))/\theta_{p^m}(\bs{\nu},1)=
\theta_{p^m}(\gamma(\bs{\nu})-\bs{\nu},1)\in\bs{\mu}_{p^{m+1}},
\end{equation}
because $\gamma(\pi_m)=\pi_m$, since $\gamma\in\mathrm{G}_{L_m}$.
Now $\gamma(\bs{\nu})-\bs{\nu}\in\O^{\varphi=1}_L$, because
$\gamma(\bs{\nu})$ is again a solution of the equation
\eqref{equation: varphi(nu)-nu=lambda witt}. By \ref{theta pass to
the quotient}, the root
$\theta_{p^m}(\gamma(\bs{\nu})-\bs{\nu},1)$ depends only on the
reduction of $\gamma(\bs{\nu})-\bs{\nu}$ in $k_L^{\mathrm{sep}}$,
and is equal to $\xi_m^{-\alpha(\gamma)}$. $\Box$

\begin{lemma}\label{varphi(nu)-nu=lambda has a solution in Lunr}
Let $L$ have discrete valuation, and let
$\widehat{L}^{\mathrm{unr}}$ be the completion of the unramified
extension of $L$. Then for all $\lambda\in
\O_{\widehat{L}^{\mathrm{unr}}}$, the equation $\varphi(\nu)-\nu =
\lambda$ has a solution in $\widehat{L}^{\mathrm{unr}}$.
\end{lemma}
\emph{Proof : } The equation $\bar{v}^p-\bar{v}=\bar{\lambda}$ has
a solution in $k_{L}^{\mathrm{sep}}$, hence
$|(\varphi(v)-v)-\lambda|<1$, for all lifts $v$ of $\bar{v}$.
Since $L$ has discrete valuation, the lemma follows from an
induction on the value of the ``error'' $\eta$, in the equation
$\varphi(\nu)-\nu=\lambda+\eta$.$\Box$

\begin{theorem}\label{L' is generate by theta!!}
Let $L$ have discrete valuation. Let $\alpha=\delta(\bar{\lb})$ be
the Artin-Schreier  character defined by $\bar{\lb}\in\W_m(k_L)$
(cf. \eqref{artin-screier-diagram}). Let $k_{\alpha}/k_{L}$ be the
separable extension of $k_L$, defined by the kernel of $\alpha$,
and let $L_{\alpha}/L$ be the corresponding unramified extension.
Then
\begin{equation}
L_{\alpha}(\pi_m)=L_m(\theta_{p^m}(\bs{\nu},1))\;,
\end{equation}
where $\lb$ is an arbitrary lifting of $\bar{\lb}$ in
$\W_m(\O_L)$, and
$\bs{\nu}\in\W_m(\O_{\widehat{L}^{\mathrm{unr}}})$ is a solution
of the equation $\varphi(\bs{\nu})-\bs{\nu}=\lb$. In other words,
up to replacing $L$ by $L_m$, the extension $L_\alpha$ is
generated by $\theta_{p^m}(\bs{\nu},1)$.
\end{theorem}
\emph{Proof : } Since both $L_\alpha(\pi_m)$ and
$L_m(\theta_{p^s}(\bs{\nu},1))$ contain $L^0$ (cf. \ref{definition
of k_L^0}), and since $\varphi$ extends uniquely to $L^0$, we can
suppose $L=L^0$. In this case $\overline{e}$ is injective, $L_m/L$
is totally ramified, and $\G_{k_L}$ can be identified with
$\mathrm{G}_{k_{m}}$. Let us show the inclusion
$L_m(\theta_{p^m}(\bs{\nu},1))\subseteq L_\alpha(\pi_m)$. If
$\mathrm{G}_{k_\alpha}:=\mathrm{Gal}(k_L^{\mathrm{sep}}/k_\alpha)=
\mathrm{Ker}(\alpha)$, then the inclusion follows from the fact
that $\theta_{p^m}(\bs{\nu},1)$ is fixed by
$\G_{k_\alpha}(\subseteq\!\G_{k_{m}}\!\!\!\xrightarrow[]{\sim}\!
\mathrm{Gal}(\widehat{L}_m^{\mathrm{unr}}/L_m))$. Indeed, for all
$\gamma\in \mathrm{Gal}(\widehat{L}_m^{\mathrm{unr}}/L_m)$, we
have, as in the proof of \ref{lifting characters diagram for L
Theorem},
\begin{equation}\label{action of gamma on theta}
\gamma\left(\theta_{p^m}(\bs{\nu},1)\right)=
\theta_{p^m}(\gamma(\bs{\nu})-\bs{\nu},1)\cdot
\theta_{p^m}(\bs{\nu},1)=
\xi_m^{-\alpha(\gamma)}\cdot\theta_{p^m}(\bs{\nu},1)\;,
\end{equation}
and if $\gamma\in\G_{k_\alpha}$, we have $\alpha(\gamma)=0$. Then
$L_m(\theta_{p^m}(\bs{\nu},1))\subseteq L_\alpha(\pi_m)$. In
particular,
\begin{equation}
[L_m(\theta_{p^m}(\bs{\nu},1)):L_m]\leq
[L_\alpha(\pi_m):L_m]=[k_{\alpha,m}:k_{m}]\;,
\end{equation}
where $k_{\alpha,m}$ is the smallest field in
$k_{m}^{\mathrm{sep}}$ containing $k_{m}$ and $k_\alpha$ (i.e. the
sub-field of $k_{m}^{\mathrm{sep}}$ fixed by $\G_{k_\alpha}$
acting on $k_{m}^{\mathrm{sep}}$). The inclusion
$L_\alpha(\pi_m)\subseteq L_m(\theta_{p^m}(\bs{\nu},1))$ follows
from the equality
$[L_m(\theta_{p^m}(\bs{\nu},1)):L_m]=[k_{\alpha,m}:k_{m}]$.
Indeed, since $L_0=L$, the map $\overline{e}$ is injective. Hence
$[L_m(\theta_{p^m}(\bs{\nu},1)):L_m]=[k_\alpha:k_L]$, because
these two degrees are equal to the \emph{cardinality} of the
cyclic Galois groups generated by $\overline{e}(\alpha)$ and
$\alpha$ respectively. On the other hand, since $k_L=k_L^0$, we
have $[k_\alpha:k_L]=[k_{\alpha,m}:k_{m}]$.$\Box$

\begin{remark}\label{F_p instead of varphi}
The hypothesis of discreteness of $L$, in Theorem \ref{L' is
generate by theta!!}, and the hypothesis of existence of $\varphi$
can be removed as follows. Let $\F_p:\W_m\to\W_m$ be the map
$(\lambda_0,\ldots,\lambda_m)\mapsto
(\lambda_0^p,\ldots,\lambda_m^p)$. Replace $\varphi$ by $\F_p$,
and define
$\theta^{(\F_p)}_d(\lb,T):=\mathrm{e}_{d}(\F_p(\lb),T^p)/\mathrm{e}_d(\lb,T)$.
Then $\F_p$ is defined for \emph{all} extensions of $L$, and
commutes with the Galois action. It is easy to recover Theorems
analogous to \ref{theta pass to the quotient}, \ref{lifting
characters diagram for L Theorem}, and \ref{L' is generate by
theta!!}. In particular the analogues of diagram \eqref{lifting
characters diagram for L} is defined and functorial, on \emph{all}
complete (or algebraic) extensions of $L$. Observe that the map
$\lb\mapsto\theta^{(\F_p)}_d(\lb,T)$ is not a group morphism, but
induces again the group morphism
$1\mapsto\xi_m^{-1}:\mathbb{Z}/p^{m+1}\mathbb{Z}\stackrel{\sim}{\to}\bs{\mu}_{p^{m+1}}$
(cf. \eqref{1 to xi_m^-1}), which is the reduction of the
\emph{set} $\W_m(\O_L^{\F_p=1}):=\{\lb\in\W_m(\O_L)\;|\;
\F_p(\lb)=\lb\}$, formed by Witt vectors whose entries are $0$ or
$p-1$ roots of $1$.
\end{remark}

\subsubsection{\textbf{Application to the field $\Ed_K$}}
\begin{remark}
These methods apply to obtain a description of the Kummer
extensions of $\mathcal{E}_K$ (resp. $\Ed_K$) coming by
henselianity from an Artin-Schreier extension of $k(\!(t)\!)$ (see
below). This description is really entirely explicit, since the
Kummer generator $\theta_{p^m}(\bs{\nu},1)$ is explicitly and
directly given by the vector $\lb$. Indeed, we will give meaning
to the expression
$\theta_{p^m}(\bs{\nu},1)=\mathrm{e}_{p^m}(\lb,1)$, and we do not
need to find a solution of the equation
$\varphi(\bs{\nu})-\bs{\nu}=\lb$ (cf. definition \ref{definition*
of et(f(T),1)}, and Theorem \ref{second main identification}-(3)).
\end{remark}
The precedent theory can be applied to the field
$L=\mathcal{E}_K$, under the following assumptions on $K$:
\begin{equation}\label{hypothesis on K for AS-Kummer theory}
\left\{\begin{array}{l} (1) \textrm{ }K \textrm{ has a discrete
valuation (used in \ref{varphi(nu)-nu=lambda
has a solution in Lunr}).}\\
(2)\textrm{ There exists an \emph{absolute} Frobenius }
\sigma:K\to
K\\
\;\,\quad\textrm{(i.e. a lifting of the }p-\textrm{th power map of
}k\textrm{)}.
\end{array}\right.
\end{equation}
Fixing an absolute Frobenius of $\mathcal{E}_K$, the theory
applies without problems. Recall that we can suppress these two
hypothesis if necessary (cf. remark \ref{F_p instead of varphi}).

The situation is slightly different for the field $\Ed_K$, because
it is not complete. Nevertheless the preceding results are still
true for $\Ed_K$. Let $K$ satisfy \eqref{hypothesis on K for
AS-Kummer theory}, and fix an absolute Frobenius
 $\varphi: \O^{\dag}_{K} \to\O^{\dag}_{K}$, extending $\sigma$, by
choosing $\varphi(T)$ in $\O^{\dag}_{K}$, lifting
$t^p\in\E=k(\!(t)\!)$ (cf. \ref{what is an absolute Frobenius}).
\begin{remark}\label{varphi descent}
Since $\varphi(T)\in\O^\dag_K$ is a lifting of $t^p\in\E$, hence
there exists $0<\varepsilon_\varphi<1$ such that
$\varphi(\a_{K_m}(I^p))\subseteq \a_{K_m}(I)$, where
$I=]1-\varepsilon_\varphi,1[$ .
\end{remark}

%The fact that $K$ is discrete valued is used also in these two
%following results.
\begin{theorem}[\protect{\cite[4.2]{Cr},\cite[2.2]{Ma}}]
\label{henselian property}\label{Henselian bijection}\label{Ed is
henselian} If $K$ has discrete valuation, then $\O^{\dag}_K$ is
Henselian, hence we have a bijection
\begin{equation*}\label{lift henselian bijection}
\{\textrm{Finite unramified extensons of } \Ed_{K} \}
\stackrel{\sim}{\leadsto}  \{\textrm{Finite separable extensions
of } \E=k(\!(t)\!)\}\;.
\end{equation*}
\end{theorem}
\begin{proposition}\label{W_s(O_{K(pi_s)}^{dag}) is stable}
Let $\bs{f}(T)\in\W_m(\O_K^{\dag})$, then
 both series $\theta_{p^m}(\bs{f}(T),1)$ and
 $\mathrm{e}_{p^m}(\bs{f}(T),1)^{p^{m+1}}$ lies in $\O_{K_m}^{\dag}$.
%\begin{equation}
%\theta_{p^m}(\bs{f}(T),1)\;,\;\mathrm{e}_{p^m}(\bs{f}(T),1)^{p^{m+1}}
%\in\O_{K_m}^{\dag}\;.
%\end{equation}
Moreover if $\bs{u}(T)=(u_0(T),\ldots,u_{s}(T))
\in\W_{s}(\O_{K_s}^{\dag})$ is such that $|u_i(T)|_1 < 1$, for all
$i$, then $\e_{p^{s}}(\bs{u}(T),1)$ makes sense, and lies in
$\O_{K_s}^{\dag}$.
\end{proposition}
\emph{Proof : } Let $\varepsilon>0$ be such that
$\bs{f}(T)\in\W_{s}(\a_{K}(]1-\varepsilon,1[))$. For all compact
$J\subset]1-\varepsilon,1[$, the algebra $\a_{K}(J)$ is complete
with respect to the absolute value $\|f(T)\|_{J}:=\sup_{\rho\in
J}|f(T)|_\rho$. Hence
$\mathrm{e}_{p^m}(\bs{f}(T),1)^{p^{m+1}}\in\W_m(\a_{K_m}(J))$, for
all compact $J\subset]1-\varepsilon,1[$, and then
$\mathrm{e}_{p^m}(\bs{f}(T),1)^{p^{m+1}}\in
\W_m(\a_{K_m}(]1-\varepsilon,1[))$. On the other hand,
$\theta_{p^m}(\bs{f}(T),Z)\in 1+\pi_mZ\O_{\mathcal{E}_{K_m}}[[Z]]$
is a series in $Z$ depending only on $\bs{f}(T)$ and
$\varphi(\bs{f}(T))$. By \ref{varphi descent}, there exists
$\varepsilon'$ such that both $\bs{f}(T)$ and $\varphi(\bs{f}(T))$
lie in $\W_m(\a_{K}(]1-\varepsilon',1[))$. Hence as before
$\theta_{p^m}(\bs{f}(T),Z)\in 1+\pi_m\a_{K_m}(J)[[Z]]$, for all
compact $J\in ]1-\varepsilon',1[$, and hence
$\theta_{p^m}(\bs{f}(T),Z)\in
1+\pi_mZ\a_{K_m}(]1-\varepsilon',1[)[[Z]]$. The assertion on
$\bs{u}(T)$ follows from \ref{reduction to k}, and the same
considerations.$\Box$

\begin{corollary}\label{lift-diagr-for Ed}
The diagram \eqref{lifting characters diagram for L} can be
computed for $\Ed_K$ instead of $L$. The other assertions of
Theorems \ref{lifting characters diagram for L Theorem} and
\ref{L' is generate by theta!!} remain true (see diagram
\eqref{big diagram for Edag}). In particular, if
$\alpha=\delta(\overline{\bs{f}}(t))$ is the Artin-Schreier
character defined by $\overline{\bs{f}}(t)\in\W_m(\E)$, and if
$\E_{\alpha}/\E$ is the separable extension defined by the kernel
of $\alpha$, then the (Kummer) unramified extension of
$\Ed_{K_m}$, whose residue field is $\E_\alpha$, is
$\Ed_{K_m}(\theta_{p^m}(\bs{\nu},1))$, where $\bs{\nu}$ is a
solution of $\varphi(\bs{\nu})-\bs{\nu}=\bs{f}(T)$, for an
arbitrary lifting $\bs{f}(T)$ of $\overline{\bs{f}}(t)$.
\end{corollary}
\emph{Proof : } Let $\F_\alpha/\E$ be the separable Artin-Schreier
extension defined by $\bs{f}(T)\in\W_m(\E)$, and let
$\mathcal{F}_\alpha^\dag$ be the corresponding unramified
extension of $\Ed_K$. Let $\bs{\nu}\in \W_m(\widehat{
\mathcal{E}}_K^{\mathrm{unr}})$ be a solution of
$\varphi(\bs{\nu})-\bs{\nu}=\bs{f}(T)$. The non trivial fact is
that $\theta_{p^m}(\bs{\nu},1)$ lies in
$\mathcal{F}_\alpha^{\dag}(\pi_m)$ and not only in its completion,
say $\mathcal{F}_\alpha(\pi_m)$. In other words, we shall show
that
$\mathcal{F}_\alpha^{\dag}(\pi_m)=\Ed_{K_m}(\theta_{p^m}(\bs{\nu},1))$.
Both $\Ed_{K_m}(\theta_{p^m}(\bs{\nu},1))$ and
$\mathcal{F}_\alpha^{\dag}(\pi_m)$ are unramified over
$\mathcal{E}_{K}^{\dag,0}=\Ed_{K_m}\cap\mathcal{E}_K^{\dag,\mathrm{unr}}$,
since their completions are unramified. Moreover, by Theorem
\ref{L' is generate by theta!!}, they have the same residue field,
since this last coincides with that of their completions. By
unicity (cf. \ref{henselian property}), they are equal. $\Box$
\begin{remark}
This corollary generalize \cite[3.8]{Ma}.
\end{remark}
\begin{remark}\label{the study can be reduced to polynomials}
The study of a generic Artin-Schreier  character, given by
$\bs{f}(T)\in\W(\O_K^{\dag})$, will be reduced to the case
$\bs{f}(T)\in\W_s(\O_K[T^{-1}])$ (cf. \ref{f=u+g}, \ref{how to
compute theta explicitely}, \ref{completeness of W_s(E)}).
\end{remark}
\begin{lemma}\label{f=u+g}
Let $\bs{f}(T)\in \W_{m}(\O_{K_m}^{\dag})$, then there exist
$\widetilde{\bs{f}}(T)\in\W_{m}(\O_{K_m}[[T]][T^{-1}])$ and
$\bs{u}(T)=(u_{0}(T),\ldots,u_{m}(T))\in\W_{m}(\O_{K_m}^{\dag})$
such that $|u_{j}(T)|_1 < 1$ for all $j=0,\ldots,m$ and
$\bs{f}(T)=\bs{u}(T)+\widetilde{\bs{f}}(T)$. In particular
$\theta_{p^m}(\bs{\nu},1)=\mathrm{e}_{p^m}(\bs{u}(T),1)\cdot\theta_{p^m}(\tilde{\bs{\nu}},1)$,
where $\tilde{\bs{\nu}}$ is a solution of
$\varphi(\tilde{\bs{\nu}})-\tilde{\bs{\nu}}=\widetilde{\bs{f}}(T)$.
\end{lemma}
\emph{Proof : } This is evident for $m\!=\!0$. By induction the
lemma follows from the following relation valid for Witt vectors
in general (\cite[ch.10,$\S 1$,Lemme$4$]{Bou}):
\begin{equation}\label{rule for induction}
(f_0(T),\ldots,f_{m}(T))=(f_0(T),0,\ldots,0)+
(0,f_1(T),\ldots,f_{m}(T))\;.\Box
\end{equation}

\specialsection{\textbf{Classification of rank one differential
equations over \protect{$\R_{K_\infty}$} }} \label{classific-rk-1}

Throughout this second part, we will not need the results of
Section \ref{Dw's theta function}. Namely, $(K,|.|)$ is only a
complete ultrametric field containing $(\mathbb{Q}_p,|.|)$, and we
will not suppose that $K$ satisfies \eqref{hypothesis on K for
AS-Kummer theory}, nor that its residue field is perfect. We fix a
Lubin-Tate group $\mathfrak{G}_P$, isomorphic to
$\widehat{\mathbb{G}}_m$, and fix a generator
$\bs{\pi}=(\pi_j)_{j\geq 0}$ of the Tate module
$T(\mathfrak{G}_P)$.

We recall that $K_s=K(\pi_s)$, and that $k_{s}$ is its residue
field (cf. \ref{L_infty}). For all algebraic extension $H/K$, we
set $H_s:=H(\pi_s)$. The residue fields of $H$ and $H_s$ are
denoted by $k_H$ and $k_{H_s}$ respectively. We set
$\E_s:=k_{s}(\!(t)\!)$.

\subsubsection{}The starting point of the classification is the
equation
\begin{equation}\label{theta^p^m+1=et^p^m+1}
\theta_{p^s}(\bs{\nu},1)^{p^{s+1}}=\mathrm{e}_{p^s}(\bs{f}(T),1)^{p^{s+1}}\;,
\end{equation}
with the notations of Corollary \ref{lift-diagr-for Ed} and
diagram \eqref{big diagram for Edag}. In some cases the symbol
$\mathrm{e}_{p^s}(\bs{f}(T),1)$ does make sense, and the
interesting ``Kummer generator'' $\theta_{p^s}(\bs{\nu},1)$ is
equal to $\mathrm{e}_{p^s}(\bs{f}(T),1)$. We will show that all
rank one solvable differential equations over $\R_{K_m}$ admit, in
some basis, such an exponential as solution.

\begin{definition}\label{definition* of et(f(T),1)}
\index{f-(T)@$\bs{f}^-(T)=(f_0^-(T),\ldots,f_s^-(T))\in\W_s(T^{-1}\O_K[T^{-1}])$}
Let $\bs{f}^-(T)\in\W_s(T^{-1}\O_K[T^{-1}])$, then we set
\begin{equation}
\index{e(f-,1)@$\mathrm{e}_{p^m}(\bs{f}^-(T),1):=
\exp(\sum_{j=0}^{m}\pi_{m-j}\cdot\phi_j^-(T)\cdot \frac{1}{p^j})$}
\mathrm{e}_{p^s}(\bs{f}^-(T),1):=
\exp\Bigl(\pi_s\phi^-_0(T)+\pi_{s-1}\frac{\phi^-_1(T)}{p}+
\cdots+\pi_0\frac{\phi^-_s(T)}{p^s}\Bigr)\;,
\end{equation}
where $\phi_{j}^-(T)$ is the $j$-th phantom component of
$\bs{f}^{-}(T)=(f_0^-(T),\ldots,f_s^-(T))$.
\end{definition}
\begin{remark}\label{how to compute theta explicitely} Clearly $\phi_{j}^-(T)$ lies in
$T^{-1}\O_K[T^{-1}]$, for all $j=0,\ldots,s$, and hence the
expression \ref{definition* of et(f(T),1)} converges
$T^{-1}$-adically. Moreover,
\begin{equation}\label{et(f(T),1) is convergent for T > 1}
\mathrm{e}_{p^s}(\bs{f}^{-}(T),1)=\prod_{j=0}^s E_{s-j}(f^-_j(T))
\in 1 + \pi_s T^{-1}\O_{K_s}[[T^{-1}]]\;.
\end{equation}
In particular, $\mathrm{e}_{p^s}(\bs{f}^-(T),1)$ is convergent for
$|T|>1$. As mentioned in Remark \ref{the study can be reduced to
polynomials}, Lemmas \ref{f=u+g}, \ref{decomposition of CW(E)},
and \ref{criteria of solvability lemma}, will be useful to reduce
the study of $\theta_{p^s}(\bs{\nu},1)$, with
$\varphi(\bs{\nu})-\bs{\nu}=\bs{f}(T)$, to the case in which
$\bs{f}(T)\in \W_s(\O_K[T^{-1}])$.
\end{remark}

\subsection{Survey of the Results} \label{Classification of rank
one differential equats}
\begin{remark} For all algebraic extensions $H/K$, the function (cf. \ref{definition* of et(f(T),1)})
 \begin{equation}
 \bs{f}^-(T)\mapsto\mathrm{e}_{p^s}(\bs{f}^-(T),1)
 \end{equation}
 defines a group morphism
 (as we can see by considering the phantom components)
\begin{equation}
\mathrm{e}_{p^s}(-,1) : \W_s(T^{-1}\O_{H}[T^{-1}])\longrightarrow
1+\pi_sT^{-1}\O_{H_s}[[T^{-1}]]\;.
\end{equation}
Indeed $\bs{f}^-(T)$ involve only a finite numbers of coefficients
of $H$, then the series $\mathrm{e}_{p^s}(\bs{f}^-(T),1)$ lies in
a finite (and hence complete) extension of $K$.

Let $\bs{f}^-(T)\in\W_s(T^{-1}\O_H[T^{-1}])$, we set
\begin{equation}
\L(0,\bs{f}^-(T)) =
\d-\partial_{T,\log}(\mathrm{e}_{p^s}(\bs{f}^-(T),1)).
\end{equation}
Observe that $1+\pi_sT^{-1}\O_{H_s}[[T^{-1}]]$ is not contained in
$\mathcal{E}_{H_s}=\mathcal{E}_K\otimes_K H_s$. However, every
series in this multiplicative group is convergent for $|T|>1$ (cf.
\ref{how to compute theta explicitely}). Then, by \ref{transfer}
and by continuity of the radius, $\L(0,\bs{f}^-(T))$ is solvable
over $\R_{H_s}$.
\end{remark}
\begin{theorem}[Main Theorem]\label{Introductive foundamental Theorem}
Let $M$ be a rank one solvable differential module over
$\R_{K_\infty}$ (i.e. over $\R_{K_m}$, for some $m\geq 0$, or over
$\R_K$ (cf. \ref{Picsol})). Then there exists a basis of $M$ such
that
\begin{enumerate}
 \item the $1\times 1$ matrix of the derivation of $M$ lies in
 $\a_K(]0,\infty])\cap \O_K[[T]][T^{-1}]$;
 \item there exist an $s\geq 0$, and a Witt vector
 $\bs{f}^-(T)\in\W_s(T^{-1}\O_{K_s}[T^{-1}])$ such that the
 Taylor solution (cf. \eqref{s_x(T)}) of $M$, at $\infty$, is
 \begin{equation}\label{solution of M(a_0,alpha)}
 T^{a_0}\cdot\mathrm{e}_{p^s}(\bs{f}^-(T),1)\;,
 \end{equation}
 with $a_0\in\mathbb{Z}_p$. In particular $M$ is defined (in this
basis) by the operator
\begin{eqnarray}\index{Laf@$\L(a_0,\bs{f}^-(T))=\d-\partial_{T,\log}(\mathrm{e}_{p^s}(\bs{f}^-(T),1))$}
\label{explicit L(a_0,f^-)}
\L(a_0,\bs{f}^-(T))&:=&\d-\partial_{T,\log}\left(T^{a_0}\cdot
\mathrm{e}_{p^s}(\bs{f}^-(T),1)\right)\nonumber\\
&=&\d-a_0+\sum_{j=0}^{s}\pi_{s-j}\sum_{i=0}^{j}f_{i}^-(T)^{p^{j-i}}
\partial_{T,\log}(f_{i}^-(T)).
\end{eqnarray}
\end{enumerate} Moreover the isomorphism class of $M$ depends
bijectively on
\begin{itemize}
 \item[-] the class of $a_0$ in $\mathbb{Z}_p/\mathbb{Z}$;
 \item[-] the Artin-Schreier  character $\alpha:=\delta(\overline{\bs{f}^-}(t))$ defined
by the reduction $\overline{\bs{f}^-}(t)\in\W_s(E_s)$ of
$\bs{f}^-(T)$.
\end{itemize}
\end{theorem}

\begin{definition}\label{M(a_0,alpha)}
We will denote indifferently by
$\M(a_0,\alpha)$,\index{M(a,f)@$\M(a_0,\bs{f}^-(T))$,
$\M(a_0,\overline{\bs{f}^-}(T))$, $\M(a_0,\alpha)$ }
$\M(a_0,\overline{\bs{f}^-}(t))$ or $\M(a_0,\bs{f}^-(T))$, the
differential module defined by $\L(a_0,\bs{f}^-(T))$.
\end{definition}
Assume the point $(1)$ and $(2)$ of the Theorem \ref{Introductive
foundamental Theorem}. Then the last assertion can be translated
in terms of $\bs{\pi}$-exponentials as follow. Recall that
$p=\mathrm{w}$ (cf. \ref{overconvergence of theta}).
\begin{theorem}\label{rewritten in terms of P-exp}
Let $\bs{f}^-(T)\in\W_s(T^{-1}\O_{K_s}[T^{-1}])$, and let
$\overline{\bs{f}^-}(t)\in \W_s(t^{-1}k_{s}[t^{-1}])$ be its
reduction. Then
\begin{enumerate}
\item[(3)] If
$\widetilde{\bs{f}}^-(T)\in\W_s(T^{-1}\O_{K_s}[T^{-1}])$ is
another lifting of $\overline{\bs{f}^-}(t)$, then
\begin{equation}
\frac{\mathrm{e}_{p^s}(\bs{f}^-(T),1)}{\mathrm{e}_{p^s}(\widetilde{\bs{f}}^-(T),1)}=
\mathrm{e}_{p^s}(\bs{f}^-(T)-\widetilde{\bs{f}}^-(T),1)
\end{equation}
is convergent for $|T|>1-\varepsilon$, for some $\varepsilon > 0$
(i.e. lies in $\R_{K_s}$).

\item[(4)] If $\bs{f}^-_{(\Fb)}(T)\in\W_s(T^{-1}\O_{K_s}[T^{-1}])$
is an arbitrary lifting of $\Fb(\overline{\bs{f}^-}(t))$, then
\begin{equation}
\frac{\mathrm{e}_{p^s}(\bs{f}^-_{(\Fb)}(T),1)}{\mathrm{e}_{p^s}(\bs{f}^-(T),1)}=
\mathrm{e}_{p^s}(\bs{f}^-_{(\Fb)}(T)-\bs{f}^-(T),1)
\end{equation}
is convergent for $|T| > 1-\varepsilon'$, for some
$\varepsilon'>0$ (i.e. lies in $\R_{K_s}$).

\item[(5)] Conversely the function
$\mathrm{e}_{p^s}(\bs{f}^-(T),1)$ lies in $\R_{K_s}$ if and only
if the equation $\Fb\left(\bar{\bs{\nu}}^-\right) -
\bar{\bs{\nu}}^- = \overline{\bs{f}^-}(t)$ has a solution
$\bar{\bs{\nu}}^-\in\W_s(t^{-1}k_{s}[t^{-1}])$.
\end{enumerate}
\end{theorem}
\begin{notation}
The point $(5)$ will be called  the Frobenius Structure Theorem.
\end{notation}
\subsubsection{}By the Main Theorem \ref{Introductive foundamental
Theorem}, definition \ref{definition* of et(f(T),1)} and by rule
\ref{d + g + gtilde}, it follows that, for all $s\geq 0$, and for
all algebraic extensions $H/K$, we have an exact sequence of
abelian groups (functorial on the algebraic extensions $H$ of $K$)
\begin{equation}
\W_s(t^{-1}k_{H}[t^{-1}])
\xrightarrow{\Fb-1}\W_s(t^{-1}k_H[t^{-1}]) \xrightarrow{\M(0,-)}
\mathrm{Pic}^{\mathrm{sol}}(\R_{H_s})\;.
\end{equation}
On the other hand, it follows from the definition \ref{definition*
of et(f(T),1)}, that we have
\begin{equation}\label{invariant by V.}
\mathrm{e}_{p^{s+1}}(\V(\bs{f}^-(T)),1) =
\mathrm{e}_{p^s}(\bs{f}^-(T),1)\;.
\end{equation}
Hence, for all $s\geq 0$, we have the following functorial
commutative diagram
\begin{equation}
\xymatrix{\W_s(t^{-1}k_H[t^{-1}])\ar[r]^{\Fb-1}\ar@{}[dr]|{\odot}\ar[d]_{\V}&
\W_s(t^{-1}k_H[t^{-1}])\ar[r]^-{\M(0,-)} \ar[d]_-{\V} &
\mathrm{Pic}^{\mathrm{sol}}(\R_{H_{s+1}}) \\
\W_{s+1}(t^{-1}k_H[t^{-1}])\ar[r]^{\Fb-1}&
\W_{s+1}(t^{-1}k_H[t^{-1}])\ar[ur]_-{\M(0,-)} & }
\end{equation}
This shows, by passing to the inductive limit, that we again have
an exact sequence
\begin{equation}
\bs{\mathrm{CW}}(t^{-1}k_H[t^{-1}])\xrightarrow{\Fb-1}
\bs{\mathrm{CW}}(t^{-1}k_H[t^{-1}])\xrightarrow{\M(0,-)}
\mathrm{Pic}^{\mathrm{sol}}(\R_{H_\infty})\;.
\end{equation}
The group $\mathbb{Z}_p/\mathbb{Z}$ has no $p$-torsion element. On
the other hand, every element of
$\bs{\mathrm{CW}}(t^{-1}k_H[t^{-1}])$ is killed by a power of $p$.
Since we are assuming that all solutions are of the form
\eqref{solution of M(a_0,alpha)}, this proves the following
\begin{lemma}
Let $H$ be an algebraic extension of $K_\infty$. The image of
$\M(0,-)$ is the sub-group of the $p$-torsion elements of
$\mathrm{Pic}^{\mathrm{sol}}(\R_{H})$, and if $H/K_\infty$ is
Galois, then $\mathrm{Pic}^{\mathrm{sol}}(\R_H)$ is isomorphic, as
$\mathrm{Gal}(H/K_\infty)$-module, to the direct sum of
$\mathbb{Z}_p/\mathbb{Z}$ with the image of $\M(0,-)$. $\Box$
\end{lemma}
\begin{corollary}\label{explicit description of Pic}
The map $(a_0,\alpha)\mapsto \M(a_0,\alpha)$ induces an
isomorphism
\begin{equation}\label{Z/Z oplus CW/CW}
 \mathbb{Z}_p/\mathbb{Z}\oplus
\frac{\bs{\mathrm{CW}}(t^{-1}k_{\infty}[t^{-1}])}{(\Fb-1)
 \bs{\mathrm{CW}}(t^{-1}k_{\infty}[t^{-1}])}\quad
 \xrightarrow[\sim]{\M(-,-)}\quad \mathrm{Pic}^{\mathrm{sol}}(\R_{K_\infty})\;.
\end{equation}
\end{corollary}
\emph{Proof : } By Galois descent $\M(-,-)$ induces an
isomorphism, with
$k_{\infty}^{\mathrm{perf}}:=(k^{\mathrm{alg}})^{\mathrm{Gal}(k^{\mathrm{alg}}/k_\infty)}$
instead of $k_\infty$. But actually, the co-vector quotient is
invariant under inseparable extension of $k_\infty$ as explained
in subsection \ref{PAS is inseparable invariant} below. $\Box$
\subsubsection{} On the other hand, it is well known that (cf.
lemma \ref{decomposition of CW(E)} and \ref{completeness of
W_s(E)})
 \begin{equation*}
 \mathrm{H}^1\left(
 \mathrm{Gal}\bigl(k_\infty(\!(t)\!)^{\mathrm{sep}}/
 k_\infty(\!(t)\!)\bigr),
 \mathbb{Q}_p/\mathbb{Z}_p\right)
 =
 \PAS(k_\infty)\oplus
 \mathrm{H}^1\left(
 \mathrm{Gal}\bigl(k_\infty^{\mathrm{sep}}/k_\infty\bigr),
 \mathbb{Q}_p/\mathbb{Z}_p\right),
 \end{equation*}
where $\PAS(k_\infty)$ is the character group of
$\mathcal{P}_{\E_\infty}$, with $\E_\infty=k_\infty(\!(t)\!)$ (cf.
\ref{first remark on PAS}). More precisely we have the following
(for a more convenient description of $\PAS(\kappa)$ see
\ref{description more handy})
\begin{lemma}
For all fields $\kappa$ of characteristic $p>0$, one has
\begin{equation}\label{PKAPPA}
\PAS(\kappa)= \frac{\bs{\mathrm{CW}}(t^{-1}
\kappa[t^{-1}])}{(\Fb-1)\bs{\mathrm{CW}}(t^{-1}\kappa[t^{-1}])}\;.
\end{equation}
\end{lemma}
\emph{Proof : } This will follow from \ref{decomposition of CW(E)}
and \ref{completeness of W_s(E)}. $\Box$
\subsubsection{}\label{PAS is inseparable invariant}
Furthermore we have
$\PAS(k_\infty^{\mathrm{perf}})=\PAS(k_\infty)$, because, by
Remark \ref{CW(R)=CW(R^p)} (or \ref{VF is the same of F}), the
Artin-Schreier  complex is stable under purely inseparable
extensions, that is $\mathrm{Gal}(k_\infty^{\mathrm{perf,sep}}/
k_{\infty}^{\mathrm{perf}})=
\mathrm{Gal}(k_{\infty}^{\mathrm{sep}}/k_\infty)$. In other words,
for all $r\geq 0$, the co-vectors
$\overline{\bs{f}^-}(t)=(\ldots,0,\overline{f^-_0}(t),\ldots,
\overline{f^-_s}(t))$ and
$\Fb^r(\overline{\bs{f}^-}(T))=(\ldots,0,
\overline{f^-_0}(t)^{p^r},\ldots,\overline{f^-_s}(t)^{p^r})$ have
the same image in the right hand quotient of equation
\eqref{PKAPPA}.

\subsection{Proofs of the statements}\label{Proofs of the stat} We
first prove the statements $(3)$, $(4)$, and $(5)$ of Theorem
\ref{rewritten in terms of P-exp}. The idea is to express
$\mathrm{e}_{p^s}(\bs{f}^-(T),1)$ as a product of
$\bs{\pi}$-exponentials of the type $\mathrm{e}_d(\lb,T^{-1})$.
The main tool will be the notion of $s$-\emph{co-monomial}
\index{lambda T@$\lb T^{n}=s$-co-monomial (cf. \ref{co-monomials
and PAS}), or monomial (cf. \ref{MONOMIAL})}which reduce the study
to $\bs{\pi}$-exponentials (see equation \eqref{decomposition as
product of P-exp}). The principal lemma will be
\ref{transformation}.
\begin{definition} \label{co-monomials and PAS}
Let $H/K$ be an algebraic extension. Let $d=np^m>0$, $(n,p)=1$.
Let $s\geq 0$. We will call $s$-\emph{co-monomial} of degree $-d$
relative to $\lb:=(\lambda_0,\ldots,\lambda_m)\in\W_m(\O_H)$ the
Witt vector in $\W_{s}(T^{-1}\O_{H}[T^{-1}])$
\begin{equation}
\begin{array}{rclcl}
\lb T^{-d}&:=&(\overbrace{0,\ldots,0}^{s-m},\lambda_0T^{-n},
\lambda_1T^{-np},\ldots, \lambda_mT^{-d})&
\textrm{ if }& m\leq s\;, \\
\lb
T^{-d}&:=&(\lambda_{m-s}T^{-np^{m-s}},\lambda_{m-s+1}T^{-np^{m-s+1}},\ldots,
\lambda_mT^{-d})& \textrm{ if }& m\geq s\;.\\
\end{array}
\end{equation}
We denote by $\W_s^{(-d)}(\O_H)$
\index{W_s^(d)@$\W_s^{(-d)}(\O_H)$, $\W_s^{(-d)}(k_H)$} the
sub-group of $\W_s(T^{-1}\O_H[T^{-1}])$ formed by $s$-co-monomials
of degree $-d$, and by $\W_s^{(-d)}(k_H)$ its image in
$\W_s(t^{-1}k_H[t^{-1}])$.
\end{definition}

\begin{remark}\label{W^d simto W_min(m,s)}
By looking at the phantom components we find an isomorphism of
groups
$\W_s^{(-d)}(\O_H)\stackrel{\sim}{\to}\W_{\min(s,m)}(\O_H)$, and
hence $\W_s^{(-d)}(k_H)\stackrel{\sim}{\to}\W_{\min(s,m)}(k_H)$.
\end{remark}

\begin{lemma} \label{transformation}
Let now $H/K$ be an algebraic extension. Let $d=np^m>0$,
$(n,p)=1$, let $s\geq 0$, and let
$\lb:=(\lambda_0,\ldots,\lambda_m)\in\W_m(\O_H)$. If $m\leq s$, we
have
\begin{equation}
\mathrm{e}_{p^s}(\lb T^{-d},1)= \mathrm{e}_{d}(\lb,T^{-1})\;.
\end{equation}
\end{lemma}
\emph{Proof : } The phantom vector of
$(\underbrace{0,\ldots,0}_{s-m},\lambda_0T^{-n},\lambda_1T^{-np},\ldots,\lambda_mT^{-d})$,
is
\begin{equation}\label{phantom of m_d}
\ph{0,\ldots,0,p^{s-m}\phi_0T^{-n},p^{s-m}\phi_1T^{-np},\ldots,p^{s-m}\phi_mT^{-d}},
\end{equation}
where $\ph{\phi_0,\cdots,\phi_m}$ is the phantom vector of
$(\lambda_0,\ldots,\lambda_m)$. The proof follows immediately from
the definitions \ref{definition* of et(f(T),1)} and
\eqref{definition of et_d(lb,T)}. $\Box$

\begin{definition}
For all algebraic extensions $H/K$ we set $\E_H:=k_H(\!(t)\!)$.
\end{definition}
\begin{lemma}\label{decomposition of CW(E)}
For all $s\geq 0$, there is a (functorial) decomposition
\begin{enumerate}
\item $\W_s(\E_H)=\oplus_{d>0} \W_s^{(-d)}(k_H)\oplus
\W_s(k_H)\oplus \W_s(tk_H[[t]])\;;$ \item
$\W_s(\O_H[[T]][T^{-1}])= \oplus_{d>0}\W_s^{(-d)}(\O_H)
\oplus\W_s(\O_H)\oplus \W_s(T\O_H[[T]]).$
\end{enumerate}
\end{lemma}
\emph{Proof : } Let $s=0$, then $k_H((t))=\oplus_{d>0}k_H
t^{-d}\oplus k_H\oplus tk_H[[t]]$. The proof follows easily by
induction from \eqref{rule for induction}.$\Box$
\begin{remark}\label{remark-decomp}
Witt vectors in $\oplus_{d>0}\W_s^{(-d)}(\O_H)$ (resp.
$\W_s(\O_H)$, $\W_s(T\O_K[[T]])$) have their phantom components in
$T^{-1}\O_H[T^{-1}]$ (resp. $\O_H$, $T\O_K[[T]]$).
\end{remark}
\begin{corollary}
We have a (functorial) decomposition
$$\W_s(T^{-1}\O_H[T^{-1}])=\oplus_{d>0}\W_s^{(-d)}(\O_H)\;, \quad
\W_s(t^{-1}k_H[t^{-1}])=\oplus_{d>0}\W_s^{(-d)}(k_H)\;.$$
\end{corollary}
\emph{Proof : } The inclusion $\subseteq$ follows by Remark
\ref{remark-decomp}. Since all monomials belong to
$\W_s(T^{-1}\O_H[T^{-1}])$ we have the inclusion $\supseteq$. The
right hand equality follows from the first one by reduction.
$\Box$

\begin{definition} \label{dec -0+}
For all $\bs{f}(T)\in\W_s(\O_H[[T]][T^{-1}])$, we will denote by
\begin{equation}
\bs{f}(T)=\bs{f}^-(T)+\bs{f}_{0}+\bs{f}^+(T)
\end{equation}\index{f+@$\bs{f}^+(T)$, $\bs{f}^-(T)$, $\bs{f}(T)=\bs{f}^-(T)+\bs{f}_{0}+\bs{f}^+(T)$}
the unique decomposition of $\bs{f}(T)$ satisfying
$\bs{f}^-(T)\in\W_s(T^{-1}\O_H[T^{-1}])$, $\bs{f}_0\in\W_s(\O_H)$,
$\bs{f}^+(T)\in\W_s(T\O_H[[T]])$ (cf. \ref{decomposition of
CW(E)}). The same notation will be used for a Witt vector
$\overline{\bs{f}}(t)\in\W_s(\E_H)$.
\end{definition}

\begin{remark}
By \eqref{artin-screier-diagram}, we then have a corresponding
decomposition of $\alpha:=\delta(\overline{\bs{f}}(t))$, i.e.
$\alpha=\alpha^-+\alpha_0$, ($\alpha^+=0$ by \ref{completeness of
W_s(E)}), with $\alpha^-=\delta(\overline{\bs{f}^-}(t))$, and
$\alpha_0=\delta(\overline{\bs{f}}_0)$. This shows that
$\mathrm{Gal}(\E_H^{\mathrm{sep}}/\E_H)^{\mathrm{ab}}
=\mathrm{Gal}(k_H^{\mathrm{sep}}/k_H)^{\mathrm{ab}}\oplus
\mathcal{I}_{\E_H}^{\mathrm{ab}}$, where $\E_H=k_H(\!(t)\!)$ .
\end{remark}
\begin{proposition} \label{completeness of W_s(E)}
$\W_s(tk_H[[t]])\subseteq (\Fb-1)\W_s(tk_H[[t]])$, for all $s\geq
0$.
\end{proposition}
\emph{Proof : } Since $\E_H$ is complete, by \ref{completeness of
W_m(R)}, $\W_s(\E_H)$ is complete. Let
$\overline{\bs{f}}^+(t)\in\W_s(tk_H[[t]])$. Then the series
$\overline{\bs{g}}^+(t):=-\sum_{i\geq
0}\Fb^i(\overline{\bs{f}}^+(t))$ is Cauchy for this topology, and
hence converges in $\W_s(\E_H)$. Moreover
$\overline{\bs{f}}^+(t)=\Fb(\overline{\bs{g}}^+(t))-
\overline{\bs{g}}^+(t)$.$\Box$

\begin{remark}
Let $H/K$ be an algebraic extension and let $\bs{f}^-(T)\in
\W_s(t^{-1}\O_H[t^{-1}])$. Let $v_p(-)$ be the $p$-adic valuation
(namely $v_p(d)=m$ if $d=np^m$, $(n,p)=1$). Let
$\bs{f}^-(T)=\sum_{d>0}\lb_dT^{-d}$, with
$\lb_d\in\W_{v_p(d)}(k_H)$ be its decomposition in
$s$-co-monomials of degree $-d$. We can suppose $s\gg 0$ (cf.
\eqref{invariant by V.}), then
\begin{equation}\label{decomposition as product of P-exp}
\mathrm{e}_{p^s}(\bs{f}^-(T),1)=\mathrm{e}_{p^s}(\sum_{d>0}\lb_dT^{-d},1)=
\prod_{d>0}\mathrm{e}_{p^s}(\lb_dT^{-d},1)\stackrel{\ref{transformation}}{=}
\prod_{d>0}\mathrm{e}_d(\lb_d,T^{-1}).
\end{equation}
Then $\mathrm{e}_{p^s}(\bs{f}^-(T),1)$ is a (finite) product of
elementary $\bs{\pi}$-exponentials. In terms of differential
modules, we have
$\M(0,\bs{f}^-(T))=\otimes_{d>0}\M(0,\lb_dT^{-d})$. Hence, by
\ref{d + g + gtilde}, the study can be reduced to
$\bs{\pi}$-exponentials.
\end{remark}

\subsubsection{\textbf{Proof of the statements $(3),(4),(5)$ of Theorem
\ref{rewritten in terms of P-exp}:}}
\begin{notation}
For all $d>0$, we set $d=np^{m}$, with $(n,p)=1$ and $v_p(d):=m$.
In the sequel the letters $n$ and $m$ will indicate always this
decomposition.
\end{notation}
By Lemma \ref{transformation}, for all $d$ appearing in the
product \eqref{decomposition as product of P-exp}, we have (cf.
\ref{eq: L_d})
 \begin{equation}\label{transl}
\L_d(\lb_d)=\L(0,\lb_d T^{-d})\qquad,
\qquad\Mt_d(\lb_d)=\M(0,\lb_d T^{-d})\;,
 \end{equation}
where  $\lb_d T^{-d}$ is the $s$-co-monomial of degree $-d$
attached to $\lb_d\in\W_{v_p(d)}(\O_H)$ (cf. \ref{co-monomials and
PAS}). Actually, by the rule \eqref{invariant by V.}, we can
suppose $s\gg v_p(d)\!=\!m$, for all $d>0$ appearing in the
(finite) product \eqref{decomposition as product of P-exp}.

The assertions $(3)$ and $(4)$ are consequences of the Reduction
Theorem \ref{reduction to k}, and the Frobenius Structure Theorem
for $\bs{\pi}$-exponentials \ref{overconvergence of theta},
respectively. Let us prove the assertion $(3)$. We decompose
$\bs{f}^-(T)-\widetilde{\bs{f}}^-(T)$ in $s$-co-monomials of
degree $-d$,
$\bs{f}^-(T)-\widetilde{\bs{f}}^-(T)=\sum_{d}\lb_{d}T^{-d}$, with
$\lb_d\in\W_{v_p(d)}(\O_H)$ (cf. \ref{decomposition of CW(E)}).
Then
\begin{equation}
\mathrm{e}_{p^s}(\bs{f}^-(T)-\widetilde{\bs{f}}^-(T),1)=
\prod_{d>0}\mathrm{e}_{p^s}(\lb_dT^{-d},1)\stackrel{\ref{transformation}}{=}
\prod_{d>0}\mathrm{e}_{d}(\lb_d,T^{-1})\;.
\end{equation}
The over-convergence of
$\mathrm{e}_{p^s}(\bs{f}^-(T)-\widetilde{\bs{f}}^-(T),1)$ will
result from the over-convergence of every
$\mathrm{e}_d(\lb_d,T^{-1})$. In order to apply the Reduction
Theorem \ref{reduction to k}, we shall show that the reduction
$\overline{\lb_d}$ of $\lb_d$ is $0$, for all $d>0$. Since the
reduction of $\bs{f}^-(T)-\widetilde{\bs{f}}^-(T)$ is $0$, it
follows from Lemma \ref{decomposition of CW(E)} that the reduction
of $\lb_dT^{-d}$ in $\W_s^{(-d)}(k_H)$ is $0$, for all $d>0$. By
Remark \ref{W^d simto W_min(m,s)}, for all $d>0$, we have an
isomorphism
$\lb_dT^{-d}\mapsto\lb_d:\W_s^{(-d)}(\O_H)\stackrel{\sim}{\to}\W_{v_p(d)}(\O_H)$.
Hence $\overline{\lb_d}=0$ in $\W_{v_p(d)}(k_H)$, for all $d>0$.

The proof of $(4)$ follows the same lines. Namely, by the
assertion $(3)$, the isomorphism class of $\M(0,\bs{f}^-(T))$
depends only on the reduction
$\overline{\bs{f}^-}(t)\in\W_s(t^{-1}k_H[t^{-1}])$ of
$\bs{f}^-(T)$. As usual, we decompose
$\overline{\bs{f}^-}(t)=\sum_{d>0}\bar{\lb}_dt^{-d}$, with
$\bar{\lb}_dt^{-d}\in\W_s^{(-d)}(k_H)$. The morphism
$\Fb:\W_s(\E_H)\to\W_s(\E_H)$ sends the monomial
$\bar{\lb}_dt^{-d}$ into $\Fb(\bar{\lb}_d)t^{-pd}$, Hence
$\Fb(\overline{\bs{f}^-}(t))=\sum_{d>0}\Fb(\bar{\lb}_d)t^{-pd}$.
Then
\begin{equation*}
\M(0,\overline{\bs{f}^{-}}(t))\stackrel{\sim}{\to}
\otimes_{d>0}\M_d(\bar{\lb}_d) \xrightarrow[\ref{overconvergence
of theta}]{\sim} \otimes_{d>0}\M_{pd}(\V\Fb(\bar{\lb}_d))
\stackrel{\sim}{\to}\M(0,\Fb(\overline{\bs{f}^-}(t)))\;,
\end{equation*}
where the last isomorphism follows from the fact that
$\V(\Fb(\bar{\lb}_d)t^{-pd})$ and $\Fb(\bar{\lb}_d)t^{-pd}$ define
the same differential module (cf. \eqref{invariant by V.}).

The proof of the assertion $(5)$ of Theorem \ref{rewritten in
terms of P-exp} follows from assertions $(3)$ and $(4)$ of Theorem
\ref{rewritten in terms of P-exp} in the following way. Suppose
that $\mathrm{e}_{p^s}(\bs{f}^-(T),1)$ is over-convergent. We want
to show that the equation $\Fb(\bar{\bs{\nu}})-\bar{\bs{\nu}}=
\overline{\bs{f}^-}(t)$ has a solution $\bar{\bs{\nu}}\in
\W_s(t^{-1}k_{H}[t^{-1}])$. In other words, we shall show that
$\overline{\bs{f}^-}(t)$ belongs to
$(\Fb-1)\W_s(t^{-1}k_{H}[t^{-1}])$. Let us write
$\bs{f}^-(T)=\sum_{d>0}\lb_dT^{-d}$ as a (finite) sum of
$s$-co-monomials. We need to replace $\bs{f}^-(T)$ by a more
convenient Witt vector.
\begin{definition}\label{pure witt vector}\index{pure@ pure Witt vector}
A Witt vector $\bs{f}_p^-(T)\in\W_s(T^{-1}\O_H[T^{-1}])$ is called
\emph{pure} if its decomposition in $s$-co-monomials is a (finite)
sum of the type
\begin{equation}\label{pure witt vector example of}
\bs{f}_p^-(T)=\sum_{n\in\J}\lb_{np^{m(n)}}T^{-np^{m(n)}}\;,
\end{equation}
where $\J=\{n\in\mathbb{Z}\;|\; (n,p)=1 \;,\; n>0 \}$, and
$\lb_{np^{m(n)}}\in\W_{m(n)}(\O_{H})$.
\end{definition}
\begin{remark}\label{remark on pure Witt vector}

$\partial_{T,\log}(\mathrm{e}_{p^s}(\bs{f}_p^-(T),1))\!=
\!\sum_{n\in\J}\!\!\!-n\!\sum_{j=0}^{m(n)}
\pi_{m(n)-j}\phi_{np^{m(n)},j}T^{-np^{j}}\!\!$, where
$\ph{\phi_{np^{m(n)},0},\ldots,\phi_{np^{m(n)},m(n)}}$ is the
phantom vector of $\lb_{np^{m(n)}}$. In this case the coefficients
of the differential equation are simpler and directly related to
$\lb_{np^{m(n)}}$ instead of $\lb_{np^{m(n)}}T^{-np^{m(n)}}$. This
will be useful for explicit computations (cf. \ref{no equation
over abs unr}).
\end{remark}
\begin{lemma}\label{Irreg of a pure witt-lemma}
Let $\bs{f}_p^-(T)\in \W_{s}(T^{-1}\O_H[T^{-1}])$ be a pure Witt
vector. The exponential $\mathrm{e}_{p^s}(\bs{f}_p^-(T),1)$ is
 over-convergent if and only if $\overline{\bs{f}_p^-}(t)=0$.
 Moreover,
 \begin{equation}\label{Irreg of a pure witt}
 \mathrm{Irr}\left(\M(0,\bs{f}_p^-(T))\right) =
 \max_{n\in\J}\;\mathrm{Irr}\left(\M_{np^{m(n)}}(\lb_{np^{m(n)}})\right).
 \end{equation}
\end{lemma}
\emph{Proof : } Write $\M(0,\bs{f}_p^-(T))=
\otimes_{n\in\J}\M(0,\lb_{np^{m(n)}}T^{-np^{m(n)}})$. The
irregularity of
$\M(0,\lb_{np^{m(n)}}T^{-np^{m(n)}})\xrightarrow[\eqref{transl}]{\sim}
\Mt_{np^{m(n)}}(\lb_{np^{m(n)}})$ is, by Theorem \ref{reduction to
k}, a number belonging to the set $\{0\} \cup \{n\cdot p^m \;|\;
m\geq 0\}$. Hence, for different values of $n$, we have different
values of the $p$-adic slope of $\M_{np^{m(n)}}(\lb_{np^{m(n)}})$.
Theorem \ref{tensor product with different radius} then implies
the equation \eqref{Irreg of a pure witt}. Suppose now that
$\mathrm{e}_{p^s}(\bs{f}_p^-(T),1)$ is over-convergent, then this
irregularity is equal to $0$. Hence all
$\M_{np^{m(n)}}(\lb_{np^{m(n)}})$ are trivial, and
$\mathrm{e}_{p^s}(\lb_{np^{m(n)}}T^{-np^{m(n)}},1)$ is
over-convergent (i.e. lies in $\R_H$), for all $n\in\J$. By
Theorem \ref{reduction to k} this implies
$\overline{\lb_{np^{m(n)}}}=0$, for all $n\in\J$. $\Box$

 Assertion $(5)$ of Theorem \ref{rewritten in terms of P-exp} follows then
 by point $(1)$ of the following
\begin{lemma}\label{isomorphic to a pure vector}
Let $\bs{f}^-(T)\in\W_s(T^{-1}\O_H[T^{-1}])$. Then
\begin{enumerate}
\item There exists a pure Witt vector
$\bs{f}_p^{-}(T)\in\W_s(T^{-1}\O_H[T^{-1}])$ such that
 \begin{equation}
 \overline{\bs{f}^-(T)-\bs{f}_p^-(T)}\in
 (\Fb-1)\W_s(t^{-1}k_H[t^{-1}])\;.
 \end{equation}
 In particular, by assertion $(4)$ of Theorem
 \ref{rewritten in terms of P-exp}\;,
 $\mathrm{e}_{p^s}(\bs{f}^-(T)-\bs{f}_p^-(T),1)$ is over-convergent, and
 $\M(0,\bs{f}^-(T))\stackrel{\sim}{\to}\M(0,\bs{f}_p^{-}(T))$, over
 $\R_{H_s}$;
\item There exists a pure Witt vector
 $\bs{h}_p^-(T)\in\W_s(T^{-1}\O_{H_\infty}[T^{-1}])$ such that
 \begin{equation}
 \mathrm{e}_{p^s}(\bs{f}^-(T),1)=\mathrm{e}_{p^s}(\bs{h}_p^-(T),1)\;.
 \end{equation}
 \end{enumerate}
\end{lemma}
\emph{Proof : } Let us write $\bs{f}^-(T)=\sum_{d>0}\lb_dT^{-d}$
as a (finite) sum of $s$-co-monomials. Write $\lb_dT^{-d}=
(0,\ldots,0,\lambda_{d,0}T^{-n},\ldots,\lambda_{d,m}T^{-np^m})
\in\W_{v_p(d)}(T^{-1}\O_H[T^{-1}])$, where, for all $d>0$, we set
$d=np^m$, $m=v_p(d)$. Now set
$$\lb_{pd}^{(\Fb)}T^{-pd}:=
(0,\ldots,0,\lambda_{d,0}^{p}T^{-np},\ldots,
\lambda_{d,m}^{p}T^{-np^{m+1}})\;,$$ then the reduction
$\overline{\lb_{pd}^{(\Fb)}T^{-pd}-\lb_dT^{-d}}$ lies in
$(\Fb-1)\W_s(k(\!(t)\!))$. Hence we can replace $\lb_dT^{-d}$ with
$\lb_{pd}^{(\Fb)}T^{-pd}$. Replacing in this way
$\lb_{np^m}T^{-np^m}$ with $\lb_{np^m}^{(\Fb)}T^{-np^{m+1}}$, step
by step, we obtain a pure Witt vector. In other words, we can
suppose that for all $n\in\J$ there exists a unique $m(n)\geq 0$
such that $\lb_{np^{m(n)}}T^{-np^{m(n)}}\neq 0$. Now let us
construct $\bs{h}_p^{-}(T)$. First we arrange the sum
$\bs{f}^-(T)=\sum_{n\in\J}\sum_{m\geq 0}\lb_{np^m}T^{-np^m}$. Then
we construct, for all $n\in\J$, a natural number $m(n)\geq 0$, and
a Witt vector $\bs{\nu}_{np^{m(n)}}\in\W_s(\O_H)$, satisfying
$\mathrm{e}_{p^s}(\bs{\nu}_{np^{m(n)}}T^{-np^{m(n)}},1)=
\mathrm{e}_{p^s}(\sum_{m\geq 0}\lb_{np^m}T^{-np^m},1)$. Let
$m(n)=\sup\{m\;|\;\lb_{np^m}\neq 0\}$. By \eqref{invariant by V.},
we can suppose $s\geq m(n)$.  Let
$\lb_{np^m}=(\lambda_{np^m,0},\ldots,\lambda_{np^m,m})$, and let
$\ph{\phi_{np^m,0},\ldots,\phi_{np^m,m}}$ be its phantom vector.
Then $\mathrm{e}_{p^s}(\sum_{m=0}^{m(n)}\lb_{np^m}T^{-np^m},1)=
\exp\Bigl(\pi_{m(n)}a_0T^{-n}+\cdots+\pi_0a_{m(n)}\frac{T^{-np^{m(n)}}}{p^{m(n)}}\Bigr)$,
where, for all $j=0,\ldots,m(n)$, we have
 \begin{equation}
 a_j=\frac{\pi_0}{\pi_{m(n)-j}}\cdot \phi_{np^j,j}+\frac{\pi_1}{\pi_{m(n)-j}}
 \cdot \phi_{np^{j+1},j}+\cdots+\frac{\pi_{m(n)-j}}{\pi_{m(n)-j}}\cdot\phi_{np^{m(n)},j}\;.
 \end{equation}
Let $P(X)$ be the chosen Lubin-Tate series. Denote by
 $P^{(1)}(X):=P(X)$, $P^{(r)}(X):=P(P(\cdots P(X)\cdots))$,
 $r$-times.
 We set $h_0(X):=1$, and $h_r(X):=P^{(r)}(X)/X$, for $r=1,\ldots,m(n)$.
The phantom vector of $[h_r(\pi_{m(n)})]\in\W_{m(n)}(\O_H)$ is as
usual $\ph{h_r(\pi_{m(n)}),h_r(\pi_{m(n)-1}),\ldots,h_r(\pi_{0})}$
and is then equal to
\begin{equation}
\ph{\frac{\pi_{m(n)-r}}{\pi_{m(n)}},\frac{\pi_{m(n)-r-1}}{\pi_{m(n)-1}},
\ldots,\frac{\pi_0}{\pi_{r}},0,\ldots,0}\in\O_H^{m(n)+1}\;,\textrm{
if }r>0\;,
\end{equation}
while $[h_0(\pi_{m(n)})]=1$, and its phantom vector is
$\ph{1,\ldots,1}$.
%So $h_r(\pi_{m(n)-j})=\frac{\pi_{m(n)-r-j}}{\pi_{m(n)-j}}$ for all
%$j=0,\ldots,r$.
Hence we have
 \begin{equation*}
 a_j=h_{m(n)}(\pi_{m(n)-j})\phi_{n,j}^*+h_{m(n)-1}(\pi_{m(n)-j})\phi_{np,j}^*+\cdots+h_0(\pi_{m(n)-j})
 \phi_{np^{m(n)},j}^*\;,
 \end{equation*}
where, for all $k=0,\ldots,m(n)$,
$\ph{\phi_{np^k,0}^*,\ldots,\phi_{np^k,m(n)}^*}$ is the phantom
vector of
$\lb_{np^k}^*:=(\lambda_{np^k,0},\ldots,\lambda_{np^k,k},*,\ldots,*)\in\W_{m(n)}(\O_H)$,
where the last $m(n)-k$ components are arbitrarily chosen. Observe
that $\phi_{np^k,j}^*=\phi_{np^k,j}$, for all $j=0,\ldots,k$,
while, if $j>k$ we have $h_{m(n)-k}(\pi_{m(n)-j})=0$. This shows
that
\begin{equation*}
\bs{\nu}_{np^{m(n)}}:=[h_{m(n)}(\pi_{m(n)})]\lb_{n}^*+[h_{m(n)-1}(\pi_{m(n)})]\lb_{np}^*+\cdots+
[h_{0}(\pi_{m(n)})]\lb_{np^{m(n)}}^*\;.\Box
\end{equation*}

\subsubsection{\textbf{Proof of $(1)$ and $(2)$}}\label{first statement} The
assertions $(1)$ and $(2)$ of Theorem \ref{Introductive
foundamental Theorem} will be a direct consequence of the
following Theorem.
%In this subsection the
%hypothesis $\mathrm{w}=p$ is not necessary.
The algorithm employed is due to Robba \cite[10.10]{RoIV} (see
also \cite[13.3]{Ch-Ro}). We translate his techniques in terms of
Witt vectors. Recall that, by \ref{Reduction to a
K[T^-1]-lattice}, every rank one solvable equation has a basis in
which the matrix is a polynomial in $T^{-1}$ with coefficients in
$\O_K$.
\begin{theorem}\label{Robba's algorithm}
Let $H/K$ be a finite extension. Let $M$ be a solvable rank one
differential module over $\R_H$, defined by an operator $\d-g(T)$,
$g(T)=\sum_{-d\leq i\leq -1} a_iT^i\in\O_H[T^{-1}]$. Then there
exists a Witt vector $\bs{f}^-(T)\in\W_s(T^{-1}\O_{H'}[T^{-1}])$,
whose coefficients lies in a finite extension $H'/H$, such that
$\d-g(T)=\L(0,\bs{f}^-(T))$. More explicitly we have $\sum_{-d\leq
i\leq
-1}a_iT^i=-\sum_{j=0}^s\pi_{s-j}\sum_{i=0}^jf_i^-(T)^{p^{j-i}}\partial_{T,\log}(f_i^-(T))$
(cf. \eqref{explicit L(a_0,f^-)}). In particular
$\exp(\sum_{-d\leq
i\leq-1}a_iT^i/i)=\mathrm{e}_{p^s}(\bs{f}^-(T),1)$.
\end{theorem}
\emph{Proof : }  We shall express $\exp(\sum_{-d\leq i\leq
-1}a_iT^i/i)$ as a product of elementary $\bs{\pi}$-exponentials,
with coefficients in $H^{\mathrm{alg}}$. Observe that solvability
does not change by scalar extension of $H$. Let $d=np^m$,
$(n,p)=1$, and let $b_{d}\in H^{\mathrm{alg}}$ be such that
%\begin{equation}
$b_d^{p^m}=a_{-d}/(n\pi_0)$. %\;.
%\end{equation}
By Lemma \ref{ray in 0}, $|a_{-d}|\leq \omega<1$, hence $|b_d|\leq
1$. We consider the Witt vector
$\lb_d:=(b_d,0,\ldots,0)\in\W_m(\O_{H^{\mathrm{alg}}})$, whose
phantom vector is $\ph{b_d,b_d^p,\ldots,b_d^{p^m}}$. By
construction, we have
\begin{equation}
L_d(\lb_d) = \d + n \cdot(\pi_0 b_d^{p^m}T^{-d} + \pi_1
b_d^{p^{m-1}} T^{-d/p} + \cdots + \pi_m b_d T^{-n}).
\end{equation} Then $M\otimes\M_d(b_d,0,\ldots,0)$
is defined by an operator of the form $\d-\sum_{-d+1\leq i\leq
-1}\tilde{a}_iT^i$, $\exists\; \tilde{a}_i\in H^{\mathrm{alg}}$
(cf. \ref{tensor product}).
%Observe that $\tilde{a}_j=a_j$, for all $j\geq 0$.
Moreover $M\otimes\M_d(b_d,0,\ldots,0)$ is again solvable, so, by
\ref{ray in 0}, we have again $|\tilde{a}_{-d+1}|\leq\omega$. This
shows that we can iterate this process. More precisely there exist
$\lb_i=(b_i,0,\ldots,0)\in\W_{v_p(i)}(\O_{H^{\mathrm{alg}}})$,
$i=1,\ldots,d$, such that
\begin{equation}\label{epsilon(T) canonical}
\epsilon(T):=\prod_{i=1,\ldots,d}\mathrm{e}_i(\lb_i,T^{-1})=
\mathrm{e}_{p^s}(\sum_{i=1}^d\lb_iT^{-i},1)\;,\quad s\gg 0,
\end{equation}
satisfies $\partial_{T,\log}(\epsilon(T))=\sum_{-d\leq i\leq
-1}a_iT^i$. Then $\bs{f}^{-}(T):=\sum_{1\leq i\leq d}\lb_iT^{-i}$
(cf. \eqref{decomposition as product of P-exp}).$\Box$

\specialsection{\textbf{Applications}}
\subsection{Description of character group}
\label{DESC-CAR-GR}
\begin{lemma}\label{description more handy}\label{J}
Let \index{Jp@$\J:=\{n\in\mathbb{Z}:(n,p)=1,n>0\}$}
$\J:=\{n\in\mathbb{Z}\;|\;(n,p)=1,n>0\}$. For all fields $\kappa$
of characteristic $p$, one has the following isomorphisms of
additive groups (cf. \ref{covectors tilde}):
\begin{equation}
\bs{\mathrm{CW}}(t^{-1}\kappa
[t^{-1}])\cong\oplus_{d>0}\W_{v_p(d)}(\kappa)\qquad;\qquad\PAS(\kappa)\cong
\CW(\kappa)^{(\J)}\;,
\end{equation}
where $\CW(\kappa)^{(\J)}$ means $\oplus_{n\in\J}\CW(\kappa)$
(direct sum of copies of $\CW(\kappa)$).
\end{lemma}
\emph{Proof : } We have
 $\bs{\mathrm{CW}}(t^{-1}\kappa[t^{-1}])=\varinjlim_{s}
 \W_s(t^{-1}\kappa[t^{-1}])\stackrel{\ref{decomposition of
 CW(E)}}{=}\varinjlim_s\oplus_{d>0}\W_s^{(-d)}(\kappa)$.
 Observe that $\W_s^{(-d)}(\kappa)=\W_{\min(s,v_p(d))}(\kappa)$ (cf.
 remark \ref{W^d simto W_min(m,s)}), hence
\begin{equation}
\bs{\mathrm{CW}}(t^{-1}\kappa[t^{-1}])=
\oplus_{d>0}\varinjlim_s\W_{\min(s,v_p(d))}(\kappa)=
\oplus_{d>0}\W_{v_p(d)}(\kappa)\;.
\end{equation}
%On the other hand, we have
%$$\frac{\bs{\mathrm{CW}}(t^{-1}\kappa
%[t^{-1}])}{(\Fb-1)\bs{\mathrm{CW}}(t^{-1}
%\kappa[t^{-1}])}=\varinjlim\left(\bs{\mathrm{CW}}(t^{-1}\kappa
%[t^{-1}])\stackrel{\Fb}{\longrightarrow}
%\bs{\mathrm{CW}}(t^{-1}\kappa
%[t^{-1}])\stackrel{\Fb}{\longrightarrow}\cdots\right)$$
Now we write $d=np^m$, $n\in\J =
\{n\in\mathbb{Z}\;|\;(n,p)=1,n>0\}$ and $m\geq 0$, then on the
right hand side we have $\oplus_{d>0}\W_{v_p(d)}(\kappa)=
\oplus_{n\in\J}(\oplus_{m}\W_{v_p(np^m)}(\kappa))$. The Frobenius
morphism $\Fb$ sends $\W_{s}^{(-d)}(\kappa)$ into
$\W_s^{(-pd)}(\kappa)$, and, under the isomorphism
$\W_s^{(-d)}(\kappa)\simto\W_{\min(s,v_p(d))}(\kappa)$ (cf.
 remark \ref{W^d simto W_min(m,s)}), it becomes the morphism $\Fb\V:
\W_{v_p(np^m)}(\kappa)\to\W_{v_p(np^{m+1})}(\kappa)$ as
illustrated in the picture
\begin{center}
\begin{scriptsize}
\begin{picture}(120,60)
% ASSEX
\put(0,5){\vector(0,1){45}} \put(0,5){\vector(1,0){120}}
% COORDINATE
\put(-10,45){$m$}\put(120,0){$n$}
% Comment
\put(125,43){$d=np^m$}

% Points
\put(17.5,3.25){$\bullet$} \put(37.5,3.25){$\bullet$}
\put(57.5,3.25){$\bullet$} \put(77.5,3.25){$\bullet$}
\put(97.5,3.25){$\bullet$}
\put(17.5,22.5){$\bullet$} \put(37.5,22.5){$\bullet$}
\put(57.5,22.5){$\bullet$} \put(77.5,22.5){$\bullet$}
\put(97.5,22.5){$\bullet$}
\put(17.5,42.5){$\bullet$} \put(37.5,42.5){$\bullet$}
\put(57.5,42.5){$\bullet$} \put(77.5,42.5){$\bullet$}
\put(97.5,42.5){$\bullet$}
%
%\put(17.5,67.5){$\bullet$} \put(37.5,67.5){$\bullet$}
%\put(57.5,67.5){$\bullet$} \put(77.5,67.5){$\bullet$}
%\put(97.5,67.5){$\bullet$}

% Arrows
\put(17.5,12.5){$\uparrow$}\put(37.5,12.5){$\uparrow$}
\put(57.5,12.5){$\uparrow$}\put(77.5,12.5){$\uparrow$}
\put(97.5,12.5){$\uparrow$}
\put(17.5,32.5){$\uparrow$}\put(37.5,32.5){$\uparrow$}
\put(57.5,32.5){$\uparrow$}\put(77.5,32.5){$\uparrow$}
\put(97.5,32.5){$\uparrow$}
\put(17.5,52.5){$\uparrow$}\put(37.5,52.5){$\uparrow$}
\put(57.5,52.5){$\uparrow$}\put(77.5,52.5){$\uparrow$}
\put(97.5,52.5){$\uparrow$}
%Name of arrows
\put(102,31){\begin{tiny}$\Fb\V$\end{tiny}}
\put(102,11){\begin{tiny}$\Fb\V$\end{tiny}}
\put(102,51){\begin{tiny}$\Fb\V$\end{tiny}}
\end{picture}
\end{scriptsize}
\end{center}
Then
\begin{eqnarray}
\qquad\PAS(\kappa)&\cong&\oplus_{n\in\J}\left(\oplus_{m\geq
0}\W_{v_p(np^m)}(\kappa)/(\Fb\V-1)(\oplus_{m\geq
0}\W_{v_p(np^m)}(\kappa))\right)\\
&\cong&\left(\oplus_{m\geq
0}\W_{m}(\kappa)/(\Fb\V-1)(\oplus_{m\geq
0}\W_{m}(\kappa))\right)^{(\J)}\;.
\end{eqnarray}
One sees that $\oplus_{m\geq
0}\W_{m}(\kappa)/(\Fb\V-1)(\oplus_{m\geq 0}\W_{m}(\kappa))$ is
isomorphic to $\CW(\kappa)=\varinjlim
(\W_{m}(\kappa)\xrightarrow[]{\Fb\V}\W_{m+1}(\kappa)\xrightarrow[]{\Fb\V}\cdots)$.
\CVD

\subsection{Equations killed by an abelian extension}\label{iiii}

\subsubsection{\textbf{Extension of the field of constants}}
\label{FROB-STR}
\begin{corollary} \label{all frobenius structure} \emph{ }
The natural morphism
$$M\mapsto M\otimes K^{\mathrm{alg}}:
\mathrm{Pic}^{\mathrm{sol}}(\R_K)\to\mathrm{Pic}^{\mathrm{sol}}(\R_{K^{\mathrm{alg}}})$$
is a monomorphism. In other words, two $\R_K$-differential modules
are isomorphic if and only if they are isomorphic over
$\R_{K^{\mathrm{alg}}}$ after scalar extension.
\end{corollary}
\emph{Proof : } We show that the kernel of
$\mathrm{Pic}^{\mathrm{sol}}(\R_K)\to\mathrm{Pic}^{\mathrm{sol}}(\R_{K^{\mathrm{alg}}})$
is equal to $0$. Let $M$ be defined by the operator $L=\d-g(T)$,
$g(T):=\sum_{i}a_iT^i\in\R_K$, and suppose that $M\otimes
K^{\mathrm{alg}}$ is trivial over $\R_{K^{\mathrm{alg}}}$. By
\ref{Reduction to a K[T^-1]-lattice}, we can suppose $a_i=0$, for
all $i\neq -d,\ldots,0$. We know that $M\otimes
K^{\mathrm{alg}}\stackrel{\sim}{\to} \M(a_0,\bs{f}^{-}(T))=
\M(a_0,0)\otimes\M(0,\bs{f}^-(T))$, for a suitable
$\bs{f}^-(T)\in\W_s(T^{-1}\O_{K_s}[T^{-1}])$. Then $M\otimes
K^{\mathrm{alg}}$ is trivial only if both $\M(a_0,0)$ and
$\M(0,\bs{f}^-(T))$ are trivial over $K^{\mathrm{alg}}$. This
implies that $a_0\in\mathbb{Z}$, and hence $\M(a_0,0)$ is trivial
also over $\R_K$. On the other hand, $\M(0,\bs{f}^{-}(T))$ is
trivial if and only if $\mathrm{e}_{p^s}(\bs{f}^-(T),1)$ lies in
$\R_{K^{\mathrm{alg}}}$. By \ref{Robba's algorithm}, the series
$\mathrm{e}_{p^s}(\bs{f}^-(T),1)$ has its coefficients in $K$, and
$\M(0,\bs{f}^-(T))\in\mathrm{Pic}^{\mathrm{sol}}(\R_K)$. Since the
convergence does not change by scalar extension of $K$, it follows
that $\M(0,\bs{f}^-(T))$ is trivial over $\R_K$. $\Box$

\begin{corollary} We have $\mathrm{Pic}^{\mathrm{sol}}(\R_K)=
\mathrm{Pic}^{\mathrm{sol}}(\R_{K_\infty})^{\mathrm{Gal}(K_\infty/K)}$.
$\Box$
\end{corollary}

\subsubsection{\textbf{Frobenius structure}}
Assume now that $K$ has an absolute Frobenius $\sigma:K\to K$
(cf.\ref{what is an absolute Frobenius}), and fix an absolute
Frobenius  $\varphi:\R_K\to\R_K$. By Theorem \ref{rewritten in
terms of P-exp}$-(5)$, for any Artin-Schreier  characters
$\alpha$, the module $\M(0,\alpha)$ has a Frobenius structure of
order $1$ over $K_\infty$ (with respect to one, and hence any
absolute Frobenius, cf.\ref{Independence on varphi}). By Lemma
\ref{all frobenius structure}, this isomorphism descends to $K$.

\begin{lemma}\label{moderate caracters with frob}
$\M(a_0,0)$ has a Frobenius structure of order $h$ (cf. \ref{frob
str}) if and only if $a_0\in\mathbb{Z}_{(p)}$. Moreover let
$a_0=a/b$, $a,b\in\mathbb{Z}$, and let $b=\prod_iq_i^{r_i}$ be the
factorization of $b$ in positive prime numbers. For all
$q,r\in\mathbb{Z}$, $r>0$, we define
$[q]_{r}:=q^{q^{\cdots^{q}}}$, $r$-times, (i.e. $[q]_1=q$ and
$[q]_{r+1}=q^{[q]_{r}}$). Then
$(\varphi^*)^h(\M(a_0,0))\simto\M(a_0,0)$, with $h=\prod_i
([q_i]_{r_i}-1)$.
\end{lemma}
\emph{Proof : }By \ref{Independence on varphi} we can suppose
$\varphi=\varphi_\sigma$. Suppose that $\M(a_0,0)$ has a Frobenius
structure of order $h$. Since $\M(a_0,0)$ is solvable (cf.
\ref{frobenius implies solvability}), hence $a_0\in\mathbb{Z}_p$.
By definition \ref{frob str}, $p^h\!\cdot\!a_0-a_0\in\mathbb{Z}$,
hence $a_0\in\mathbb{Q}$. Conversely, let
$a_0=a/b\in\mathbb{Z}_{(p)}$, $b>0$. We have $p^{[q]_r-1}\equiv
1\pmod {q^r}$. Then if $h=\prod_i ([q_i]_{r_i}-1)$ we have
$(p^{h}-1)a_0\in\mathbb{Z}$.\CVD

\begin{remark} \emph{ }\label{final remarks}
Let $L=\d+\sum_{i\in\mathbb{Z}} a_i T^i$, be an operator over
$\R_K$ with Frobenius structure. The order $h$ of the Frobenius
structure depends only on the exponent $a_0\in\mathbb{Z}_{(p)}$.
Explicitly, if $a_0=a/b$, $a,b\in\mathbb{Z}$, $(b,p)=1$, and if
$b=\prod_iq_i^{r_i}>0$, $q_i>0$, is a factorization of $b$ in
prime numbers, then, by \ref{moderate caracters with frob}, we
have $h \leq \prod_i([q_i]_{r_i}-1)$.
\end{remark}

\begin{definition}
We denote by
$\mathrm{Pic}^{\mathrm{Frob}}(\R_{K_\infty})\subseteq\mathrm{Pic}^{\mathrm{sol}}(\R_{K_\infty})$
the sub-group of differential modules having a Frobenius structure
of some order $h$.
\end{definition}
\begin{corollary}
$\mathrm{Pic}^{\mathrm{Frob}}(\R_{K_\infty})\cong
\mathbb{Z}_{(p)}/\mathbb{Z}\oplus \PAS(k_\infty)$. $\Box$
\end{corollary}

\subsubsection{\textbf{Artin-Schreier extensions}} In order to apply
Theorems \ref{henselian property}, and \ref{mathcal F dag are
still series}, in this section $K$ have a discrete valuation, and
$k$ will be perfect.
\begin{proposition}[\protect{\cite[3.4]{Ma}, \cite[2.2.2]{Ts}}]
\label{mathcal F dag are still series} Let $\mathrm{F}/k(\!(t)\!)$
be a finite separable extension. Let $\mathcal{F}^{\dag}$ be the
corresponding unramified extension of $\Ed_{K,T}$. Then

$(1)$ There exist a finite unramified extension
       $\widetilde{K}/K$, a new variable $\widetilde{T}$ and an
       isometric isomorphism
$\tau:(\mathcal{F}^{\dag},|.|)\stackrel{\sim}{\to}
           (\Ed_{\widetilde{K},\widetilde{T}},
           |.|_{\widetilde{T},1})$,
       where $|.|_{\widetilde{T},1}$ is the Gauss norm with
       respect to $\widetilde{T}$. In particular, for all
       $f(T)\in\Ed_{K,T}$, one has
$           |f(T)|_{T,1}=|f(T)|_{\widetilde{T},1}$.

$(2)$ Let $\tilde{t}$ and $t$ be the reductions of
       $\widetilde{T}$ and $T$ respectively.
       Let $F=\tilde{k}(\!(\tilde{t})\!)$. Let $r$ be the
       ramification index of $\F/k(\!(t)\!)$. Write
       $t=\bar{a}_r\tilde{t}^r+\bar{a}_{r+1}
       \tilde{t}^{r+1}+\cdots$, with
       $\bar{a}_i\in \tilde{k}$. Then $\widetilde{T}$ can be
       chosen such that $\tau(T) = a_r\widetilde{T}^r +
       a_{r+1}\widetilde{T}^{r+1}+\cdots$, $a_i\in\O_{\widetilde{K}}$,
       where the $a_i$'s are liftings in $\O_{\widetilde{K}}$ of the
       $\bar{a}_i$'s.
\end{proposition}
\emph{Proof : } Let $Q(\widetilde{T}):=
a_r\widetilde{T}^r+a_{r+1}\widetilde{T}^{r+1} +\cdots$. The proof
consists in showing that $f(T)\mapsto \tau(f(T))\!:=\!
f(Q(\widetilde{T})):\Ed_{K,T}\to\Ed_{\widetilde{K},\widetilde{T}}$
is \'etale (cf. \cite[3.4]{Ma}).$\Box$
\begin{notation}
We denote by $\R_{\widetilde{K},\widetilde{T}}$ the corresponding
Robba ring.
\end{notation}
\begin{remark}
We have $(\partial_{\widetilde{T}}\circ\tau)(f(T))=
\partial_{\widetilde{T},\log}(Q(\widetilde{T}))
\cdot(\tau\circ\d)(f(T))$, where as usual
$\partial_{\widetilde{T},\log}(Q(\widetilde{T}))= \frac{\partial_{
\widetilde{T}}(Q(\widetilde{T}))}{Q(\widetilde{T})}$. Then, after
scalar extension, a generic differential operator $\d-g(T)$
becomes $\partial_{\widetilde{T}}-
\partial_{\widetilde{T},\log}(Q(\widetilde{T}))
\cdot g(Q(\widetilde{T}))$. Indeed the unique $K_\infty$
derivation of the \'etale extension $\R_{K_\infty,\widetilde{T}}$
extending $\d$ is
$\partial_{\widetilde{T},\log}(Q(\widetilde{T}))^{-1}
\cdot\partial_{\widetilde{T}}$. The solutions of this operator are
the same as those of $\d-g(T)$.
\end{remark}

\begin{corollary}\label{killed by extension}
Let $\E=k(\!(t)\!)$. Let $\F/\E$ be the Artin-Schreier  extension
defined by the kernel of $\alpha=\delta(\overline{\bs{f}}(t))$,
with $\overline{\bs{f}}(t)\in \W_s(\E)$. Let
$\R_{K,T}\to\R_{\widetilde{K},\widetilde{T}}$ be the corresponding
\'etale extension. Then the kernel of the scalar extension map
\begin{equation}\label{scalar extension morphism}
\mathrm{Res} :
\mathrm{Pic}^{\mathrm{sol}}(\R_{K_\infty,T})\longrightarrow
\mathrm{Pic}^{\mathrm{sol}}(\R_{\widetilde{K}_\infty,\widetilde{T}})
\end{equation}
is the (finite and cyclic) sub-group of
$\mathrm{Pic}^{\mathrm{sol}}(\R_{K_\infty,T})$, formed by
(isomorphism classes of) modules of the type
$\M(0,\bs{f}^-(T))^{\otimes k}$, $k\geq 0$, (cf.\ref{dec -0+}),
where $\bs{f}(T)\in\W_s(\O_{K}[[T]][T^{-1}])$ is an arbitrary
lifting of $\overline{\bs{f}}(t)$. This kernel has order
$[\F:\E]$.
\end{corollary}
\emph{Proof : } By \ref{all frobenius structure}, we can suppose
$K=K^{\mathrm{alg}}$. We decompose
$\overline{\bs{f}}(t)=\overline{\bs{f}^-}(t)+\overline{\bs{f}}_0+\overline{\bs{f}^+}(t)$
(cf. \ref{dec -0+}). Since $k=\tilde{k}$, we have
$\delta(\overline{\bs{f}}_0)=0$ (cf.
\eqref{artin-screier-diagram}). On the other hand, by
\ref{completeness of W_s(E)}, we always have
$\delta(\overline{\bs{f}^+}(t))=0$. Hence we can suppose
$\overline{\bs{f}}(t)=\overline{\bs{f}^-}(t)=(\overline{f_0^-}(t),\ldots,\overline{f_s^-}(t))$.
Since the Artin-Schreier  complex is invariant by $\V$ (cf.
\eqref{artin-screier-diagram}), we can suppose
$\overline{f_0^-}(t)\neq 0$ (i.e. the degree $[\F:\E]$ is
$p^{s+1}$). By Corollary \ref{explicit description of Pic}, the
morphism \eqref{scalar extension morphism} can be viewed as a map
\begin{equation}\label{scalar extension morphism in terms of covectors}
\mathbb{Z}_p/\mathbb{Z}\oplus
\frac{\bs{\mathrm{CW}}(t^{-1}k[t^{-1}])}{ (\Fb-1)\bs{\mathrm{CW}}(
t^{-1}k[t^{-1}])} \xrightarrow[]{\;\mathrm{Res}\;}
\mathbb{Z}_p/\mathbb{Z}\oplus
\frac{\bs{\mathrm{CW}}(\tilde{t}^{-1} k[\tilde{t}^{-1}])}{
(\Fb-1)\bs{\mathrm{CW}}(\tilde{t}^{-1} k[\tilde{t}^{-1}])}\;,
\end{equation}
where $\tilde{t}$ is the reduction of $\widetilde{T}$. We start by
studying the term $\mathbb{Z}_p/\mathbb{Z}$. By \ref{mathcal F dag
are still series}, $T = Q(\widetilde{T})$, with $Q(\widetilde{T})
= a_{p^{s+1}}\widetilde{T}^{p^{s+1}} + \cdots $, with $a_{i}\in
\O_{K_\infty}$. The differential operator $\d-a_0$,
$a_0\in\mathbb{Z}_p$ is sent to $\partial_{\widetilde{T}}-
\partial_{\widetilde{T},\log}(Q(\widetilde{T}))\cdot a_0$.
Observe that
\begin{equation}
\partial_{\widetilde{T},\log}(Q(T)) =
p^{s+1} + Q_1(\widetilde{T})\quad,\qquad
Q_1(\widetilde{T})\in\widetilde{T}\cdot\O_K[[\widetilde{T}]]\;.
\end{equation}
Hence the new operator is $\partial_{\widetilde{T}}- p^{s+1}\cdot
a_0-Q_1(\widetilde{T})\cdot a_0$. By \ref{positive queue}, this
operator is isomorphic to $\partial_{\widetilde{T}}- p^{s+1}a_0$.
Then the morphism \eqref{scalar extension morphism in terms of
covectors} sends $\mathbb{Z}_p/\mathbb{Z}$ into itself by
multiplication by $p^{s+1}=[\F:\E]$, and so is bijective on
$\mathbb{Z}_p/\mathbb{Z}$.

On the co-vectors quotient, the morphism \eqref{scalar extension
morphism in terms of covectors} is the usual functorial map
corresponding to the inclusion $t^{-1}k[t^{-1}]\longrightarrow
\tilde{t}^{-1}k[\tilde{t}^{-1}]$. The module
$\M(0,\bs{f}(T))\stackrel{\sim}{\to}\M(0,\bs{f}^-(T))$ then lies
in the kernel. Indeed, by definition of $\F/\E$, there exists
$\bs{\nu}(\widetilde{t})\in
\W_s(\widetilde{t}^{-1}k[\widetilde{t}^{-1}])$ such that (cf.
remark \ref{R(nu_0,...,nu_m)})
\begin{equation}\label{eq: 10}
\Fb(\overline{\bs{\nu}}(\widetilde{t}))-
\overline{\bs{\nu}}(\widetilde{t})= \overline{\bs{f}^-}(t),
\end{equation}
hence, by Theorem \ref{rewritten in terms of P-exp},
$\mathrm{e}_{p^s}(\bs{f}(T),1)$ lies in $\R_{K,\widetilde{T}}$. In
other words, this exponential is over-convergent in the new
variable $\widetilde{T}$. Conversely, a module $\M(0,\bs{g}^-(T))$
lies in the kernel, if and only if the exponential
$\mathrm{e}_{p^s}(\bs{g}^-(T),1)$ belongs to
$\R_{K,\widetilde{T}}$. By Theorem \ref{rewritten in terms of
P-exp}, this happens if and only if the equation
$\Fb(\bs{\nu})-\bs{\nu}=\overline{\bs{g}^-}(t)$ has a solution
$\bs{\nu}\in\W_{s}(k(\!(\tilde{t})\!))$. This happens if and only
if the kernel of $\delta(\overline{\bs{g}^-}(t))$ contains the
kernel of $\alpha=\delta(\overline{\bs{f}^-}(t))$. Since the
quotient $\G_{\E}/\mathrm{Ker}(\delta(\overline{\bs{f}^-}(t)))$ is
cyclic, this implies that $\delta(\overline{\bs{g}^-}(t))=
m\cdot\delta(\overline{\bs{f}^-}(t))$, for some $m\geq 0$. Hence
$\M(0,\bs{g}^-(T))\stackrel{\sim}{\to}\M(0,\bs{f}^{-}(T))^{\otimes
m}$.$\Box$
\begin{remark}
As suggested by the referee, this corollary is in relation with
the proposition $4.11$ of \cite{Cr}.
\end{remark}

\subsubsection{\textbf{Kummer extensions}}

\begin{corollary}
Let $\F/\E$ be an abelian totally ramified extension of degree
$[\F:\E]=n$, with $(n,p)=1$. Let
$\R_{K,T}\mapsto\R_{K,\widetilde{T}}$ be the corresponding \'etale
extension. Then the scalar extension morphism $\mathrm{Res} :
\mathrm{Pic}^{\mathrm{sol}}(\R_{K_\infty,T})\longrightarrow
\mathrm{Pic}^{\mathrm{sol}}(\R_{K_\infty,\widetilde{T}})$ is
multiplication by $n$, and so its kernel is
$(\frac{1}{n}\mathbb{Z})/\mathbb{Z}$.
\end{corollary}
\emph{Proof : } Indeed, in this case we can choose $\widetilde{t}$
satisfying $t=Q(\tilde{t})=\tilde{t}^n$.$\Box$

%               F I N O   A   Q U I   S O N O   A R R I V A T O   !  !  !  !  !  !  !  !  !  !  !  !

\subsection{A criterion of solvability}\label{OVER-CAL(E)}

This sub-section is devoted to proving the corollary \ref{crit fo
solv}. The aim of this result is to characterize the solvability
of the differential equation $\d-g(T)$, with $g(T) = \sum a_iT^i$
giving an explicit condition on the coefficients ``$a_i$''.
Roughly this Theorem shows that every solvable differential
equation over $\mathcal{E}_K$ has, \emph{without change of basis},
a solution which can be represented by the symbol
\begin{equation}
E(\bs{f}^-(T),1)\cdot T^{a_0}\cdot E(\bs{f}^+(T),1)\;,
\end{equation}
where $\bs{f}^-(T)\in\W(T^{-1}\O_K[[T^{-1}]])$ and
$\bs{f}^-(T)\in\W(T\O_K[[T]])$ are certain (\emph{infinite}) Witt
vectors, satisfying some convergence properties which ensure that
the series $E(\bs{f}^-(T),1)$ makes sense (cf. \ref{E(lb(T),1)}).
Similarly to the previous situation, this Witt vector will be a
sum of \emph{monomials} (dual notion of $s$-co-monomial, cf.
\ref{MONOMIAL}). If a Lubin-Tate group $\mathfrak{G}_P$ is chosen,
then this classification is a generalization of Theorem
\ref{Introductive foundamental Theorem}, because
$\W(T^{-1}\O_{K_\infty}[[T^{-1}]])$ contains
$\bs{\mathrm{CW}}(T^{-1}\O_{K_\infty}[[T^{-1}]])$, via the choice
of a generator $\bs{\pi}\in\mathrm{T}(\mathfrak{G}_P)$ (cf.
diagram \eqref{W(T^-1O_K[[T^-1]]) contains CW(T^-1O_K[[T^-1]])}),
and the exponential $E(\bs{f}^-(T),1)$ becomes
$\mathrm{e}_{p^s}(-,1)$ if applied to the image of a co-vector
(cf. \eqref{E(f,1)=et_p^s(pr(f),1)}).

We maintain the notations of Section \ref{Construction of Witt
vectors}. In the sequel we will work both with $T\O_K[[T]]$ and
$T^{-1}\O_K[[T^{-1}]]$. Almost all assertions have a dual meaning.

\begin{lemma}\label{E(lb(T),1)}
Let $E(-,Y):\W(\O_K[[T]])\to 1+Y\O_K[[T]][[Y]]$ be the Artin Hasse
exponential (cf. \ref{prod_j geq 0 E(lambda_j
T^p^j)=exp(phi_0T+phi_1T^p/p+cdots)}). Let $v_T$ be the $T$-adic
valuation. Let
$\bs{f}(T)=(f_0(T),f_1(T),\ldots)\in\W(T\O_K[[T]])$, and let
$\phi_j(T)$ be its $j$-th phantom component. If
$\lim_{j\to\infty}v_T(f_j(T))=+\infty$, then
$\lim_{j\to\infty}v_T(\phi_j(T))=+\infty$, and  $E(\bs{f}(T),Y)$
converges $T$-adically at $Y=1$. $\Box$
\end{lemma}
\begin{definition}
We denote by $\W^{\downarrow}(T\O_K[[T]])$
\index{Wd@$\W^{\downarrow}(T\O_K[[T]])$, $\W^{(d)}(\O_K)$} the
ideal of $\W(\O_K[[T]])$ satisfying the condition of Lemma
\ref{E(lb(T),1)}.
\end{definition}
\begin{remark}
For all
$\bs{f}^+(T)=(f_0(T),f_1(T),\ldots)\in\W^{\downarrow}(T\O_K[[T]])$,
we have
\begin{equation}\index{E(f(T),1)@$E(\bs{f}(T),1):=\prod_{j\geq 0} E(f_j(T)) =
\exp(\phi_0(T)+\frac{\phi_1(T)}{p}+\frac{\phi_2(T)}{p^2}+\cdots).
$} E(\bs{f}^+(T),1):=\prod_{j\geq 0} E(f_j(T)) =
\exp(\phi_0(T)+\frac{\phi_1(T)}{p}+\frac{\phi_2(T)}{p^2}+\cdots)\;,
\end{equation}
where $\phi_j^+(T)$ is the $j$-th phantom component of
$\bs{f}^+(T)$. The $T$-adic convergence of this product is
guaranteed by Lemma \ref{E(lb(T),1)}.
\end{remark}
\begin{definition}[Monomials]\label{MONOMIAL}
Let $\lb=(\lambda_0,\lambda_1,\ldots)\in\W(\O_K)$ and $d$ a
positive integer. We will call $\lb
T^d:=(\lambda_0T^d,\lambda_1T^{dp},
\lambda_2T^{dp^2},\ldots)\in\W^{\downarrow}(T\O_K[[T]])$ the
\emph{monomial}  of degree $d$ relative to the Witt vector
$\lb$.\footnote{Observe that if $\lb T^{-d}$ is a monomial in
$\W^{\downarrow}(T^{-1}\O_K[[T^{-1}]])$, its reduction in
$\W_m(T^{-1}\O_K[[T^{-1}]])$ is NOT a co-monomial of degree $-d$,
but it is a co-monomial of degree $-dp^m$.} In analogy with
\ref{co-monomials and PAS}, we call $\W^{(d)}(\O_K)$ the sub-group
of $\W^{\downarrow}(T\O_K[[T]])$, formed by monomials of degree
$d$.
\end{definition}
\begin{lemma}\label{prod W^(n) subset W^dag}
Let $\J:=\{n\in\mathbb{Z}\;|\;(n,p)=1, n>0\}$. There is an
injection
\begin{equation}
\prod_{n\in\J}\W^{(n)}(\O_K)\subset\W^{\downarrow}(T\O_K[[T]])\;,
\end{equation}
given by $(\lb_n T^{n})_{n\in\J}\mapsto\sum_{n\in\J}\lb_nT^{n}$.
\end{lemma}
\emph{Proof : } If $\phi_{n}=(\phi_{n,0},\phi_{n,1},\ldots)$ is
the phantom vector of $\lb_n$, then the phantom vector of
$\lb_{n}T^n$ is
$(\phi_{n,0}T^n,\phi_{n,1}T^{np},\phi_{n,2}T^{np^2},\cdots)$.
Hence all terms have different degree and they do not ``blend''
when we sum the phantom components. $\Box$

\begin{remark}
%\begin{enumerate}
%\item
$\bs{1.}$ Let $\bs{f}(T)\in \W^{\downarrow}(T\O_K[[T]])$, %(resp. $\bs{f}(T)\in\W^{\downarrow}(\T^{-1}\O_K[[T^{-1}]])$)
let $\lb,\lb_d\in\W(\O_K)$, $d>0$. Then we have
\begin{eqnarray}\label{remarkable invariant under verschiebung}
 E(\V(\bs{f}(T)),1)&=&E(\bs{f}(T),1)\;,\\
 E(\lb,T^d)&=&E(\lb T^d,1)\;,\\
 \prod_{d\geq 1}E(\lb_d,T^d)&=&E\bigl(\sum_{d\geq 1}\lb_d T^d,1\bigr)\;.
\end{eqnarray}
%\item
$\bs{2.}$ If $\phi_{-n}=(\phi_{-n,0},\phi_{-n,1},\cdots)$ is the
phantom vector of $\lb_{-n}$, then we have
$$E(\sum_{n\in\J}\lb_{-n}T^{-n},1)=\exp(\sum_{n\in\J}\sum_{m\geq
0}\phi_{-n,m}\frac{T^{-np^m}}{p^m})$$
%\item
$\bs{3.}$ If
$\bs{f}^-(T)=(f_0^-(T),f_1^-(T),\ldots)\in\W^{\downarrow}(T^{-1}\O_K[[T^{-1}]])$
and if $\mathrm{pr}_m(\bs{f}^-(T))$
\index{pr_m@$\mathrm{pr}_m:\W(-)\to\W_m(-)$}is the image of
$\bs{f}^-(T)$ in $\W_m(T^{-1}\O_K[[T^{-1}]])$, then (cf.
\eqref{definition of et_d(lb,T)})
\begin{equation}\label{E(f,1)=et_p^s(pr(f),1)}
 E([\pi_m]\cdot\bs{f}^-(T),1)=
 \mathrm{e}_{p^m}\bigl(\mathrm{pr}_m(\bs{f}^-(T)),1\bigr)=
 \prod_{j\geq 0}^mE_{m-j}(f_j^-(T)).
\end{equation}
The exponentials used in the preceding section are then a
particular case of $E(-,1)$.
%\end{enumerate}
\end{remark}
\begin{remark}
Recall that $\W_m(T^{-1}\O_{K_m}[T^{-1}])
\stackrel{\sim}{\to}[\pi_m]\W(T^{-1}\O_{K_m}[T^{-1}])
\subset\W(T^{-1}\O_{K_m}[[T^{-1}]])$ (see \eqref{W_m simto
[pi_m]W}). We have the following commutative diagram
\begin{equation}\label{W(T^-1O_K[[T^-1]]) contains CW(T^-1O_K[[T^-1]])}
 \xymatrix{\W_{m+1}(T^{-1}\O_{K_{m+1}}[T^{-1}])
 \ar@{^(->}[r]_{[\pi_{m+1}]\cdot}
 &
 \W^{\downarrow}(T^{-1}\O_{K_{m+1}}[[T^{-1}]])\ar[r]^-{E(-,1)}
 &1+T^{-1}\O_{K_\infty}[[T^{-1}]]\\
 \W_{m}(T^{-1}\O_{K_m}[T^{-1}])\ar[u]^{\V}
 \ar@{^(->}[r]_{[\pi_m]\cdot}&
 \W^{\downarrow}(T^{-1}\O_{K_m}[[T^{-1}]])\ar[u]^{\V}\ar[ur]_-{E(-,1)}& \;.}
\end{equation}
Indeed, we see, looking at the phantom components, that
\begin{equation}
[\pi_m](f_0^-,\ldots,f_{m}^-,f_{m+1}^-,\cdots)
=[\pi_m](f_0^-,\ldots,f_{m}^-,0,0,\cdots)\;,
\end{equation}
for all
$\bs{f}^-(T)=(f_0^-,\ldots,f_{m}^-,f_{m+1}^-,\cdots)\in\W(T^{-1}\O_K[T^{-1}])$.
Hence $[\pi_m]\bs{f}^-(T)$ lies in
$\W^{\downarrow}(T^{-1}\O_K[T^{-1}])$, for all
$\bs{f}^-(T)\in\W(T^{-1}\O_K[T^{-1}])$.
\end{remark}
\begin{remark}\label{we can reduce the study to monomials}
By definition one has
$\W^{(d)}(\O_K)\subset\W^{\downarrow}(T\O_K[[T]])$, for all $d\geq
1$. The group $\W^{\downarrow}(T\O_K[[T]])$ is not generated by
the family $\{\W^{(d)}(\O_K)\}_{d\geq 0}$ of sub-groups. Indeed,
for example, the $m$-th phantom component $\phi_m(T)$ of a Witt
vector of the form $\sum_{d>0}\lb_d T^d$ is always of the type
$\phi_m(T)=h(T^{p^m})$, for some $h(T)\in\O_K[[T]]$.

However, the basic fact is that, for all
$\bs{f}^+(T)\in\W^{\downarrow}(T\O_K[[T]])$ (resp.
$\bs{f}^-(T)\in\W^{\downarrow}(T^{-1}\O_K[[T^{-1}]])$), there
exists an (infinite) family of monomials $\{\lb_{n}
T^{n}\}_{n\in\J}\in\prod_{n\in\J}\W^{(n)}(\O_K)$ (resp.
$\{\lb_{-n} T^{-n}\}_{n\in\J}\in\prod_{n\in\J}\W^{(-n)}(\O_K)$)
satisfying
\begin{equation*}
E(\bs{f}^+(T),1)=E(\sum_{n\in\J}\lb_{n}T^{n},1)\quad;\quad
E(\bs{f}^-(T),1)=E(\sum_{n\in\J}\lb_{-n}T^{-n},1).
\end{equation*}
In other words, a general Witt vector is not an infinite sum of
monomials, but the Artin-Hasse exponential of this Witt vector is
always equal to the Artin-Hasse exponential of an infinite sum of
monomials with support in $\J$.
%This follows essentially from
%\eqref{remarkable invariant under verschiebung}.
\end{remark}

\begin{lemma}\label{criteria of solvability lemma}
The differential equation $\d-g^+(T)$, $g^+(T)=\sum_{i\geq
1}a_iT^i\in\R_K$ is solvable if and only if there exists a family
$\{\lb_{n}\}_{n\in\J}$, $\lb_n\in\W(\O_K)$, with phantom
components $\phi_{n}=(\phi_{n,0},\phi_{n,1},\ldots)$ satisfying
\begin{equation}\label{a_np^m=n phi_n,m}
a_{np^m}=n\phi_{n,m}\;,\qquad\textrm{ for all }n\in\J,\; m\geq
0\;.
\end{equation}
In other words we have $\exp(\sum_{i\geq 1}a_i
\frac{T^i}{i})=E(\sum_{n\in\J}\lb_nT^n,1)$.
\end{lemma} \dem The formal series
$E(\sum_{n\in\J}\lb_nT^{n},1)\in 1+T\O_K[[T]]$ is a solution of
the equation $L:=\d-\sum_{n\in\J}\sum_{m\geq
0}n\phi_{n,m}T^{np^m}$. Since this exponential converges in the
unit disk, then $Ray(L,\rho)=\rho$, for all $\rho<1$ and $L$ is
solvable. Conversely if $\d-g^+(T)$ is solvable, then the Witt
vectors $\lb_n=(\lambda_{n,0},\lambda_{n,1},\ldots)$ are defined
by the relation \eqref{a_np^m=n phi_n,m} (cf. \ref{phi is inj}).
For example for all $n\in\J$ we have
\begin{equation}
\lambda_{n,0} = \frac{a_{n}}{n} \;\;,\qquad \lambda_{n,1} =
\frac{1}{p}\left(\frac{a_{np}}{n} - \bigl(\frac{a_n}{n}\bigr)^p
\right)\;.
\end{equation}
We must show that $|\lambda_{n,m}|\leq 1$ for all $n\in\J$, $m\geq
0$.\\
\emph{Step 1:} By the Small Radius Lemma \ref{small radius}, we
have $|a_i|\leq 1$, for all $i\geq 1$. Hence, for all $n\in\J$, we
have $|\lambda_{n,0}|\leq 1$. Then the exponential
$$E(\sum_{n\in\J}(\lambda_{n,0},0,0,\ldots)T^{n},1)=
\exp(\sum_{n\in\J}\sum_{m\geq
0}\lambda_{n,0}^{p^m}\frac{T^{p^m}}{p^m})$$ converges in the unit
disk and is a solution of the operator $Q^{(0)}:=\d - h^{(0)}(T)$,
with $h^{(0)}(T)=\sum_{n\in\J}\sum_{m\geq
0}\lambda_{n,0}^{p^m}T^{p^m}$, which is therefore solvable.\\
\emph{Step 2:} The tensor product operator $\d - (g^+(T) -
h^{(0)}(T))$ is again solvable and satisfies
$g^+(T)-h^{(0)}(T)=p\cdot g^{(1)}(T^p)$, for some $g^{(1)}(T)\in T
K[[T]]$. In other words the ``antecedent by ramification''
$\varphi_p^*$ (cf. \ref{varphi_p^*}) of the equation $\d - (g^+(T)
- h^{(0)}(T))$ is given by $\d-g^{(1)}(T)$,
which is therefore solvable.\\
\emph{Step 3: } We observe that $g^{(1)}(T) =\frac{1}{p}
\sum_{n\in\J}\sum_{m\geq
0}(a_{np^{m+1}}\!-n(\frac{a_n}{n})^{p^{m+1}}) T^{np^{m}}$, and
again by the Small Radius Lemma we have
$|a_{np}\!-n(\frac{a_n}{n})^{p}|\leq 1$, which implies
$|\lambda_{n,1}|\leq 1$. The process can be iterated indefinitely.
\CVD

\begin{remark}\label{discussion}We shall now consider the general
case of an equation $\d-g(T)$, with
$g(T)=\sum_{i\in\mathbb{Z}}a_iT^i\in\R_K$, and get a criterion of
solvability. Suppose that $\d-g(T)$ is solvable. We know that
$\d-g^-(T)$, $\d-a_0$ and $\d-g^+(T)$ are all solvable (cf.
\ref{division of the problem}). We can then consider $\d-g^-(T)$
as an operator on $]1,\infty]$ (instead of
$]1-\varepsilon,\infty]$) and Lemma \ref{criteria of solvability
lemma} gives us the existence of a family of Witt vector
$\{\lb_{-n}\}_{n\in\J}$ satisfying $a_{-np^m} = -n\phi_{-n,m}$,
for all $n\in\J$, and all $m\geq 0$. Conversely suppose that we
are given two families $\{\lb_{-n}\}_{n\in\J}$ and
$\{\lb_{n}\}_{n\in\J}$, with $\lb_n\in\W(\O_K)$. Since the phantom
components of $\lb_n$ are bounded by $1$, then $g^+(T)$ is bounded
and belongs to $\R_K$. Now we need a condition on the family
$\{\lb_{-n}\}_{n\in\J}$ in order that the series
$g^-(T):=\sum_{n\in\J}\sum_{m\geq 0}-n \phi_{-n,m}T^{-np^m}$
belongs to $\R_K$.
\end{remark}

\begin{lemma}\label{phi_i/p^i leq c rho^np^i <==> lambda_i leq c rho^np^i}
Let $c\leq\omega=|p|^{\frac{1}{p-1}}$, $n\in\J$, $\rho\leq 1$ be
fixed. Let $(\lambda_0,\lambda_1,\ldots)\in\W(\O_K)$ and let
$\phi=(\phi_{0},\phi_1,\ldots)$ be its phantom vector. Then
$|\phi_i/p^i|\leq c\rho^{np^i}$ for all $i\geq 0$ if and only if
$|\lambda_i|\leq c\rho^{np^i}$ for all $i\geq 0$.
\end{lemma}
\dem Recall that $c^{p^i}\leq |p|^ic$, for all $i\geq 0$. Suppose
that $|\phi_i/p^i|\leq c\rho^{np^i}$ for all $i\geq 0$. Then
$|\lambda_0|=|\phi_0|\leq c\rho^n$. By induction suppose that
$|\lambda_j|\leq c\rho^{np^j}$ for all $j=0,\ldots,i-1$, then
$|\lambda_i|=|\frac{1}{p^i}(\phi_i-\lambda_0^{p^i}-p\lambda_{1}^{p^{i-1}}-\ldots-p^{i-1}\lambda_{i-1}^p)|$.
By induction $|\phi_{i}|\leq |p|^ic\rho^{np^i}$ and
$|p^k\lambda_{k}^{p^{i-k}}|\leq
|p|^k(c\rho^{np^{k}})^{p^{i-k}}=|p|^k
c^{p^{i-k}}\rho^{np^i}\leq|p^i|c\rho^{np^i}$, hence
$|\lambda_i|\leq c\rho^{np^i}$. Conversely suppose that
$|\lambda_i|\leq c\rho^{np^i}$ for all $i\geq 0$. Then
$|\phi_i|=|\lambda_0^{p^i}+p\lambda_1^{p^{i-1}}+\cdots+p^i\lambda_i|\leq
\sup((c\rho^n)^{p^i},|p|(c\rho^{np})^{p^{i-1}},\cdots,|p|^i(c\rho^{np^i}))\leq
|p|^ic\rho^{np^i}$.\CVD

\begin{definition}
Let $c\leq \omega$ and $\rho\leq 1$. We denote by
$$\W^{\downarrow}_{c,\rho}(T^{-1}\O_K[[T^{-1}]])\subset
\prod_{n\in\J}\W^{(-n)}(\O_K)\stackrel{\ref{prod W^(n) subset
W^dag}}{\subset} \W^{\downarrow}(T^{-1}\O_K[[T^{-1}]])$$ the
sub-group formed by the sums $\sum_{n\in\J}\lb_{-n}T^{-n}$ such
that $\lb_{-n}=(\lambda_{-n,0},\lambda_{-n,1},\ldots)\in\W(\O_K)$
verify the conditions of Lemma \ref{phi_i/p^i leq c rho^np^i <==>
lambda_i leq c rho^np^i}. In other words $|\lambda_{-n,m}|\leq
c\rho^{np^m}$.
\end{definition}

\begin{remark}
Observe that, by Lemma \ref{phi_i/p^i leq c rho^np^i <==> lambda_i
leq c rho^np^i}, a Witt vector $\sum_{n\in\J}\lb_{-n}T^{-n}$
belongs to the subgroup
$\W^{\downarrow}_{c,\rho}(T^{-1}\O_K[[T^{-1}]])$ if and only if
the argument of the exponential
$E(\sum_{n\in\J}\lb_{-n}T^{-n},1)=\exp(\sum_{n\in\J}\sum_{m\geq
0}\phi_{-n,m}\frac{T^{-np^m}}{p^m})$ satisfies
\begin{equation}\label{small argument of the exponential}
\left|\sum_{n\in\J}\sum_{m\geq
0}\phi_{-n,m}\frac{T^{-np^m}}{p^m}\right|_\rho:=\sup_{n\in\J,m\geq
0}\Bigl(\frac{|\phi_{-n,m}|}{|p|^m}\rho^{-np^m}\Bigr)\leq c\leq
\omega\;.
\end{equation}
\end{remark}
\begin{definition}
Let
$\W^{\dag}(T^{-1}\O_K[[T^{-1}]])\subset\W^{\downarrow}(T^{-1}\O_K[[T^{-1}]])$
be the subgroup defined as the sum of the sub-group
$\bigcup_{c<\omega,\rho<1}\W_{c,\rho}^{\downarrow}(T^{-1}\O_K[[T^{-1}]])$
with the sub-group $\Bigl(\bigcup_{j\geq
0}[\pi_j]\cdot\W^{\downarrow}(T^{-1}\O_{K^{\mathrm{alg}}}[T^{-1}])\Bigr)
\bigcap \W^{\downarrow}(T^{-1}\O_K[T^{-1}])\:.$
\end{definition}
\begin{corollary}[Solvability Criterion]
\label{crit fo solv} Let $\d-g(T)$,
$g(T):=\sum_{i\in\mathbb{Z}}a_iT^i\in\R_K$ be a solvable equation.
Then $a_0\in\mathbb{Z}_p$ and there exist two families
$\{\lb_{-n}\}_{n\in\J}$ and $\{\lb_{n}\}_{n\in\J}$ such that for
all $n\in\J$, and all $m\geq 0$ we have
$$a_{-np^m}=-n\phi_{-n,m}\quad;\qquad a_{np^m}=n\phi_{n,m}\;,$$
where $(\phi_{-n,0},\phi_{-n,1},\ldots)$ (resp.
$(\phi_{n,0},\phi_{n,1},\ldots)$) is the phantom vector of
$\lb_{-n}$ (resp. $\lb_n$). Moreover $\sum_{n\in\J}\lb_{-n}T^{-n}$
belongs to $\W^{\dag}(T^{-1}\O_K[[T^{-1}]])$.

Conversely given a triplet
$(\sum_{n\in\J}\lb_{-n}T^{-n},a_0,\sum_{n\in\J}\lb_{n}T^{n})$,
with
$\sum_{n\in\J}\lb_{-n}T^{-n}\in\W^{\dag}(T^{-1}\O_K[[T^{-1}]])$,
$a_0\in\mathbb{Z}_p,\sum_{n\in\J}\lb_{n}T^{n}\in\W(T\O_K[[T]])$,
then
$$g(T):=\sum_{n\in\J}\sum_{m\geq
0}-n\phi_{-n,m}T^{-np^m}+a_0+\sum_{n\in\J}\sum_{m\geq
0}n\phi_{n,m}T^{np^m}$$ belongs to $\R_K$, and the equation
$\d-g(T)$ is solvable.
\end{corollary}
\begin{remark}\label{remark on the pure witt of finite queue}
This corollary asserts that $\sum_{n\in\J}\lb_{-n}T^{-n}$ is a sum
of a ``small'' vector, i.e. verifying the relation \eqref{small
argument of the exponential}, and a vector of ``type Robba'', i.e.
of the type $[\pi_j]\bs{f}^-(T)$, $\bs{f}^-(T) \in
\W(T^{-1}\O_{K^{\mathrm{alg}}}[[T^{-1}]])$, for some $j\geq 0$ and
such that the product $[\pi_j]\bs{f}^-(T)$ lies in
$\W(T^{-1}\O_{K}[[T^{-1}]])$ i.e. has its coefficients in $K$.
Actually the proof will show that $\bs{f}^-(T)$ can be chosen pure
(see below).
\end{remark}
\emph{Proof of \ref{crit fo solv}: } Let $\d-g(T)$ be solvable. By
\ref{criteria of solvability lemma}, we know the existence of
$\{\lb_{-n}\}_{n\in\J}$ and $\{\lb_{n}\}_{n\in\J}$ (cf.
\ref{discussion}). We must show that $\sum_{n\in\J}\lb_{-n}T^{-n}$
lies in $\W^{\dag}(T^{-1}\O_K[[T^{-1}]])$. Let $d>0$ be such that
$|\sum_{i<-d}a_iT^i/i|_\rho<\omega$ for some $\rho<1$ (cf.
\ref{algebricity}). Write $g^-(T)=\sum_{i<-d}a_iT^i + \sum_{-d\leq
i\leq -1}a_iT^i$. By \ref{algebricity} we know that
$\exp(\sum_{i<-d}a_iT^i/i)\in\R_K$, hence the equation
$\d-\sum_{i<-d}a_iT^i$ is solvable (and actually trivial). In
particular $\d-\sum_{-d\leq i\leq -1}a_iT^i$ is solvable and
hence, again by \ref{discussion}, there exists a family
$\{\lb_{-n}'\}_{n\in\J}$, such that
$$-n\phi_{-n,m}'=\left\{\begin{array}{rcl}
a_{-np^m}&\textrm{if}& \phantom{-d\leq}-np^m<-d\\\
0&\textrm{if} & -d\leq -np^m\leq -1
\end{array}\right.$$
where $(\phi_{-n,0},\phi_{-n,1},\cdots)$ is the phantom vector of
$\lb_{-n}$. Since, by construction
$|\sum_{i<-d}a_iT^i/i|_\rho<\omega$, this implies
$|\sum_{n\in\J}\sum_{m\geq
0}\phi_{-n,m}T^{-np^m}/p^m|_\rho<\omega$, hence
$\sum_{n\in\J}\lb_{-n}'T^{-n}$ lies in
$\W_{c,\rho}^{\downarrow}(T^{-1}\O_K[[T^{-1}]])$ for some
$c<\omega$. Now we consider $\lb_{-n}'' := \lb_n - \lb_n'$, the
family $\{\lb_{-n}''\}_{n\in\J}$ then satisfies
$$-n\phi_{-n,m}''=\left\{\begin{array}{rcl}
0&\textrm{if}& \phantom{-d\leq}-np^m<-d\\\
a_{-np^m}&\textrm{if} & -d\leq -np^m\leq -1
\end{array}\right.$$
By \ref{Robba's algorithm}, and by \ref{isomorphic to a pure
vector} there exists a \emph{pure} Witt vector
$\bs{f}^-(T)=(f_0^-(T),\ldots,f_s^-(T))\in
\W_s(T^{-1}\O_{K^{\mathrm{alg}}}[T^{-1}])$ such that
$\L(0,\bs{f}^-(T))=\d-\sum_{-d\leq i\leq -1}a_iT^i$. Hence
$[\pi_s]\bs{f}^-(T)$ and $\sum_{n\in\J}\lb_{-n}''T^{-n}$ have the
``same'' phantom vector because $\bs{f}^-(T)$ is pure. Then
$\sum_{n\in\J}\lb_n''T^{-n}$ lies in the image of the morphism
$\W_s(T^{-1}\O_{K^{\mathrm{alg}}}[T^{-1}]) \stackrel{\sim}{\to}
[\pi_s]\cdot\W(T^{-1}\O_{K^{\mathrm{alg}}}[T^{-1}])$. \CVD
%\begin{corollary}
%Let $\d-g(T)$, $g(T)=\sum_{-d\leq -1}a_iT^i$ be a solvable
%equation. Then there exists a \emph{pure} Witt vector
%$\bs{f}^-_p(T)\in\W_s(\O_{K^{\mathrm{alg}}})$
%\end{corollary}
%\dem It follow from the proof of \ref{crit fo solv}.\CVD
\begin{remark}\label{g_|n}
Let $L:=\d-\sum_{i\in\mathbb{Z}}a_iT^i$, $g(T)\in\R_K$ be a given
equation. Then $L$ is solvable if and only if
$a_0\in\mathbb{Z}_p$, and, for all $n\in\J$, both the operators
$\d-\sum_{m\geq 0}a_{np^m}T^{np^m}$ and $\d-\sum_{m\geq
0}a_{-np^m}T^{-np^m}$ are solvable.
\end{remark}

\begin{corollary}\label{no equation over abs unr}
If $K$ is unramified over $\mathbb{Q}_p$, then every solvable
differential module over $\R_K$ of rank one is isomorphic to a
moderate module (cf. \ref{moderate characters}). In other words,
\begin{equation}
\mathrm{Pic}^{\mathrm{sol}}(\R_K)=\left\{
\begin{array}{lcl}
\mathbb{Z}_p/\mathbb{Z}&\textrm{ if }&p>
2\\
&&\\
\mathbb{Z}_p/\mathbb{Z}\oplus
k(\!(t)\!)/(\Fb-1)k(\!(t)\!)&\textrm{ if }&p=2
\end{array}\right.
\end{equation}
\end{corollary}
\emph{Proof : } We must show that all $\bs{\pi}$-exponential
$\mathrm{e}_{p^s}(\bs{f}^-(T),1)$ whose logarithmic derivative has
its coefficients in $K$ is trivial. Actually we can suppose that
the co-monomial $\bs{f}^-(T)$ is \emph{pure} (cf. \ref{isomorphic
to a pure vector}). Write
$\partial_{T,\log}(\mathrm{e}_{p^s}(\bs{f}^-(T),1))=\sum_{-d\leq
i\leq -1}a_iT^i$ with $a_i\in\O_K$ for all $i=-d,\ldots,-1$. Write
$\bs{f}^-(T) = \sum_{n\in\J}\lb_{-np^{m(n)}}T^{-np^{m(n)}}$,
$\lb_{-np^{m(n)}}=(\lambda_{-np^{m(n)},0},\ldots,\lambda_{-np^{m(n)},m(n)})
\in\W_{m(n)}(\O_{K^{\mathrm{alg}}})$. Since $\bs{f}^-(T)$ is pure,
one has (cf. \ref{remark on pure Witt vector})
\begin{equation}\label{first equality}
a_{-np^j}=-n\pi_{m(n)-j}\phi_{-np^{m(n)},j}\;,
\end{equation}
for all $j=0,\ldots,m(n)$, where
$\ph{\phi_{-np^{m(n)},0},\ldots,\phi_{-np^{m(n)},m(n)}}$ is the
phantom vector of $\lb_{-np^{m(n)}}$. On the other hand the
Criterion of Solvability \ref{crit fo solv} asserts the existence
of a family $\{\lb_{-n}'\}_{n\in\J}$ with phantom vector
$\{\bs{\phi}_{-n}'\}_{n\in\J}$, with
$\bs{\phi}_{-n}':=\ph{\phi_{-n,0},\phi_{-n,1},\ldots}$, such that
$a_{-np^m}=-n\phi_{-n,m}'$ for all $n\in\J$, $m\geq 0$. Observe
that $\phi_{-n,m}'\in \O_K$. Since $K$ is unramified over
$\mathbb{Q}_p$, then we can employ Lemma \ref{congruece on the
phantom components}. Then
$\phi_{-n,m}'\equiv\sigma(\phi_{-n,m-1}')\mod p^{m}\O_K$, for all
$n\in\J$, $m\geq 0$, that is
\begin{equation}
a_{-np^{j}}\equiv \sigma(a_{-np^{j-1}})\mod
p^j\O_K\quad\textrm{for all }j\geq 0.
\end{equation}
Since $a_{-np^{m(n)+1}}=0$ we obtain, by \eqref{first equality},
the estimate $|\pi_{m(n)-j}\phi_{-np^{m(n)},j}|\leq|p|^{j+1}$, for
all $j=0,\ldots,m(n)$. Then we have a system of conditions
$$|\lambda_{-np^{m(n)},0}^{p^j}+p\lambda_{-np^{m(n)},1}^{p^{j-1}}+
\cdots+p^j\lambda_{-np^{m(n)},j}|\leq
|p|^{j+1}|\pi_{m(n)-j}|^{-1}\;,$$ which easily gives
$|\lambda_{-np^{m(n)},j}|\leq |p|^j|\pi_{m-j}|^{-1}<1$. This last
is $\leq 1$, and is $=1$ if and only if $p=2$, and $m(n)=j=0$. If
$p\neq 2$, then by \ref{reduction to k}, and corollary \ref{all
frobenius structure}, $\mathrm{e}_{d}(\lb_{-np^{m(n)}},T^{-1})$
lies in $\R_K$, for all $n\in\J$, and $\L_d(0,\bs{f}^-(T))$ is
trivial.\CVD

\subsection{Explicit computation of the Irregularity in some cases}
\label{COMP-IRR} Let $v_t$ be the $t$-adic valuation of
$k(\!(t)\!)$.
\begin{lemma}\label{irreg of (0...0f0...0) with (n,p)=1}
Let $f(T)\in T^{-1}\O_K[T^{-1}]$ be a polynomial in $T^{-1}$, and
let $\overline{f}(t)\in t^{-1}k[t^{-1}]$ be the reduction of
$f(T)$. Let $n:=-v_t(\overline{f}(t))>0$. If $(n,p)=1$, then
\begin{equation}\label{shiffttt}
 \mathrm{Irr}\Bigl( \M\bigl(0,(0,\ldots,0,
 \underbrace{f(T),0,\ldots,0}_{\ell+1}\;)\bigr)
 \Bigr)= n\cdot p^{\ell}\;,
\end{equation}
where $\ell=\ell(0,\ldots,0,f(T),0,\ldots,0)$ (cf. \ref{ell}).
\end{lemma}
\emph{Proof : } We have
$\M\bigl(0,(0,\ldots,0,f(T),0,\ldots,0)\bigr) =
\M(0,(f(T),0,\ldots,0))\;.$ (cf. \eqref{invariant by V.}).
Moreover, the isomorphism class of this module depends only on
$\overline{f}(t)$, hence we can suppose that
$f(T)=a_{-n}T^{-n}+\cdots+a_{-1}T^{-1}$, with $|a_{-n}|=1$. Then
 %\begin{equation*}
$\L(0,(f(T),0,\ldots,0))=\d +
 \partial_{T,\log}(f(T))\cdot\left[\pi_sf(T)+\pi_{s-1}
 f(T)^p+\cdots+\pi_{0}f(T)^{p^s}\right]$.
 %\end{equation*}
We have $\partial_{T,\log}(f(T))=-n+TQ(T)$, with
$Q(T)\in\O_K[[T]]$, so that
\begin{equation}
g(T)=-\pi_0\cdot n\cdot a_{-n}^{p^\ell}\cdot T^{-n\cdot
p^\ell}+(\textrm{terms of degree} > -np^\ell)\;.
\end{equation}
Since $(n,p)=1$, we can apply \ref{graphic} and
$\mathrm{Irr}(\d+g(T))=\mathrm{Irr}_F(\d+g(T))=np^\ell$. $\Box$

\begin{corollary}\label{garuti} Let
$\overline{\bs{f}^-}(t)=
(\overline{f^-_0},\ldots,\overline{f^-_s})
\in\W_s(t^{-1}k[t^{-1}])$. Let $n_j:=
-v_t\left(\overline{f^-_j}\right)$. If $(n_j,p)=1$, or $n_j=0$,
for all $j=0,\ldots,s$ (cf. \ref{J}), then
\begin{equation}
\mathrm{Irr}\left(\M(0,\overline{\bs{f}^-}(t))\right)=\max_{0\leq
j\leq s}\left(n_j\cdot p^{s-j}\right).
\end{equation}
\end{corollary}
\emph{Proof : } Let $M_j$ be the differential module defined by
$(0,\ldots,0,\overline{f^-_j}(t),0,\ldots,0)$. By \ref{irreg of
(0...0f0...0) with (n,p)=1}, $\mathrm{Irr}(M_j)=n_jp^{s-j}$. Since
$\M(0,\bs{f}^{-}(T))=\otimes_jM_j$ (cf. \eqref{rule for
induction}), and since $n_jp^{s-j}$ are all different, then, by
\ref{min tensor radius}, we have the desired conclusion. $\Box$

\subsection{Tannakian group}\label{tannakian group section} In this
section we study the category of solvable differential modules
over $\R_K$ which are extensions of rank one sub objects. We
remove the hypothesis ``K is spherically complete'', present in
the literature. Let $H/K$ be an arbitrary algebraic extension. We
set
\index{Hdag@$\mathcal{H}^{\dag}_H:=\cup_\varepsilon\;\a_H(]1-\varepsilon,\infty[)$}
$\mathcal{H}^{\dag}_H:=\cup_\varepsilon\;\a_H(]1-\varepsilon,\infty[)$.
Let $S$ \index{S@$S\subseteq\mathbb{Z}_p$} be a sub-group of
$\mathbb{Z}_p$ without Liouville numbers and containing
$\mathbb{Z}$.
\begin{definition}
Let $\mathcal{C}$ be an additive category. If there exists a
function \emph{rank} on $\mathcal{C}$, then we denote by
\index{MLSplus@$\MLCS_{\oplus\textrm{-}1}(\mathcal{H}^{\dag}_L)$,
$\MLCS_{\mathrm{ext}\textrm{-}1}(\mathcal{H}^{\dag}_L)$,
$\MLCS_{\mathrm{ext}\textrm{-}1}^{\mathrm{reg}}(\mathcal{H}^{\dag}_L,S)$}
$\mathcal{C}_{\oplus\textrm{-}1}$ (resp.
$\mathcal{C}_{\mathrm{ext}\textrm{-}1}$) the full sub-category of
$\mathcal{C}$ whose objects are finite direct sum (resp. finite
successive extension) of rank one objects.
\end{definition}
\begin{definition}
An object is said to be \emph{simple} if it has no non trivial
sub-objects. It is said \emph{indecomposable} if it is not a
direct sum of non trivial objects.
\end{definition}
\begin{definition}
Let $\MLCS(\mathcal{H}^{\dag}_H)$ be the category of (free)
differential modules over $\mathcal{H}^{\dag}_H$ solvable at $1$
(i.e. $Ray(N,1)=1$, cf. \ref{solvability}). Recall that, by
definition, such a module comes, by scalar extension, from a
module over $\mathcal{H}^{\dag}_L$, for some finite extension
$L/K$ (cf. \eqref{a_H=a_K otimes H}). Let
$N\in\MLCS_{\mathrm{ext}\textrm{-}1}(\mathcal{H}^{\dag}_H)$ be
extension of rank one modules, say $\{N_i\}_{i=1,\ldots,k}$. We
will say that $N$ is regular at $\infty$, write
$N\in\MLCS_{\mathrm{ext}\textrm{-}1}^{\mathrm{reg}}(\mathcal{H}^{\dag}_H)$,
if, for all $i$, the module $N_i$ is defined, in some basis, by an
operator $\d+g_i(T)$, satisfying
\begin{equation}
g_i(T)=\sum a_{i,j}T^{j}\;,\quad \textrm{ with } a_{i,j}=0\;,
\textrm{ for all }j \geq 1\;.
\end{equation}
We will say that
$N\in\MLCS_{\mathrm{ext}\textrm{-}1}^{\mathrm{reg}}(\mathcal{H}^{\dag}_H,S)$,
if
$N\in\MLCS_{\mathrm{ext}\textrm{-}1}^{\mathrm{reg}}(\mathcal{H}^{\dag}_H)$,
and if $a_{i,0}\in S,\forall i$.
\end{definition}
\begin{lemma}[Schur's lemma] \label{Schur}
Let $M_1,M_2$ be two rank one objects in $\MLCS(\R_{H})$ (resp.
$\MLCS_{\mathrm{ext}\textrm{-}1}^{\mathrm{reg}}(\mathcal{H}^{\dag}_H)$).
Every non zero morphism $\varrho:M_1\to M_2$ is an isomorphism.
\end{lemma}
\emph{Proof : } Let $a_{0,i}\in \mathbb{Z}_p$ be the exponent of
$M_i$. In Theorem \ref{Robba's algorithm} we have seen that $M_i$
has a basis $\e_i\in M_i$ in which the solution is of the type
$T^{a_{0,i}}\epsilon_i(T)$, where
$\epsilon_i(T)\in\a_H(]1,\infty])$ is a series with coefficients
in $H$. We have then
$\varrho(\bs{\mathrm{e}}_1)=h(T)\bs{\mathrm{e}}_2$, with $h(T)=
T^{a_{0,2}-a_{0,1}}\epsilon_2(T)\epsilon_1(T)^{-1} \in\R_H$. Then
$a_{0,2}-a_{0,1}\in\mathbb{Z}$, and
$\epsilon_2(T)\epsilon_1(T)^{-1}\in \R_H$. Since $\epsilon_i(T)$
is a product of $\bs{\pi}$-exponentials, both $h(T)$ and its
inverse lie in $\R_H$. If $M_1,M_2\in
\MLCS_{\oplus\textrm{-}1}^{\mathrm{reg}}(\mathcal{H}^{\dag}_H)$,
then, by the proof of \ref{algebricity}, the base change necessary
to obtain $\bs{\mathrm{e}}_i$ lies in
$(\mathcal{H}^{\dag}_H)^{\times}$.$\Box$
\begin{remark}
By \ref{Schur}, rank one objects in $\MLCS(\R_H)$ are simple in
$\mathrm{Mod}\textrm{-}\R_H[\d]$. Then, by the Jordan-H\"older
Theorem in $\mathrm{Mod}\textrm{-}\R_H[\d]$, the categories
$\MLCS_{\mathrm{ext}\textrm{-}1}(\R_{H})$, and
$\MLCS_{\mathrm{ext}\textrm{-}1}(\R_{H},S)$ are abelian, and, for
all objects $M$, the set of rank one objects appearing in a
decomposition series of $M$ does not depend, up to the order, on
the chosen decomposition. Moreover the sub-categories
$\MLCS_{\oplus\textrm{-}1}(\R_H)$ and
$\MLCS_{\oplus\textrm{-}1}(\R_H,S)$ are abelian and semi-simple.
The same facts are true for
$\MLCS_{\mathrm{ext}\textrm{-}1}^{\mathrm{reg}}(\mathcal{H}^{\dag}_H)$,
$\MLCS_{\mathrm{ext}\textrm{-}1}^{\mathrm{reg}}(\mathcal{H}^{\dag}_H)$
and
$\MLCS_{\oplus\textrm{-}1}^{\mathrm{reg}}(\mathcal{H}^{\dag}_H,S)$.
\end{remark}
\begin{theorem} \label{M otimes U_m}
Let $N\in\MLCS_{\mathrm{ext}\textrm{-}1}(\R_H,S)$ (resp.
$N\in\MLCS_{\mathrm{ext}\textrm{-}1}^{\mathrm{reg}}(\mathcal{H}^{\dag}_H,S)$).
There exists a basis of $N$ in which the matrix of the derivation
is in Jordan canonical form. In other words, $N$ is a direct sum
of objects of the form $M\otimes U_m$, where $M$ is a rank one
object and $U_m$ is defined by the operator $\d^m$.
\end{theorem}
\emph{Proof : } Let
$M_1,M_2\in\MLCS_{\mathrm{ext}\textrm{-}1}(\R_H)$. By the Robba's
index Theorem for rank one operators whose matrix is a rational
fraction \cite{RoIII}, we have
$\dim\mathrm{Hom}(M_1,M_2)=\dim\mathrm{Ext}^{1}_{\R_K[\d]}(M_1,M_2)$.
This fact does not need the ``spherically complete'' hypothesis on
the field $K$. The Theorem results then by classical
considerations. $\Box$
\begin{remark}
 $\d:\Ed_K\to\Ed_K$ has a big co-kernel, hence
$\mathrm{Ext}_{\Ed_K[\d]}^{1}(\Ed_K,\Ed_K)$ is not one dimensional
(see \cite[pp.133-134]{Astx}). While the theory of rank one
equation over $\R_K$ coincide with the theory over $\Ed_K$, this
is false for rank $\geq 2$.
\end{remark}
\begin{theorem}[Canonical extension]
The canonical restriction functor
$\mathrm{Res}:\MLCS_{\mathrm{ext}\textrm{-}1}^{\mathrm{reg}}(\mathcal{H}^{\dag}_H,S)
\to \MLCS_{\mathrm{ext}\textrm{-}1}(\R_H,S)$ is an equivalence.
\end{theorem}
\emph{Proof : } By \ref{Introductive foundamental Theorem},
$\mathrm{Res}:\MLCS_{\mathrm{ext}\textrm{-}1}^{\mathrm{reg}}(\mathcal{H}^{\dag}_H,S)\to
\MLCS_{\mathrm{ext}\textrm{-}1}(\R_H,S)$ is essentially
surjective. Indeed $\L(a_0,\bs{f}^-(T))$ has its coefficients in
$\mathcal{H}^{\dag}_H$. By \ref{Schur}, two rank one modules in
$\MLCS_{\mathrm{ext}\textrm{-}1}^{\mathrm{reg}}(\mathcal{H}^{\dag}_H,S)$
are isomorphic if and only if they are isomorphic over $\R_H$,
because the base change is given by an over-convergent exponential
in $1+T^{-1}\O_H[[T^{-1}]]$. Hence, by the Schur Lemma
\ref{Schur}, $\mathrm{Res}$ is also fully-faithful. $\Box$
\begin{corollary}
The Tannakian category $\MLCS_{\mathrm{ext}\textrm{-}1}(\R_H,S)$
is neutral.
\end{corollary}
\emph{Proof : } Let
$\omega_S:\MLCS_{\oplus\textrm{-}1}^{\mathrm{reg}}(\mathcal{H}^{\dag}_H,S)
\to \Vect(H)$ be the fiber functor sending a rank one object in
its Taylor solution at $1$ (cf. \eqref{s_x(T)}). An $H$-linear
fiber functor of $\MLCS_{\oplus\textrm{-}1}(\R_H,S)$ is given by
composing $\omega_S$ with a quasi-inverse of $\mathrm{Res}$.$\Box$
\begin{definition}
An affine group scheme $\mathcal{H}$ (over $H$) is \emph{linear}
if there exists a closed immersion
$\mathcal{H}\to\mathrm{GL}_H(V)$, for some finite dimensional
vector space $V$.
\end{definition}
\begin{definition}
Let $\omega_S:\MLCS_{\oplus\textrm{-}1}(\R_H,S)\to\Vect(H)$
\index{omegaS@$\omega_S:\MLCS_{\oplus\textrm{-}1}(\R_H,S)\to\Vect(H)$}
be a fiber functor. We denote
by\index{GL@$\GT_{H}:=\mathrm{Aut}^{\otimes}(\omega_S)$,
$\bs{T}_{H_\infty}$, $\bs{\mathcal{I}}_{H_\infty}$}
$\GT_{H}:=\mathrm{Aut}^{\otimes}(\omega_S)$ the Tannakian group of
$\MLCS_{\oplus\textrm{-}1}(\R_H,S)$.
\end{definition}
\begin{remark} By \ref{M otimes U_m}, the Tannakian group of
$\MLCS_{\mathrm{ext}\textrm{-}1}(\R_H,S)$ is
$\GT_H\times\mathbb{G}_a$.
\end{remark}
\subsubsection{\textbf{Study of $\GT_{H}$}} \label{galois tannakian} For all
(finite dimensional) representations $\rho_V:\GT_H \to
\mathrm{GL}_{H}(V)$, we set $\GT_{H,V} := \rho_V(\GT_H)$. The
group $\GT_{H,V}$ is then linear and affine. Moreover, $\GT_{H,V}$
is diagonalizable (i.e. closed subgroup of the group of diagonal
matrices). The group $\GT_H$ is the inverse limit of its linear
(compact) quotients $\GT_{H,V}$, and is endowed with the limit
topology. Hence $\GT_H$ is abelian, because every $V$ is a direct
sum of rank one objects, and $\GT_{H,V}$ is abelian.
\begin{remark}\label{passage to the limit}
Let $I$ be a non empty directed set. The functor $\varinjlim_{i\in
I}$ is exact if applied to exact sequences of compact algebraic
groups (see  \cite[ch.3,$\S 7$,Cor.1]{Bou-GT}). All exact
sequences in the sequel will be studied at level $\GT_{H,V}$.
\end{remark}
\begin{definition}
We set
$\mathbf{X}(\GT_H):=\Hom_{\mathrm{gr}}^{\mathrm{cont}}(\GT_H,
\mathbb{G}_m\otimes H)$, where ``cont'' means that such a morphism
$\GT_H\to\mathbb{G}_m\otimes H$ factors on a linear quotient
$\GT_{H,V}$.
\end{definition}
Let $\mathrm{Pic}^{\mathrm{sol}}_S(\R_{H_\infty})$ be the
sup-group of $\mathrm{Pic}^{\mathrm{sol}}(\R_{H_\infty})$ formed
by modules whose residue lies in $S$. By Tannakian equivalence, we
have an isomorphism of groups
\begin{equation}
\mathbf{X}(\GT_{H_\infty})\cong
\mathrm{Pic}^{\mathrm{sol}}_{S}(\R_{H_\infty})
\stackrel{\sim}{\longrightarrow} S/\mathbb{Z} \oplus
\PAS(k_{H_\infty}).
\end{equation}
This leads us to recover the group $\GT_{H_\infty}$ itself (cf.
\cite[3.2.6]{Spr}). Let us write
\begin{equation}
S/\mathbb{Z} = S/\mathbb{Z}_{(p)} \oplus
\mathbb{Z}_{(p)}/\mathbb{Z} \;.
\end{equation}

\begin{theorem}
$\GT_{H_\infty}$ is the product of a torus $\bs{T}_{H_\infty}$
(dual of $S/\mathbb{Z}_{(p)}$) with a pro-finite group
$\bs{\mathcal{I}}_{H_\infty}$ (dual of
$\mathbb{Z}_{(p)}/\mathbb{Z}\oplus \PAS(k_{H_\infty})$). This last
is isomorphic to the Galois group
$\mathcal{I}_{E_{H_\infty}}^{\mathrm{ab}}:=\mathrm{Gal}(\E_{H_\infty}^{\mathrm{sep}}/\E_{H_\infty})^{\mathrm{ab}}$,
where $\E_{H_\infty}=k_{H_\infty}((t))$ (cf. \ref{first remark on
PAS})
\begin{equation}
\mathcal{I}_{\E_{H_\infty}}^{\mathrm{ab}}
\stackrel{\sim}{\longrightarrow} \bs{\mathcal{I}}_{H_\infty}\;.
\end{equation}
\end{theorem}
\emph{Proof : } \label{n-th Kummer coverings} The proof is
standard. These two groups have the same character groups. Namely,
by Tannakian equivalence and by Corollary \ref{explicit
description of Pic}, the character group of
$\bs{\mathcal{I}}_{H_\infty}$ is
$\mathbb{Z}_{(p)}/\mathbb{Z}\oplus\PAS(k_{H_\infty})$. By
Artin-Schreier theory, and Kummer theory, this last is also the
character group of $\mathcal{I}_{H_\infty}$. $\Box$
\begin{remark}
We will see in the next section that this isomorphism is induced
by the Fontaine-Katz functor $\bs{\M}^{\dag}$. Actually, this
isomorphism exists even without the hypothesis required in the
definition of this functor.
\end{remark}
%
%
%
%
%
%
%
%
%
%
%
%
%
%
%
%
%
%        ------  T H I R D   A P P L I C A T I O N  ------
%
%
%
%
%
%
%
%
%
%
%
%
%
%
%
\subsection{Differential equations and $\varphi$-modules over
$\Ed_K$ in the abelian case}\label{PHI-MOD}

 In this section the notations, and hypotheses, will follow \cite{Ts}.
 \label{k perfect!!!} We recall that $\mathrm{w}=p$
(cf.\ref{hypothesis uni=p}). We suppose $k$ perfect (used in
\ref{mathcal F dag are still series}). Let $\Lambda/\mathbb{Q}_p$
be a finite extension containing $\mathbb{Q}_p(\xi_s)$. Let
$\mathbb{F}_q$, $q:=p^r$, be the residue field of $\Lambda$.
\begin{hypothesis}
We assume the existence of an \emph{absolute} Frobenius
$\sigma_0:\Lambda\to\Lambda$ (i.e. lifting of the $p$-th power map
$x\mapsto x^p$ of $\mathbb{F}_q$), satisfying
$\sigma_0^{r}=\mathrm{Id}_{\Lambda}$ and $\sigma_0(\pi_s)=\pi_s$.
This is always possible if $\Lambda/\mathbb{Q}_p$ is Galois.
\end{hypothesis}
We let $K:=\Lambda\otimes_{\W(\mathbb{F}_q)}\W(k)$ and
$\sigma:=\mathrm{Id}_{\Lambda}\otimes\F^{r}$. We denote again by
$\sigma_0$ the morphism $(\sigma_0\otimes \F)$ on $K$, then
$\sigma=\sigma_0^r$. We fix a continuous absolute Frobenius
$\varphi_0$ on $\O_{K}^{\dag}$, by setting $\varphi_0(\sum
a_iT^i):= \sum\sigma_0(a_i)\varphi_0(T)^i$, where
$\varphi_0(T)\in\O_{\Lambda}^{\dag}$ is a lifting of $t^p\in
k(\!(t)\!)$ (see definition \ref{what is an absolute Frobenius}).
Then $\varphi_0$ verifies $\varphi_0(\pi_s)=\pi_s$, and
$\varphi_0(\Ed_{\Lambda})\subseteq\Ed_{\Lambda}$. We set
$\varphi=\varphi_0^r$. Both $\varphi$ and $\varphi_0$ extend
uniquely to all unramified extensions of $\Ed_K$, hence they
commute with the action of
$\mathrm{G}_{\E}:=\mathrm{Gal}(\E^{\mathrm{sep}}/\E)$.
\begin{definition}
Let $\Rep_\Lambda^{\mathrm{fin}}(\GE)$
\index{Repfin@$\Rep_\Lambda^{\mathrm{fin}}(\GE)$} be the category
of continuous (finite dimensional) representations $\alpha: \GE
\to \mathrm{GL}_\Lambda(\V)$, such that $\alpha(\IE)$ is finite.
\end{definition}

\begin{definition} \label{V_chi}\index{Valpha@$\V_\alpha=$repr. given by the char. $\alpha$}
Let $\alpha:\GE\to\Lambda^{\times}$ be a character such that
$\alpha(\IE)$ is finite. Then we denote by $\V_\alpha\in
\Rep_\Lambda^{\mathrm{fin}}(\GE)$ the rank one representation of
$\GE$ given by
$$\gamma(\bs{\mathrm{e}}):=
\alpha(\gamma)\cdot\bs{\mathrm{e}}\;,\quad\textrm{for all }
\gamma\in\GE\;,$$ where $\bs{\mathrm{e}}\in\V_\alpha$ is a basis.
We denote by $\bs{\mathrm{D}}^{\dag}(\V_\alpha)$ (resp.
$\bs{\M}^\dag(\V_\alpha)$) the $\varphi-\nabla$-module over
$\Ed_{K}$ (resp. $\nabla$-module over $\R_K$) attached to
$\V_\alpha$.
Namely\index{D(V)@$\bs{\mathrm{D}}^{\dag}(\V_{\alpha})$,
$\bs{\M}^{\dag}(\V_{\alpha})$}
\begin{equation}\label{explicit expression of Mdag}
\bs{\mathrm{D}}^\dag(\V_\alpha)=
(\V_\alpha\otimes_{\Lambda}\mathcal{E}_K^{\dag,\mathrm{unr}})^{\GE}
\quad,\quad \bs{\M}^\dag(\V_\alpha)=
\bs{\mathrm{D}}^\dag(\V_\alpha)\otimes_{\Ed_K}\R_{K}\;.
\end{equation}
We recall that $\bs{\M}^\dag(\V_\alpha)$ is endowed with the
unique derivation ``commuting" with $\varphi$ (cf.
\cite[2.2.4]{Fo}). This derivation is $\nabla=1\otimes\d$. By
\ref{frobenius implies solvability},
$\bs{\M}^\dag(\V_\alpha)\in\MLCS(\R_{K})$.

\end{definition}
\begin{definition}\label{V_chi imbrognlio}
We will identify $\mathbb{Z}/p^{s+1}\mathbb{Z}$ with
$\bs{\mu}_{p^{s+1}}$, by sending
$$1\mapsto \xi_s\;:\;\mathbb{Z}/p^{s+1}\mathbb{Z}
\stackrel{\sim}{\longrightarrow}\bs{\mu}_{p^{s+1}},$$ where
$\xi_s$ is the unique $p^{s+1}$-th root of $1$ verifying (cf.
\ref{relations on equivalence classe})
\begin{equation}
|(\xi_s-1)-\pi_s|<|\pi_s|\;.
\end{equation}
If $\alpha \in
\Hom^{\mathrm{cont}}(\GE,\mathbb{Z}/p^{s+1}\mathbb{Z})$, we again
denote by $\V_\alpha$ the representation given by
$$\gamma(\bs{\mathrm{e}}):=\xi_s^{\alpha(\gamma)}
\cdot \bs{\mathrm{e}}\;, \quad\textrm{for all }\gamma\in\GE.$$
\end{definition}
\begin{remark}
This definition is chosen ``ad hoc'' to be the inverse of the
action of $\GE$ described in \eqref{action of G on xi...}.
\end{remark}

\begin{remark} \label{diff is indep from GE/IE}
Let $\alpha:\GE\to\Lambda^{\times}$ be a continuous character,
then $\alpha$ factors on the abelianized $\GE^{\mathrm{ab}}$. Let
$\IE^{\mathrm{ab}}$ be the inertia of $\GE^{\mathrm{ab}}$, and
$\mathrm{G}_{k}^{\mathrm{ab}}$ be the abelianized of
$\mathrm{Gal}(k^{\mathrm{sep}}/k)$. Since $k$ is perfect, then, by
\ref{mathcal F dag are still series}, the exact sequence $1\to
\IE^{\mathrm{ab}}\to \GE^{\mathrm{ab}}\to
\mathrm{G}_{k}^{\mathrm{ab}}\to 1$ is split, hence
\begin{equation}
\GE^{\mathrm{ab}}=\IE^{\mathrm{ab}}\oplus\G_{k}^{\mathrm{ab}}\;,
\end{equation}
and $\alpha=\alpha^{-} \cdot \alpha_{0}$, where
$\alpha^{-}:\IE^{\mathrm{ab}}\to\Lambda^{\times}$ and
$\alpha_{0}:\mathrm{G}_{k}^{\mathrm{ab}} \to \Lambda^{\times}$.
Then
$$\V_\alpha=\V_{\alpha^{-}}\otimes\V_{\alpha_{0}}.$$
 We observe that
$\bs{\M}^{\dag}(\V_{\alpha_{0}})\stackrel{\sim}{\to}\R_{K}$,  is
trivial because its solution is a constant. Indeed, the extension
of $\O^\dag_{K}$ defined by $\alpha_{0}$ is
$\O^\dag_{K}\otimes_{K} H$, for some unramified extensions $H/K$.
In the sequel we will treat only characters
$\alpha:\GE\to\Lambda^{\times}$ with \emph{finite image}, this
will be restrictive in terms of $\varphi$-modules but not in terms
of differential modules (cf. \ref{diff is indep from GE/IE}).
Indeed,
\begin{equation}
\bs{\mathrm{D}}^\dag(\V_\alpha)=\bs{\mathrm{D}}^\dag(\V_{\alpha^{-}})\otimes
\bs{\mathrm{D}}^\dag(\V_{\alpha_{0}})\qquad;\qquad
\bs{\M}^\dag(\V_\alpha)=\bs{\M}^\dag(\V_{\alpha^{-}}).
\end{equation}
\end{remark}
\begin{remark}
Points $(4)$ and $(5)$ of the following theorem have been already
proved in \cite{Ma} in the case $p\neq 2$ and rank one, and in
\cite{Crew-can-ext}, \cite{Matsuda-unipotent}, \cite{Tsu-swan} in
the general case. Moreover we thank the referee to pointed out to
us that the explicit form of the differential operator (answer to
(5) of \ref{presentation of the problems}) was written in the
proof of Lemma $5.2$ of \cite{Ma}, in the case $p\neq 2$.
\end{remark}
\begin{theorem} \label{second main identification}
Let $\overline{\bs{f}}(t)\in\W_s(\E)$ and let
$\alpha=\delta(\overline{\bs{f}}(t))$ be the Artin-Schreier
character defined by $\overline{\bs{f}}(t)$ (cf.
\eqref{artin-screier-diagram}). Let $(\Ed_K)'$ be the unramified
extension of $\Ed_K$ corresponding, by henselianity, to the
separable extension of $k(\!(t)\!)$ defined by $\alpha$. Then
\begin{enumerate}
\item a basis of $\bs{\mathrm{D}}^{\dag}(\V_\alpha)$ is given by
\begin{equation}
y:=\bs{\mathrm{e}}\otimes \theta_{p^s}(\bs{\nu},1)\;,
\end{equation}
where $\bs{\mathrm{e}}\in\V_\alpha$ is the basis of \ref{V_chi
imbrognlio} and
$\bs{\nu}\in\W_s(\widehat{\mathcal{E}}_K^{\mathrm{unr}})$ is a
solution of
\begin{equation}
\varphi_0(\bs{\nu})-\bs{\nu}=\bs{f}(T)\;,
\end{equation}
where $\bs{f}(T)\in\W_s(\O_K[[T]][T^{-1}])$ is an arbitrary
lifting of $\overline{\bs{f}}(t)$;

\item the Frobenius $\varphi_0$ acts on $\V_{\alpha}$, moreover
$\varphi_0(y)= \theta_{p^s}(\bs{f}(T),1)\cdot y$.
%\begin{equation}
%\varphi_0\bigl(\bs{\mathrm{e}}\otimes
%\theta_{p^s}(\bs{\nu},1)\bigr)= \theta_{p^s}(\bs{f}(T),1)\cdot
%\bigl(\bs{\mathrm{e}}\otimes\theta_{p^s}(\bs{\nu},1)\bigr)\;.
%\end{equation}
Hence, if $\mathrm{Tr}(\bs{f}(T)):=
\bs{f}(T)+\varphi_0(\bs{f}(T))+\cdots
+\varphi_0^{r-1}(\bs{f}(T))$, we have
%\begin{equation}
%\varphi\bigl(\bs{\mathrm{e}}\otimes \theta_{p^s}(\bs{\nu},1)\bigr)
%= \theta_{p^s}(\mathrm{Tr}(\bs{f}(T)),1)\cdot
%\bigl(\bs{\mathrm{e}}\otimes\theta_{p^s}(\bs{\nu},1)\bigr)\;;
%\end{equation}
%
\begin{equation}
\varphi(y) = \theta_{p^s}(\mathrm{Tr}(\bs{f}(T)),1)\cdot y;
\end{equation}

\item By \ref{lift-diagr-for Ed}, one has
$(\Ed_K)'(\pi_s)=\Ed_{K_s}[\theta_{p^s}(\bs{\nu}^-,1)]$ (cf.
definition \ref{dec -0+}). This extension can be identified with
the extension
\begin{equation}
\Ed_{K_s}[\theta_{p^s}(\bs{\nu}^-,1)]\xrightarrow[]{\;\;\sim\;\;}\Ed_{K_s}[\mathrm{e}_{p^s}(\bs{f}^-(T),1)]
\end{equation}
by sending $\theta_{p^s}(\bs{\nu},1)$ into
$\mathrm{e}_{p^s}(\bs{f}^-(T),1)$. In particular, if $\pi_s\in K$,
one has
\begin{equation}
\widetilde{y}=\e\otimes \mathrm{e}_{p^s}(\bs{f}^-(T),1)\;.
\end{equation}
Moreover
$\varphi_0(\widetilde{y})=\theta_{p^s}(\bs{f}^-(T),1)\cdot
\widetilde{y}$, and
$\varphi(\widetilde{y})=\theta_{p^s}(\mathrm{Tr}(\bs{f}^-(T)),1)\cdot
\widetilde{y}$.

\item The isomorphism class of $\bs{\M}^{\dag}(\V_\alpha)$ depends
only on $\alpha^-$ and
\begin{equation}
\bs{\M}^{\dag}(\V_\alpha) \stackrel{\sim}{\longrightarrow}
\M(0,\alpha^-)\;;
\end{equation}

\item the irregularity of $\bs{\M}^{\dag}(\V_\alpha)$ is equal to
the Swan conductor of $\V_\alpha$.
\end{enumerate}
\end{theorem}
\emph{Proof : } Let $\E=k(\!(t)\!)$. For all $\gamma\in\G_{\E}=
\mathrm{Gal}(\E^{\mathrm{sep}}/\E)$, we have (cf. \eqref{action of
gamma on theta})
\begin{equation}
\gamma(\bs{\mathrm{e}} \otimes
\theta_{p^s}(\bs{\nu},1))=(\xi_s^{\alpha(\gamma)}\cdot\bs{\mathrm{e}})
\otimes (\xi_s^{-\alpha(\gamma)}\cdot
\theta_{p^s}(\bs{\nu},1))=\bs{\mathrm{e}} \otimes
\theta_{p^s}(\bs{\nu},1),
\end{equation}
hence $\bs{\mathrm{e}} \otimes
\theta_{p^s}(\bs{\nu},1)\in\bs{\mathrm{D}}^\dag(\V_\alpha)$.
Moreover, $\varphi_0(\bs{\mathrm{e}} \otimes
\theta_{p^s}(\bs{\nu},1))=
\bs{\mathrm{e}}\otimes\varphi_0(\theta_{p^s}(\bs{\nu},1))$ and
\begin{equation}
\varphi_0(\theta_{p^s}(\bs{\nu},1)) \stackrel{(*)}{=}
\theta_{p^s}(\varphi_0(\bs{\nu}),1)=
\theta_{p^s}(\varphi_0(\bs{\nu})-\bs{\nu},1) \cdot
\theta_{p^s}(\bs{\nu},1)\;.
\end{equation}
The equality $(*)$ is true because
$\varphi_0(\pi_s)=\sigma_0(\pi_s)=\pi_s$, for all $s\geq 0$. The
derivation on $\bs{\M}^{\dag}(\V_\alpha)$ arises from the
derivation on $\mathcal{E}_K^{\dag,\mathrm{unr}}$ (cf.
\eqref{explicit expression of Mdag}), hence the operator attached
to the basis $\bs{\mathrm{e}}\otimes
\theta_{p^s}(\bs{\nu},1)\in\bs{\M}^{\dag}_K(\V_\alpha)$ is (cf.
\ref{matrix of derivation})
$\d-\partial_{T,\log}(\theta_{p^s}(\bs{\nu},1))$. As explained in
\ref{diff is indep from GE/IE}, the isomorphism class of
$\bs{\M}^{\dag}(\V_\alpha)$ depends only on
$\alpha^-=\delta(\overline{\bs{f}^-}(t))$. Hence, we can suppose
$\alpha=\alpha^-$ and
$\bs{f}(T)=\bs{f}^-(T)\in\W_s(T^{-1}\O_K[T^{-1}])$ (cf. \ref{dec
-0+}). Let us write $\bs{\nu}^-$ instead of $\bs{\nu}$. By
equation \eqref{theta^p^m+1=et^p^m+1}, we have
\begin{equation}
\theta_{p^s}(\bs{\nu}^-,1)^{p^{s+1}}=
\mathrm{e}_{p^s}(\bs{f}^-(T),1)^{p^{s+1}},
\end{equation}
hence $\partial_{T,\log}(\theta_{p^s}(\bs{\nu}^-,1))=
\partial_{T,\log}(\mathrm{e}_{p^s}(\bs{f}^-(T),1))$. This establishes point
$(3)$ and $(4)$.

Both the Swan conductor and the irregularity are stable under
extension of the constant field $K$, hence we can suppose
$K=K^{\mathrm{alg}}$. We can suppose that $\bs{f}^-(T)$ is pure
(cf. \ref{pure witt vector}), because both the Swan conductor and
the irregularity depend only on $\alpha^-$. Write
$\bs{f}^-(T)=\sum_{n\in\J}\lb_{np^{m(n)}}T^{-np^{m(n)}}$. Since
the irregularities of the $\lb_{np^{m(n)}}T^{-np^{m(n)}}$'s are
all different we can suppose
$\bs{f}^-(T)=\lb_{np^{m(n)}}T^{-np^{m(n)}}$. Now write explicitly
(cf. definition \ref{co-monomials and PAS})
\begin{eqnarray*}
\lb_{np^{m(n)}}T^{-np^{m(n)}}&=&
(\lambda_{r}T^{-np^r},\lambda_{r+1}T^{-np^{r+1}},\ldots,
\lambda_{m}T^{-np^{m(n)}})\\
&=&(\lambda_{r}T^{-np^r},0\ldots,0)+\cdots+
(0,\ldots,0,\lambda_{m}T^{-np^{m(n)}})\;.
\end{eqnarray*}
Since, by Reduction Theorem \ref{reduction to k}, the
irregularities of the these vectors are all different, and since
both the irregularity and the Swan conductor are invariant by $\V$
(cf. \eqref{shiffttt}), we can suppose $\bs{f}^-(T)=(\lambda
T^{-n},0,\ldots,0)$, with $|\lambda|=1$. Moreover, since
$K=K^{\mathrm{alg}}$, the residue field is perfect and, replacing
$\lambda T^{-n}$ with $\lambda^{1/p^k}T^{-n/p^k}$, we can suppose
$(n,p)=1$ (cf. \ref{VF is the same of F}). The irregularity is
then $np^s$ (cf. \ref{irreg of (0...0f0...0) with (n,p)=1}), and
it is equal to the Swan conductor (see for example
\cite{Brylinski}). This Theorem is the analogue of \ref{garuti}
for Artin-Schreier characters of $\mathrm{G}_{\E}$). $\Box$
\begin{remark}
Suppose that the character is totally ramified, and choose
$\bs{f}^-(T)$ in $\W_s(T^{-1}\O_K[T^{-1}])$, then
$\theta_{p^s}(\bs{\nu},1)=\mathrm{e}_{p^s}(\bs{f}^-(T),1)$.
\end{remark}
\begin{remark}
Let $a_0=\frac{m}{n}\in\mathbb{Z}_{p}\cap\mathbb{Q}$. Suppose that
$\bs{\mu}_n\subset k$. Let $\beta_{a_0}:
\mathrm{G}\to\bs{\mu}_{n}\subset \Lambda^{\times}$ be the Kummer
character defined by $t^{a_0}$. We have $\beta_{a_0}(\gamma)=
\gamma(t^{a_0})/t^{a_0}$. As before, a basis of
$\bs{\mathrm{D}}^{\dag}(\V_{\beta_{a_0}})$ is given by $\e\otimes
T^{-a_0}\in\V_{\beta_{a_0}}\otimes\mathcal{E}^{\dag,\mathrm{unr}}_{K}$,
because $\gamma(\e):=\beta_{a_0}(\gamma)\e$. Then $\varphi(\e
\otimes T^{-a_0}) = T^{a_0}\varphi(T^{-a_0})\cdot (\e \otimes
T^{-a_0})$, and $\bs{\M}^{\dag}(\V_{\beta_0})=\M(a_0,0)$. We do
not necessarily have an action of $\varphi_0$, because $\sigma_0$
does not fix the $n$-th root of $1$.
\end{remark}

\printindex

\bibliographystyle{alpha}
\bibliography{Rk1-Arxiv}

\end{document}